%% file: draft.tex
\documentclass[lettersize,journal]{IEEEtran}
\IEEEoverridecommandlockouts
\usepackage{cite}
\usepackage{amsmath, amssymb, amsfonts}
\usepackage{algorithm}
\usepackage{graphicx}
\usepackage[outdir=./]{epstopdf}
\usepackage{textcomp}
\usepackage{xcolor}
\usepackage{CJKutf8}
\usepackage{multirow} 
\usepackage{makecell}
\usepackage{hyperref}
\usepackage{algorithmic}

\def\BibTeX{{\rm B\kern-.05em{\sc i\kern-.025em b}\kern-.08em
T\kern-.1667em\lower.7ex\hbox{E}\kern-.125emX}}
\usepackage{subfigure}
\usepackage{stfloats}
\usepackage{subcaption}

\newtheorem{Thm}{Theorem}
\newtheorem{Lem}{Lemma}

\newtheorem{Prob}{Problem}

\newcommand{\red}{\textcolor{black}}
\newcommand{\blue}{\textcolor{black}}

\begin{document}
\title{Fast MLE and MAPE-Based Device Activity Detection for Grant-Free Access via PSCA and PSCA-Net}

\author{Bowen~Tan
        and~Ying~Cui
\thanks{Manuscript received 5 September 2024; revised 25 December 2024; accepted 4 February 2025. Date of publication xxx; date of current version xxx. This work was supported in part by the National Natural Science Foundation of China (62371412), the National Key Research and Development Program of China (2024YFE0200603), the Guangdong Basic and Applied Basic Research Natural Science Funding Scheme (2024A1515011184), the Guangzhou-HKUST (GZ) Joint Funding Scheme (2024A03J0539), the Guangdong Provincial Key Lab of Integrated Communication, Sensing, and Computation for Ubiquitous Internet of Things (2023B1212010007), and Guangzhou Municipal Science and Technology Project (2023A03J0011, 024A03J0623). This paper was presented in part at IEEE GLOBECOM 2024 \cite{tan2024mle}. The associate editor coordinating the review of this article and approving it for publication was D. Niyato. (Corresponding author: Ying Cui.)}
\thanks{Bowen Tan and Ying Cui are with the IoT Thrust, The Hong Kong University of Science and
Technology (Guangzhou), Guangzhou 511400, China. Ying Cui is also with the Department of Electronic Engineering, The Hong Kong University of Science and Technology, Hong Kong, China (e-mail: \href{mailto:yingcui@ust.hk}{yingcui@ust.hk}).}
\thanks{\red{Digital Object Identifier 10.1109/TWC.2025.3545649}}
}
\maketitle     
\input{abstract}

\input{introduction.tex}
\input{system.tex}
\input{known.tex}
\input{unknown.tex}
\input{numerical.tex}

\input{appendix.tex}
\bibliographystyle{IEEEtran}
\bibliography{ref}
\begin{IEEEbiography}[{\includegraphics[height=1.25in]{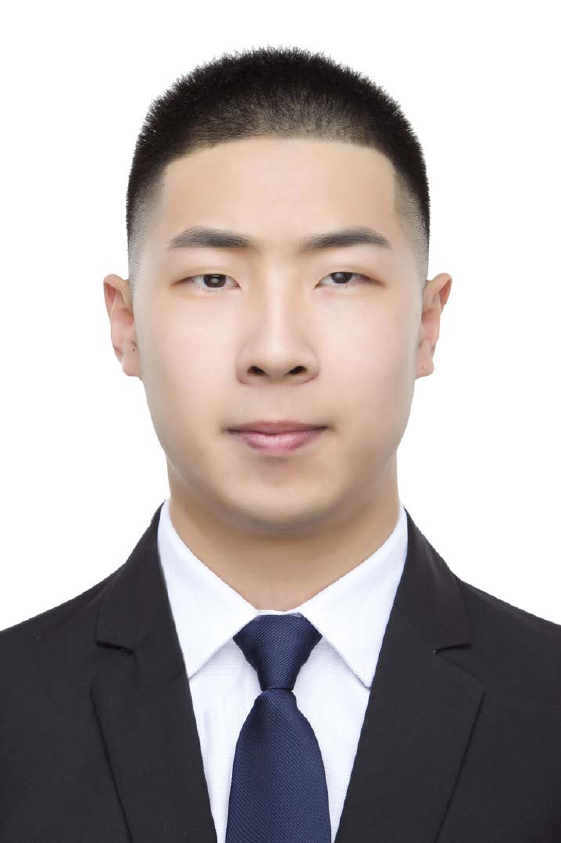}}]
{Bowen Tan} received his B.S. degree in Statistics from Nankai University, China, in 2022 and his M.phil. degree from the Hong Kong University of Science and Technology (Guangzhou), China, in 2024. He is currently an algorithm engineer in Cardinal Operations. His research interests include optimization, learning, wireless communication, power systems, and bioinformatics.
\end{IEEEbiography}
\begin{IEEEbiography}[{\includegraphics[height=1.25in]{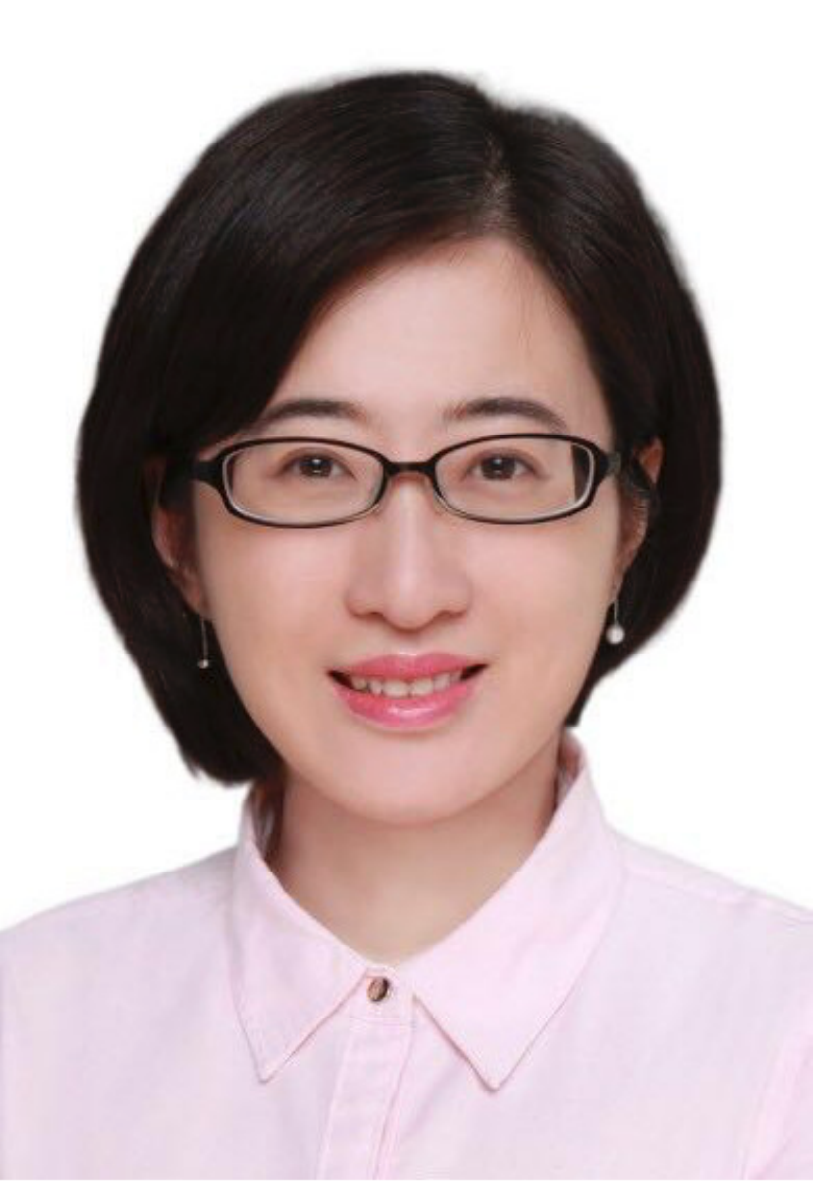}}]
{Ying Cui} received her B.Eng. degree in Electronic and Information Engineering from Xi'an Jiao Tong University, China, in 2007 and her Ph.D. degree from the Hong Kong University of Science and Technology, Hong Kong, in 2012. She held visiting positions at Yale University, US, in 2011 and Macquarie University, Australia, in 2012. From June 2012 to June 2013, she was a postdoctoral research associate at Northeastern University, US. From July 2013 to December 2014, she was a postdoctoral research associate at the Massachusetts Institute of Technology, US. From January 2015 to July 2022, she was an Associate Professor at Shanghai Jiao Tong University, China. Since August 2022, she has been an Associate Professor with the IoT Thrust at The Hong Kong University of Science and Technology (Guangzhou) and an Affiliate Associate Professor with the Department of ECE at The Hong Kong University of Science and Technology. Her current research interests include optimization, learning, IoT communications, mobile edge caching and computing, and multimedia transmission. She was selected to the National Young Talent Program in 2014. She received Best Paper Awards from IEEE ICC 2015 and IEEE GLOBECOM 2021. She serves as an Editor for IEEE Transactions on Communications and IEEE Transactions on Machine Learning in Communications and Networking. She served as an Editor for the IEEE Transactions on Wireless Communications from 2018 to 2024.
\end{IEEEbiography}
\end{document}

%% file: abstract.tex
\begin{abstract}
	Fast and accurate device activity detection is the critical challenge in grant-free access for supporting massive machine-type communications (mMTC) and ultra-reliable low-latency communications (URLLC) in 5G and beyond. The state-of-the-art methods have unsatisfactory error rates or computation times. To address these outstanding issues, we propose new maximum likelihood estimation (MLE) and maximum a posterior estimation (MAPE) based device activity detection methods for known and unknown pathloss that achieve superior error rate and computation time tradeoffs using optimization and deep learning techniques. Specifically, we investigate four non-convex optimization problems for MLE and MAPE in the two pathloss cases, with one MAPE problem being formulated for the first time. For each non-convex problem, we develop an innovative parallel iterative algorithm using the \red{parallel successive convex approximation (PSCA) method}. Each PSCA-based algorithm allows parallel computations, uses up to the objective function’s second-order information, converges to the problem’s stationary points, and has a low per-iteration computational complexity compared to the state-of-the-art algorithms. Then, for each PSCA-based iterative algorithm, we present a deep unrolling neural network implementation, called PSCA-Net, to further reduce the computation time. Each PSCA-Net elegantly marries the underlying PSCA-based algorithm’s parallel computation mechanism with the parallelizable neural network architecture and effectively optimizes its step sizes based on vast data samples to speed up the convergence. Numerical results demonstrate that the proposed methods can significantly reduce the error rate and computation time compared to the state-of-the-art methods, revealing their significant values for grant-free access.
 \end{abstract}
\begin{IEEEkeywords}
	Massive machine-type communications (mMTC), ultra-reliable low-latency communications (URLLC), grant-free access, maximum likelihood estimation (MLE), maximum a posterior estimation (MAPE), parallel successive convex approximation (PSCA), deep unrolling neural network.
\end{IEEEkeywords}

%% file: introduction.tex
\section{Introduction}
\label{Introduction}
Massive machine-type communications (mMTC) and ultra-reliable low-latency communications (URLLC) have become the critical use cases of 5G and beyond 5G (B5G) wireless networks to support the Internet of Things (IoT) in many fields \cite{huawei2016discussion}. The main challenge in realizing mMTC is to effectively support IoT devices that are power-constrained and intermittently active with small packets to send. In contrast, the primary challenge of the URLLC scenario is to provide ultra-high reliability of 99.999\% within a very low latency in an order of 1 ms. Grant-free (non-orthogonal multiple) access, which enables devices to transmit data in an arrive-and-go manner without waiting for the base station (BS)'s grant, has been identified as a promising solution to address the challenges in mMTC and URLLC by the Third Generation Partnership Project (3GPP) because it can reduce transmission latency, lower signaling overhead, and improve device energy efficiency. Under grant-free access, devices with pre-assigned non-orthogonal pilots send their pilots and data once activated, and the multi-antenna BS estimates device activities and channel conditions from the received pilot signals and detects active users' transmitted data from the received data signals.
\subsection{Related Works} 

Early studies on grant-free access focus on a simple scenario, i.e., single-cell narrow-band \red{networks} under flat Rayleigh fading and perfect time and frequency synchronization \cite{chen2018sparse,liu2018massive,senel2018grantfree,shao2019unified,fengler2021nonbayesian,chen2019covariance,wang2021efficient,dong2022faster,lin2023sparsity,jiang2020mapbased,jiang2022ml,shao2020covariancebased,wang2021accelerating,jiang2023statistical}. Later studies extend to more complex scenarios, e.g., multi-cell networks \cite{chen2019multicell,jiang2022ml,shao2020covariancebased,wang2021accelerating}, flat Rician fading \cite{liu2024mle}, frequency selective Rayleigh fading \cite{jiang2023statistical}, imperfect time or frequency synchronization \cite{liu2024mle}, etc. There are generally two classes of statistical methods for device activity detection and channel estimation in the pilot transmission phase. One is to estimate device activities directly \cite{fengler2021nonbayesian,chen2019covariance,wang2021efficient,dong2022faster,lin2023sparsity,jiang2020mapbased,jiang2022ml,shao2020covariancebased,wang2021accelerating,jiang2023statistical,liu2024mle}. Given very accurate device activity detection results, the channel estimations of active devices are readily obtained via channel estimation methods such as \red{least squares (LS) estimation and minimum mean square error (MMSE) estimation} \cite{cui2021jointly,jiang2022ml}. The other is to first estimate effective channels\footnote{The effective channel of each device is the channel if the device is active and zero otherwise.}  and then obtain device activity \red{detections} and the channel estimations of active devices via a thresholding approach \cite{chen2018sparse,liu2018massive,senel2018grantfree,shao2019unified,yuan2020iterative,chen2019multicell, hara2022blinda}.
Notably, existing work investigates known \cite{dong2022faster,lin2023sparsity,jiang2020mapbased,jiang2022ml,shao2020covariancebased,jiang2023statistical,chen2018sparse,liu2018massive,senel2018grantfree,yuan2020iterative,chen2019multicell} and unknown \cite{fengler2021nonbayesian,wang2021efficient,wang2021accelerating} pathloss cases. \red{In the known pathloss case, each device's activity is estimated. The known pathloss case is suitable for static or slowly moving IoT devices. On the contrary, in the unknown pathloss case, each device's effective pathloss, which is the pathloss if the device is active and zero otherwise, is first estimated, and the device's activity is deduced based on a thresholding method. The unknown pathloss case generally applies to slowly or fast moving IoT devices.}

This paper mainly concentrates on
the simple scenario and systematically studies device activity
detection in the known and unknown pathloss cases, which is
significant for guaranteeing excellent performance for mMTC
and URLLC under grant-free access.\footnote{\red{We start with the simple scenario to gain first-order design insights and highlight the core concepts. The proposed methods for the simple scenario can be readily extended to the complex scenarios.}} The following introduction focuses on statistical estimation methods, including classical statistical estimation methods (which assume parameters of interest to be deterministic but unknown), e.g., maximum likelihood estimation (MLE), and Bayesian estimation methods (which view the parameters of interest as realizations of random variables with known prior distributions), e.g., maximum a posterior estimation (MAPE) and MMSE.

First, MLE \cite{fengler2021nonbayesian,chen2019covariance,wang2021efficient,dong2022faster,lin2023sparsity,jiang2020mapbased,jiang2022ml,shao2020covariancebased,wang2021accelerating,jiang2023statistical,liu2024mle} and MAPE \cite{jiang2022ml,jiang2020mapbased,jiang2023statistical} methods are developed to estimate device activity states in the known \cite{fengler2021nonbayesian,chen2019covariance,wang2021efficient,dong2022faster,jiang2020mapbased,jiang2022ml,shao2020covariancebased,wang2021accelerating,jiang2023statistical} and unknown \cite{fengler2021nonbayesian,wang2021efficient,lin2023sparsity} pathloss cases. MLE-based device activity detection methods, proposed for the known \cite{fengler2021nonbayesian,chen2019covariance,wang2021efficient,dong2022faster,jiang2020mapbased,jiang2022ml,shao2020covariancebased,wang2021accelerating,jiang2023statistical} and unknown \cite{fengler2021nonbayesian,lin2023sparsity,wang2021efficient} pathloss cases typically have closed-form per-iteration updates and can achieve high detection accuracy with a long computation time. Specifically, the block coordinate descent (BCD) method is adopted to obtain stationary points of the MLE problems for device activity detection in the known \cite{fengler2021nonbayesian,chen2019covariance,wang2021efficient,dong2022faster,lin2023sparsity,jiang2020mapbased,jiang2022ml,shao2020covariancebased,wang2021accelerating,jiang2023statistical} and unknown \cite{fengler2021nonbayesian,lin2023sparsity} pathloss cases. The resulting algorithms for the known and unknown cases are called BCD-ML-K and BCD-ML-UD, respectively. Moreover, the projection gradient (PG) method is used to obtain stationary points of the MLE problem in the unknown pathloss case \cite{wang2021efficient}, leading to an algorithm termed PG-ML-UD. PG-ML-UD can be readily extended to the known pathloss case, and the extension is called PG-ML-K. \red{MAPE-based device activity detection methods, }proposed only for the known pathloss case, can further improve the detection accuracy with increased computation times compared to the MLE-based counterparts by utilizing the prior distributions of device activities, including independent \cite{jiang2020mapbased,jiang2022ml,jiang2023statistical} and correlated \cite{jiang2022ml,jiang2023statistical} device activities.
Similarly, BCD is utilized to obtain stationary points of the MAPE problem in the known pathloss case \cite{jiang2022ml,jiang2020mapbased,jiang2023statistical}.
The resulting algorithms are called BCD-MAP-K. Noteworthily, MAPE-based device activity detection methods for unknown pathloss cases have never been studied. 
The BCD-based algorithms analytically solve the coordinate descent optimization problems sequentially in each iteration. Although the BCD-based algorithms approximately utilize up to the second-order information in the coordinate updates \red{(which will be shown in Section \ref{Known})}, their sequential updates yield long computation times. The PG-based \red{algorithms} update the estimates of all device activities in parallel in each iteration, reducing the computation time. However, as first-order algorithms, PG-based algorithms still have unsatisfactory convergence speeds and computation times.

Next, MMSE is adopted to estimate the effective channel of each device in the known \cite{chen2018sparse,liu2018massive,senel2018grantfree,shao2019unified,yuan2020iterative,chen2019multicell} and unknown \cite{hara2022blinda} pathloss cases.\footnote{\red{The MMSE problems are less tractable than the MLE and MAPE problems due} to the integrals in the objective functions. No iterative algorithms can converge to the stationary points of the nonconvex MMSE problems.} In the known pathloss case, the approximate message passing (AMP) method is adopted to approximately solve the MMSE problems via \red{low-complexity parallel updates} in each iteration \cite{chen2018sparse,liu2018massive,senel2018grantfree,shao2019unified,yuan2020iterative,chen2019multicell}. The resulting algorithms are referred to as AMP-K. In the unknown pathloss case, \cite{hara2022blinda} first solves an MLE problem for the pathloss via the expectation maximization (EM) algorithm \cite{hara2022blinda} \red{and} then approximately solves the MMSE problem under the estimated pathloss using the AMP method. The resulting algorithm is called EM-AMP-UD.  Compared to MLE-based and MAPE-based methods, MMSE-based methods can achieve short computation times at the sacrifice of detection accuracy due to the inherent approximation mechanism for enforcing computational complexity reduction in solving \red{the MMSE optimization problems}.

Last, data-driven and model-driven deep learning \red{methods} have recently been proposed for device activity detection \red{mainly} in the known pathloss case to further improve detection accuracy and reduce computation time \cite{cui2024deep}. Specifically, data-driven deep learning approaches, which purely rely on basic neural network modules, including fully connected neural networks \cite{li2020jointly,li2019joint,zhang2020jointly, zhu2021deep}, convolutional neural networks (CNNs) \cite{wu2021cnn}, recurrent neural networks (RNNs) {\cite{ahn2020deep,shi2022algorithm}, \red{residual} connected modules \cite{kim2020deep}, etc., can achieve short computation \red{times} and high accuracy for a small or moderate number of devices, but fail for massive devices. Model-driven deep learning \red{methods} \cite{cui2021jointly,zhu2021deep,shi2022algorithm}, obtained by unrolling Bayesian statistical estimation-based methods, such as MAPE-based methods via BCD-MAP-K and MMSE-based methods via AMP-MAP-K, successfully cope with the scenario of massive devices. However, the resulting deep unrolling neural network methods are still underperforming due to the inherent defects of BCD-MAP-K and AMP-K \red{in} computation time and detection accuracy, respectively. 

\subsection{Contributions}
In summary, existing methods for known and unknown pathloss cases cannot achieve satisfactory detection accuracy and computation time tradeoffs. In addition, MAPE-based device activity detection in the unknown pathloss case remains open. In this paper, we propose new MLE and MAPE-based device activity detection methods for the known and unknown pathloss cases that can achieve superior error rate and computation time tradeoffs using optimization and deep learning techniques. Specifically, in the known pathloss case, we propose new methods for solving the existing non-convex MLE \cite{fengler2021nonbayesian,jiang2022ml} and MAPE \cite{jiang2022ml,jiang2023statistical} problems where device activities are treated as deterministic unknown constants and realizations of random variables, respectively. In the unknown pathloss case, we propose new methods for solving the existing non-convex MLE \red{problem} \cite{fengler2021nonbayesian,wang2021efficient} where effective pathloss values are treated as deterministic unknown constants, and we formulate a novel MAPE problem by viewing effective pathloss values as realizations of random variables with location information and propose methods for solving it. The main contributions of this paper are \red{summarized} as follows. 

(i) We assume that device activities follow independent Bernoulli distributions and device locations follow independent uniform distributions. By elegant approximation, we build a tractable model for \red{random} effective pathloss. Based on this, we formulate the MAPE problem of effective pathloss as a non-convex problem, which captures the prior information on device activities and locations.

(ii) For each of the four non-convex estimation problems, i.e., the MLE and MAPE problems for device activities and effective pathloss in the two pathloss cases, we propose an iterative algorithm to obtain stationary points using the parallel successive convex approximation (PSCA) method. The resulting algorithms for MLE in the known and unknown pathloss cases are called PSCA-ML-K and PSCA-ML-UD, respectively, and the resulting algorithms for MAPE in the known and unknown pathloss cases are termed PSCA-MAP-K and PSCA-MAP-UR, respectively. We show that PSCA-ML-K, PSCA-ML-UD, and PSCA-MAP-K can be viewed as the parallel counterparts of BCD-ML-K \cite{fengler2021nonbayesian,jiang2022ml}, BCD-ML-UD \cite{fengler2021nonbayesian}, and BCD-MAP-K \cite{jiang2022ml,jiang2023statistical}\red{, respectively}. Besides, we show that PSCA-ML-K and PSCA-ML-UD use higher-order information of the objective functions than PG-ML-K and PG-ML-UD \cite{wang2021efficient}, respectively. Last, we show that the computational complexity of each PSCA-based algorithm is lower than the BCD and PG-based counterparts.

(iii) For each of the four PSCA-based iterative algorithms, we present a deep unrolling neural network implementation to further reduce the computation time. The resulting neural networks are called PSCA-ML-K-Net, PSCA-MAP-K-Net, PSCA-ML-UD-Net, and PSCA-MAP-UR-Net. Specifically, each block of a PSCA-Net implements one iteration of the corresponding PSCA-based algorithm with the step size as a tunable parameter. \red{Each PSCA-Net elegantly marries the parallel computation mechanism of the PSCA-based algorithm} with the parallelizable neural network architecture. Besides, using neural network training, \red{each PSCA-Net effectively optimizes the step sizes of the PSCA-based algorithm} based on vast data samples. Therefore, the PSCA-Nets can successfully reduce the computation times of the PSCA-based algorithms.

(iv) Numerical results show that \red{the} PSCA-based algorithms can reduce the error rates and computation times by up to 70\%$\sim$80\% and 80\%$\sim$95\%, respectively, compared to the state-of-the-art algorithms under the same condition. Besides, compared to PSCA-ML-K (UD), PSCA-MAP-K (UR) can reduce the error rate but increase the computation time by up to 21.4\% and 33.3\%, respectively. Last, each PSCA-Net can further reduce the error rate and computation time by up to 90.8\% and 50.1\%, respectively, compared to the counterpart PSCA-based algorithm. The promising gains of the proposed methods reveal their signiﬁcant practical values in grant-free access.

The rest of this paper is organized as follows. \red{In Section~\ref{system},} we present the system model. In Section \ref{Known}, we consider the known pathloss case and investigate the MLE and MAPE to estimate each device's activity for device activity detection. 
In Section \ref{Unknown}, we consider the unknown pathloss case and investigate the MLE and MAPE to estimate each device's effective pathloss for device activity detection. In Section~\ref{Numerical_Results}, we numerically demonstrate the advantages of the proposed methods. In Section~\ref{conclusion}, we conclude this paper. We summarize the notation in Table~\ref{table_abbr}.

\begin{table}[!ht]
	\caption{{Methods and Abbreviations.}} 
	\label{table_abbr}
	\centering
	\begin{tabular}{|c|c|c|c|}
	\hline
	\makecell{\textbf{Pathloss} \\ \textbf{Cases}} & \makecell{\textbf{Estimation} \\ \textbf{Problems}} & \textbf{Methods} & \textbf{Abbreviations} \\
	\hline
	\multirow{8}{*}{\makecell{Known\\ pathloss}} 
	& \multirow{4}{*}{MLE} & BCD  & BCD-ML-K \\
	\cline{3-4}
	& & PG & PG-ML-K \\
	\cline{3-4}
	& & PSCA & PSCA-ML-K \\
	\cline{3-4}
	& & PSCA-Net & PSCA-ML-K-Net \\
	\cline{2-4}
	& \multirow{3}{*}{MAPE} & BCD  & BCD-MAP-K \\
	\cline{3-4}
	& & PSCA & PSCA-MAP-K \\
	\cline{3-4}
	& & PSCA-Net & PSCA-MAP-K-Net \\
	\cline{2-4}
	& MMSE & AMP  & AMP-K \\
	\hline
	\multirow{7}{*}{\makecell{Unknown\\ pathloss}} 
	& \multirow{4}{*}{MLE} & BCD  & BCD-ML-UD \\
	\cline{3-4}
	& & PG  & PG-ML-UD \\
	\cline{3-4}
	& & PSCA & PSCA-ML-UD \\
	\cline{3-4}
	& & PSCA-Net & PSCA-ML-UD-Net \\
	\cline{2-4}
	& \multirow{2}{*}{MAPE} & PSCA & PSCA-MAP-UR \\
	\cline{3-4}
	& & PSCA-Net & PSCA-MAP-UR-Net \\
	\cline{2-4}
	& MMSE & EM-AMP & EM-AMP-UD \\
	\hline
	\end{tabular}
	\color{black} 
\end{table}

\textit{Notations:} 
We denote scalars, vectors, and matrices by \red{lower-case letters, lower-case bold letters, and upper-case bold letters,} respectively. We use $(\cdot)^T$ and $(\cdot)^H$ to denote the transpose and conjugate transpose, respectively. We use $(\cdot)^{-1}$, $|\cdot|$, and $\operatorname{tr}\left(\cdot\right)$ to denote the inverse, determinant, and trace, respectively. We use $\mathbb{R}, \mathbb{R}_+, \mathbb{C}$, and $\mathbb{Z}_+$ to denote the set of real, positive real, complex, and natural numbers, respectively. We use $\Re (\cdot)$ and $\Im (\cdot)$ to denote the real and imaginary parts, respectively. \red{We denote the diagonal matrix generated by a vector by $\mathbf{Diag}(\cdot)$.  We use $\mathbf{I}_N$ and $\mathbf{0}_N$ to denote the $N \times N$ identity matrix and $N \times 1$ all-zero vector, respectively.}  
\red{We use $\nabla(\cdot)$ to denote the gradient of a function of several variables, and we use $\partial_i(\cdot)$ to denote the partial derivative of a function of several variables with respect to (w.r.t.) the $i$-th coordinate. We use $(\cdot)'$ and $(\cdot)''$ to denote the first-order and second-order derivatives of a function of one variable, respectively.}
We use $\left\Vert \cdot \right\Vert_p$ to denote the $l_p$-norm. We denote the indicator function by $\mathbb{I}[\cdot]$. \red{We use $\textrm{Pr}[\cdot]$ to denote the probability.}

%% file: system.tex
\section{System Model}
\label{system}
\begin{figure}[tp]
\centering
\begin{minipage}{0.49\textwidth}
\centering
\subfigure
		{\resizebox{4.0cm}{!}{\includegraphics{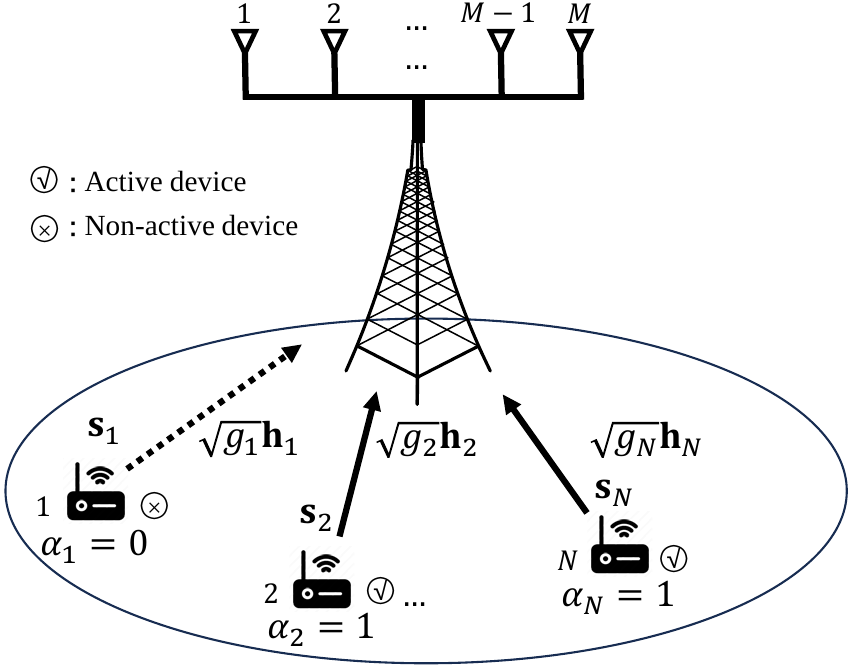}}}
  \caption{\small{Illustration of the system model.}}
	\label{sys}
\end{minipage}
\end{figure}

This paper investigates the uplink of a wireless network, which consists of an $M$-antenna BS and $N$ single-antenna IoT devices and operates on a narrow band system, as illustrated in Fig.~\ref{sys}. Denote $\mathcal M\triangleq\{1,2,..., M\}$ and $\mathcal N\triangleq\{1,2,..., N\}$, respectively.
In mMTC and URLLC, IoT devices are occasionally active and send small packets to the BS once activated. For all $n \in \mathcal{N}$, let $\alpha_n \in\{0,1\}$ denote the activity state of device $n$, where $\alpha_n=1$ indicates that device $n$ is active, and $\alpha_n=0$ otherwise. Denote $\boldsymbol{\alpha} \triangleq \left(\alpha_n\right)_{n \in \mathcal{N}} \in\{0,1\}^N$. Assume $\boldsymbol{\alpha}$ is unknown to the BS. In mMTC, $N$ is large, and $\boldsymbol{\alpha}$ is sparse, i.e., $\sum_{n \in \mathcal{N}} \alpha_n \ll N$. In contrast, in URLLC, $N$ is not necessarily large, and $\boldsymbol{\alpha}$ is not necessarily sparse. For all $n \in \mathcal{N}$, let $g_n \in \mathbb{R}_{+}$ denote the pathloss of the channel between device $n$ and the BS. Denote $\mathbf{g} \triangleq\left(g_n\right)_{n \in \mathcal{N}} \in \mathbb{R}^N_{+}$. In this paper, we consider two pathloss cases, i.e., \red{the} known pathloss case and unknown pathloss case. We consider the block flat fading model for small-scale fading. For all $n \in \mathcal{N}$, let $\mathbf{h}_n \triangleq\left(h_{n, m}\right)_{m \in \mathcal{M}} \in \mathbb{C}^M$ denote the small-scale fading coefficients of the channel between the BS and device $n$ in a coherence block, where $h_{n, m}$ represents the small-scale fading coefficient of the channel between the BS's $m$-th antenna and device $n$. Let $\mathbf{H} \triangleq\left(\mathbf{h}_n\right)_{n \in \mathcal{N}} \in \mathbb{C}^{M \times N}$. Assume $\mathbf{H}$ is unknown to the BS. Consequently, the overall channel between device $n$ and the BS is expressed as $\sqrt{g_n} \mathbf{h}_n \in \mathbb{C}^M$.

We consider a typical grant-free access scheme, which has two consecutive phases. Specifically, each device \red{$n \in \mathcal{N}$} is assigned a specific pilot sequence $\mathbf{s}_n\in \mathbb C^{L}$ of length $L$. Denote
$\mathbf{S} \triangleq\left(\mathbf{s}_n\right)_{n \in \mathcal{N}}\in \mathbb{C}^{L\times N}$, which is known to the BS. Since $L \ll N$, the pilot sequences are non-orthogonal. For simplicity, the time and frequency are assumed to be perfectly synchronized \cite{fengler2021nonbayesian,chen2019covariance,jiang2023statistical}. 
In the pilot transmission phase, active devices simultaneously send their pilot signals over the same spectrum, and the BS detects the activity states of all devices and estimates the active devices' channel states from the received pilot signals. In the following data transmission phase, active devices simultaneously send their data over the same spectrum, and the BS detects all active devices' data from the received data signals using received beamforming.

\red{Specifically, the received pilot \red{signals} at the BS in the pilot transmission phase, represented as $\mathbf{Y} \in \mathbb C^{L\times M}$, are given by:
\begin{align}
	\label{signal}
	\mathbf{Y} = \sum_{n \in \mathcal{N}}\alpha_n \sqrt{g_n} \mathbf{s}_n \mathbf{h}_n^T + \mathbf Z,
\end{align}
where $\mathbf Z  \in \mathbb C^{ L  \times M} $ is the additive white Gaussian noise (AWGN) at the BS with all elements following i.i.d. $\mathcal {CN}\left(0, \sigma^2\right)$.\footnote{Let $B$ denote the bandwidth of the (passband) narrow-band system. The power of each discrete-time pilot symbol depends on the continuous-time pilot signal's power and the sampling period $1/B$. Besides, the noise power is determined by the bandwidth $B$ and noise power spectral density. The pilot and noise powers influence the likelihood function and serve as the MLE and MAPE problems' parameters, affecting their stationary points. Therefore, $B$ \red{implicitly} influences the detection accuracy.} \red{In the known pathloss case, we rewrite $\mathbf{Y}$ in (\ref{signal}) in a compact form as:
\begin{align}
	\label{signal1}
	\mathbf{Y} = \mathbf{S} \boldsymbol\Gamma \mathbf{W}^{\frac{1}{2}}\mathbf{H} +\mathbf{Z},
\end{align}
where $\boldsymbol\Gamma \triangleq  \mathbf{Diag}\left(\boldsymbol \alpha\right) \in \mathbb R^{N\times N}$ and $\mathbf{W} \triangleq  \mathbf{Diag}\left(\mathbf{g}\right) \in \mathbb R^{N\times N}$, and estimate the device activity states $\boldsymbol \alpha$ from $\mathbf{Y}$ in (\ref{signal1}) for device activity detection. In the unknown pathloss case, we let $\gamma_n \triangleq \alpha_n g_n$ represent the effective pathloss between device $n$ and the BS, for all $n \in \mathcal{N}$, rewrite $\mathbf{Y}$ in (\ref{signal}) in another compact form as:\footnote{\red{Note that $\alpha_n \sqrt{g_n} = \sqrt{\alpha_n g_n}$ due to $\alpha_n \in\{0,1\}$, for all $n \in \mathcal{N}$.}}
\begin{align}
    	\label{signal2}
    	\mathbf{Y} =\mathbf{S} \boldsymbol\Lambda^{\frac{1}{2}}\mathbf{H} +\mathbf{Z},
    \end{align}
where $\boldsymbol\Lambda \triangleq  \mathbf{Diag}\left(\boldsymbol \gamma\right) \in \mathbb R^{N\times N}$ with $\boldsymbol{\gamma} \triangleq \left(\gamma_n\right)_{n \in \mathcal{N}}$, and estimate the effective pathloss $\boldsymbol{\gamma}$ from $\mathbf{Y}$ in (\ref{signal2}) for device activity detection.}}
    

Throughout the paper, we model the unknown $\mathbf{H}$ by Rayleigh fading, i.e., all elements of $\mathbf{h}_n, n \in \mathcal{N}$ are independently and identically distributed (i.i.d.) according to $\mathcal{C N}\left(0,1\right)$ \cite{fengler2021nonbayesian,chen2019covariance,dong2022faster,lin2023sparsity,jiang2020mapbased,jiang2022ml,shao2020covariancebased,wang2021accelerating,jiang2023statistical,liu2018massive}. We adopt different models for the unknown $\boldsymbol \alpha$ and $\mathbf{g}$, when designing different device activity detection methods.

The device activities $\boldsymbol{\alpha} \in\{0,1\}^N$ are unknown and to be estimated. There exist two models for unknown $\boldsymbol{\alpha}$.
One is to model $\boldsymbol{\alpha}$ as deterministic but unknown constants, and the other is to model $\boldsymbol{\alpha}$ as realizations of Bernoulli random variables $\mathbf{A} \triangleq \left(A_n\right)_{n \in \mathcal{N}} \in \{0,1\}^N$ to incorporate the prior knowledge of \red{device activities} \cite{jiang2022ml,jiang2023statistical}. Specifically, we consider two distribution models for $\mathbf{A}$, i.e., the general model \red{(or its approximation)} and the independent model. We adopt the multivariate Bernoulli (MVB) model as the general model for $\mathbf{A}$ \cite{ding2011learning}, and the probability mass function (p.m.f.) of $\mathbf{A}$, denoted by $p_{\mathbf{A}\textrm{\rm-G}}\left(\boldsymbol{\alpha}\right)$, is given by (\ref{MVB}) \cite{jiang2022ml,jiang2023statistical}, as shown at the top of the next page.
\begin{figure*}[!ht]
	\begin{align}
		\label{MVB}
		p_{\mathbf{A}\textrm{\rm-G}}\left(\boldsymbol{\alpha}\right)&=\exp\left(\sum_{\omega\in\Psi}c_{\omega}\prod_{n\in\omega}\alpha_n-\log \left(\sum_{\boldsymbol{\alpha'}\in\{0,1\}^{N}}\exp \left (\sum_{\omega\in\Psi} c_{\omega}\prod_{n\in\omega}\alpha_n'\right)\right)\right).
  \\ \label{second}p_{\mathbf{A}\textrm{\rm-G}}\left(\boldsymbol{\alpha}\right)&=\exp\left(\sum_{\substack{n_1<n_2\\n_1,n_2\in\mathcal{N}}}c_{\{n_1,n_2\}}\alpha_{n_1}\alpha_{n_2} +\sum_{\substack{n\in\mathcal{N}}}c_{\{n\}}\alpha_{n} -\log \left(\sum_{\boldsymbol{\alpha}'\in\{0,1\}^{N}}\exp \left (\sum_{\substack{n_1<n_2\\n_1,n_2\in\mathcal{N}}}c_{\{n_1,n_2\}}\alpha_{n_1}'\alpha_{n_2}' +\sum_{\substack{n\in\mathcal{N}}}c_{\{n\}}\alpha_{n}'\right)\right)\right).
	\end{align}
	\hrulefill
\end{figure*}
In (\ref{MVB}), $\Psi$ denotes the set of nonempty subsets of $\mathcal{N}$, and \red{for all $\omega\in\Psi, c_{\omega}$} is the coefficient reflecting the correlation among $\alpha_n$, $n\in\mathcal{N}$. The parameters for (\ref{MVB}) are $c_{\omega}, \omega \in \Psi$. In the independent model for $\mathbf{A}$, we assume that $\mathbf{A}$ are independent Bernoulli random variables, and the p.m.f. of $\mathbf{A}$, denoted by $p_{\mathbf{A}\textrm{\rm-I}}\left(\boldsymbol{\alpha}\right)$, is given by \cite{jiang2022ml,jiang2023statistical}:
\begin{align}
	\label{independent}
	p_{\mathbf{A}\textrm{\rm-I}}\left(\boldsymbol{\alpha}\right)=\prod_{n \in\mathcal{N}}\left(p_n\right)^{\alpha_n}\left(1-p_n\right)^{1-\alpha_n},
\end{align} where $p_n \triangleq \textrm{Pr}\left[\alpha_n=1\right]$ denotes the active probability of device $n$. The model parameters for the independent model given by (\ref{independent}) are $p_n, n \in\mathcal{N}$.   
In the MVB model, \( c_\omega, \omega \in \Psi \) can be estimated based on the historical device activity data using existing methods \cite{ding2011learning}.  In addition, given the p.m.f. of $\mathbf{A}$ in any
 form, the coefficients \( c_\omega, \omega \in \Psi \) can be calculated according to \cite[Lemma 2.1]{ding2011learning}.
When $c_{\omega} = 0$ for all $|\omega| > 1$, from (\ref{MVB}), we have:
\begin{align*}
	& \quad p_{\mathbf{A}\textrm{\rm-G}}\left(\boldsymbol{\alpha}\right)  = \frac{\exp\left(\sum_{n \in  \mathcal{N}}c_{\{n\}}\alpha_n\right)}{\prod_{n \in  \mathcal{N}}\left(\exp\left(c_{\{n\}}\right)+1\right)} \\
	& = \prod_{n \in\mathcal{N}}\left(\frac{\exp\left(c_{\{n\}}\right)}{\exp\left(c_{\{n\}}\right)+1}\right)^{\alpha_n} \left(1-\frac{\exp\left(c_{\{n\}}\right)}{\exp\left(c_{\{n\}}\right)+1}\right)^{1-\alpha_n}.
\end{align*}
That is, the general model in  (\ref{MVB}) reduces to the independent model in (\ref{independent}) with $p_n = \frac{\exp\left(c_{\{n\}}\right)}{\exp\left(c_{\{n\}}\right)+1}$, i.e., $c_{\{n\}} = \log\left(\frac{p_n}{1-p_n}\right), n \in \mathcal{N}$. 

For some sophisticated distributions of $\mathbf{A}$ whose general model $p_{\mathbf{A}\textrm{\rm-G}}\left(\boldsymbol{\alpha}\right)$ is not tractable (i.e., most $c_{\omega}, \omega \in \Psi$ are non-zero), we consider the first-order or \red{second-order} approximation, instead.  Specifically, by setting $c_{\omega} = 0$ for all $\omega$ with  $|\omega| > 1$, we obtain the first-order approximation, i.e., the independent model $p_{\mathbf{A}\textrm{\rm-I}}\left(\boldsymbol{\alpha}\right)$ in (\ref{independent}), which captures the marginal active probability of every single device. By setting $c_{\omega} = 0$ for all $\omega$ with  $|\omega| > 2$, we obtain the second-order approximation, which is given by (\ref{second}) \red{(with a slight abuse of notation)}, as shown at the top of the page, capturing the active probabilities of every single device and every two devices. 

The pathloss $\mathbf{g}$ can be categorized into \red{the} known and unknown cases. There exist two models for unknown $\mathbf{g}$. One is to regard $\mathbf{g}$ as deterministic but unknown constants \cite{fengler2021nonbayesian,wang2021efficient}. The other is to view $\mathbf{g}$ as realizations of  random variables $\mathbf{G} \triangleq \left(G_n\right)_{n \in \mathcal{N}} \in \mathbb{R}^N_{+}$ \cite{zhang2024activity}.  
Specifically, we assume that for all $n \in \mathcal{N}$, device $n$ is independently and uniformly located in the annulus with an inner radius \red{$d_{l,n}$ and an outer radius $d_{u,n}$} \cite{haenggi2012} and centered at the BS's location. For all $n \in \mathcal{N}$, let $D_n$ denote the distance between device $n$ and the BS. Accordingly, $D_n \in [d_{l, n}, d_{u, n}], n \in \mathcal{N}$ are independent. \red{We assume that $d_{l, n}, d_{u, n}, n \in \mathcal{N}$ are known to the BS.}
The probability density function (p.d.f.) of $D_n$, denoted by $p_{D_n}\left(d_n\right)$,  is given by \cite[Sec. 3.1]{bertsekas2008introduction}:
\begin{align}
	\label{dispdf}
p_{D_n}\left(d_n\right) = \frac{2d_n}{d_{u, n}^2-d_{l, n}^2}, \ d_n \in [d_{l, n}, d_{u, n}], \ n \in \mathcal{N}.
\end{align} 
 Denote $\mathbf{D} \triangleq \left(D_n\right)_{n \in \mathcal{N}} \in \mathbb{R}^N_{+}$. By the independence of $\mathbf{D}$ and (\ref{dispdf}), the p.d.f. of $\mathbf{D}$, denoted by $p_{\mathbf{D}}\left(\mathbf{d}\right)$, is given by:
\begin{align}
p_{\mathbf{D}}\left(\mathbf{d}\right) = \prod_{n \in\mathcal{N}}\frac{2d_n}{d_{u, n}^2-d_{l, n}^2},\  d_n \in [d_{l, n},d_{u, n}], \ n \in \mathcal{N}.
\end{align}
Furthermore, for all $n \in \mathcal{N}$, we express the pathloss of device $n$, i.e., $G_n$, as a function of $D_n$:
\begin{align}
\label{gn}
    G_n = \phi\left(D_n\right),
\end{align}
where $\phi: \mathbb{R}_{+} \to \mathbb{R}_{+}$ is a monotonically decreasing function mapping the distance to pathloss and is assumed to be second-order differentiable.
By (\ref{dispdf}) and (\ref{gn}), we derive the p.d.f. of $G_n$, denoted by $p_{G_n}\left(g_n\right)$, is given by \cite[Sec. 4.1]{bertsekas2008introduction}:
\begin{align}
	\label{pg}
	&p_{G_n}\left(g_n\right) = -\frac{2\phi^{-1}\left(g_n\right)\left(\phi^{-1}\right)'\left(g_n\right)}{d_{u, n}^2-d_{l, n}^2}, \nonumber \\ 
	&g_n \in [g_{l,n},g_{u,n}],\  n \in \mathcal{N},
\end{align}
where $g_{l,n} \triangleq \phi\left(d_{u, n}\right)$ and $g_{u,n} \triangleq \phi\left(d_{l, n}\right)$, respectively.
See Appendix~\ref{appendix:pdf} for the proof.
By the independence of $\mathbf{D}$ ($\mathbf{G}$) and (\ref{pg}), the p.d.f. of $\mathbf{G}$, denoted by $p_{\mathbf{G}}\left(\mathbf{g}\right)$, is given by:
\begin{align}
	\label{pG}
p_{\mathbf{G}}\left(\mathbf{g}\right) = \prod_{n \in\mathcal{N}}p_{G_n}\left(g_n\right), \ g_n \in [g_{l,n},g_{u,n}], \ n \in \mathcal{N},
\end{align} 
where $p_{G_n}\left(g_n\right)$ is given by (\ref{pg}), and $\phi^{-1}$ denotes the inverse function of $\phi$.

\begin{figure*}[!hbp]
	\hrulefill
	\begin{align}
\label{pgamma}
p_{\Upsilon_n}\left(\gamma_n\right) & = \left(1-p_n\right)^{\mathbb{I}\left[\gamma_n \leq 0\right]}\left(-\frac{2p_n}{d_{u, n}^2-d_{l, n}^2}\phi^{-1}\left(\gamma_n\right)\left(\phi^{-1}\right)'\left(\gamma_n\right)\right)^{1-\mathbb{I}\left[\gamma_n \leq 0\right]}, \ \gamma_n \in \{0\} \cup \left[g_{l,n},g_{u,n}\right], \ n \in \mathcal{N}.\\
 \label{H2}
H_2\left(\gamma_n\right) & \triangleq  -\left(\frac{2p_n d_{u, n}}{d_{u, n}^2-d_{l, n}^2}\left(\phi^{-1}\right)'\left(g_{l,n}\right)\right)\left(1-\frac{2\left(\gamma_n-g_{l,n}\right)}{\epsilon}\right)\left(\frac{\gamma_n-g_{l,n}+\epsilon}{\epsilon}\right)^2  \nonumber \\
& \quad 
-\left(\frac{2p_n}{d_{u, n}^2-d_{l, n}^2}\left(\left(\left(\phi^{-1}\right)'\left(g_{l,n}\right)\right)^2 +\left(\phi^{-1}\right)''\left(g_{l,n}\right)d_{u, n}\right)\right)\left(\gamma_n-g_{l,n}\right)\left(\frac{\gamma_n-g_{l,n}+\epsilon}{\epsilon}\right)^2.
\end{align}
\end{figure*}


Given the two models for unknown $\boldsymbol{\alpha}$ and $\mathbf{g}$, and $\boldsymbol{\gamma} \triangleq \left(\gamma_n\right)_{n \in \mathcal{N}}, \gamma_n \triangleq \alpha_n g_n, n \in \mathcal{N}$, we consider two models for $\boldsymbol{\gamma}$. One is to view $\boldsymbol{\gamma}$ as deterministic but unknown constants. The other is to regard unknown $\boldsymbol{\gamma}$ as realizations of independent random variables $\boldsymbol{\Upsilon} \triangleq \left(\Upsilon_n\right)_{n \in \mathcal{N}} \in \mathbb{R}^N_{+}$, with\footnote{In this model, we do not consider the dependence among $A_n, n \in \mathcal{N}$ and $G_n, n \in \mathcal{N}$ for simplicity. The distributions for i.i.d. random variables $\mathbf{G}$, can be extended to other distributions for i.i.d. random variables with twice differentiable p.d.f.s.}
$\Upsilon_n \triangleq A_nG_n, n \in \mathcal{N}.$
By (\ref{independent}) and (\ref{pG}), the p.d.f. of $\Upsilon_n$, denoted by $p_{\Upsilon_n}\left(\gamma_n\right)$, is given by (\ref{pgamma}) \cite[Section 2.1]{bertsekas2008introduction}, as shown at the bottom of the next page. The details are presented in Appendix~\ref{appendix:pdf}. Notice that $p_{\Upsilon_n}\left(\gamma_n\right)$ in (\ref{pgamma}) is non-differential and has a disconnected domain, \red{hence intractable.} For tractability, we adopt a first-order approximation and domain relaxation via Hermite interpolation \cite{nocewrig2006numerical}. The approximate p.d.f. of $\Upsilon_n$, denoted by $p^{\epsilon }_{\Upsilon_n}\left(\gamma_n\right)$, is given by:
\begin{align}
\label{peps}
    &p^{\epsilon }_{\Upsilon_n}\left(\gamma_n\right) = \nonumber \\
    & 
 \begin{cases}
     H_1\left(\gamma_n\right), & 0 \leq \gamma_n < \epsilon, \\ 
     0, & \epsilon \le \gamma_n < g_{l,n} - \epsilon, \\
     H_2\left(\gamma_n\right), & g_{l,n} - \epsilon \le \gamma_n < g_{l,n}, \\
         - \frac{2p_n}{d_{u, n}^2 - d_{l, n}^2}\phi^{-1}\left(\gamma_n\right)\left(\phi^{-1}\right)'\left(\gamma_n\right), & g_{l,n} \le \gamma_n \le g_{u,n},
 \end{cases}
\end{align}
where $\epsilon \in (0,1]$ is a constant determining the approximate level, $H_1\left(\gamma_n\right) \triangleq  \left(1-p_n\right)\left(1+\frac{2\gamma_n}{\epsilon}\right)\left(\frac{\gamma_n-\epsilon}{\epsilon}\right)^2$, \red{and} $H_2\left(\gamma_n\right)$ is given by (\ref{H2}), as shown at the bottom of the next page.
The properties of $p^{\epsilon }_{\Upsilon_n}\left(\gamma_n\right), n \in \mathcal{N}$ are presented below.
\begin{Thm}[Properties of $p^{\epsilon }_{\Upsilon_n}\left(\gamma_n\right)$]
\label{properties}
(i) $p^{\epsilon }_{\Upsilon_n}\left(\gamma_n\right) =  p_{\Upsilon_n}\left(\gamma_n\right)$, for all $\gamma_n \in \{0\} \cup \left[g_{l,n},g_{u,n}\right]$. (ii) $p^{\epsilon }_{\Upsilon_n}\left(\gamma_n\right)$ is differentiable in $\left[0, g_{u,n}\right]$. (iii) $p^{\epsilon }_{\Upsilon_n}\left(\gamma_n\right)$ converges to $\widehat{p}_{\Upsilon_n}\left(\gamma_n\right)$, as $\epsilon\to0$, for all $\gamma_n \in \left[0, g_{u,n}\right]$,  where $\widehat{p}_{\Upsilon_n}\left(\gamma_n\right)$ is given by:
\begin{align}
\label{hatp}
\widehat{p}_{\Upsilon_n}\left(\gamma_n\right) = 
\begin{cases}
p_{\Upsilon_n}\left(\gamma_n\right), & \{0\} \cup \left[g_{l,n},g_{u,n}\right], \\
0, & \gamma_n \in \left(0, g_{l,n}\right).
\end{cases}
\end{align}
\end{Thm}
\begin{IEEEproof}
    Please refer to Appendix \ref{appendix:approximate}.
\end{IEEEproof}

By Thoerem~\ref{properties}, $p^{\epsilon }_{\Upsilon_n}\left(\gamma_n\right)$ is more tractable than $p_{\Upsilon_n}\left(\gamma_n\right)$ and can well approximate $p_{\Upsilon_n}\left(\gamma_n\right)$.
\red{Therefore, we consider the approximate p.d.f. of $\boldsymbol{\Upsilon}$, denoted by $p^{\epsilon }_{\boldsymbol{\Upsilon}}\left(\boldsymbol{\gamma}\right)$ in the rest of the paper. By the independence of $\Upsilon_n, n \in \mathcal{N}$ and (\ref{peps}), $p^{\epsilon }_{\boldsymbol{\Upsilon}}\left(\boldsymbol{\gamma}\right)$ is given by:}
\begin{align}
    \label{pGamma}	p^{\epsilon }_{\boldsymbol{\Upsilon}}\left(\boldsymbol{\gamma}\right) = \prod_{n \in\mathcal{N}}p^{\epsilon }_{\Upsilon_n}\left(\gamma_n\right), \ \gamma_n \in \left[0, g_{u,n}\right], \ n \in \mathcal{N},
\end{align}
where $p^{\epsilon }_{\Upsilon_n}\left(\gamma_n\right)$ is given by (\ref{peps}).

By leveraging the deterministic and random models of unknown $\boldsymbol{\alpha}$ in the known pathloss case and unknown $\boldsymbol{\gamma}$ in the unknown pathloss case, we can design different device activity detection methods.
In Section \ref{Known}, we consider the known pathloss case and investigate the MLE and MAPE for $\boldsymbol{\alpha}$ for device activity detection in Sections \ref{MLEK} and \ref{MAPEK}, respectively. 
In Section \ref{Unknown}, we consider the unknown pathloss case and investigate the MLE and MAPE for $\boldsymbol{\gamma}$ for device activity detection in Sections \ref{MLEU} and \ref{MAPEU}, respectively.\footnote{\red{Problem~\ref{Prob_ML} \cite{fengler2021nonbayesian}, Problem~\ref{Prob_MAP} \cite{jiang2022ml,jiang2020mapbased}, and Problem~\ref{Prob_ML_U} \cite{fengler2021nonbayesian,wang2021efficient} have been investigated.} We present them here for completeness. Problem~\ref{Prob_MAP_U} has not yet been studied.}

%% file: known.tex
\section{MLE and MAPE-Based Device Activity Detection for Known Pathloss}
\label{Known}
In this section, we investigate MLE and MAPE-based device activity detection in the known pathloss case.\footnote{In the known pathloss case, the MLE-based and MAPE-based device activity estimates converge to the real device activity states in probability, as $M$ goes to infinity \cite{chen2019covariance,kay1993fundamentalsOS}.}

\subsection{MLE-Based Device Activity Detection for Known Pathloss}
\label{MLEK}
In this subsection, we view $\boldsymbol{\alpha}$ as deterministic but unknown constants and investigate the MLE of $\boldsymbol{\alpha}$ for known $\mathbf{g}$.
\subsubsection{MLE Problem for Known Pathloss}
\label{MLE_Problem}
In this part, we present the MLE problem of $\boldsymbol{\alpha}$ for known $\mathbf{g}$.
Given device activities $\boldsymbol\alpha, \mathbf{Y}_{:, m}, m \in \mathcal{M}$ follow i.i.d. $C \mathcal{N}\left(0, \boldsymbol{\Sigma}_{\boldsymbol{\alpha}}\right)$ with the covariance matrix \cite{fengler2021nonbayesian,jiang2022ml}: 
\begin{align}
    \label{cov_a}
    \boldsymbol{\Sigma}_{\boldsymbol{\alpha}} \triangleq \mathbf{S}\boldsymbol{\Gamma}\mathbf{W}\mathbf{S}^H+\sigma^2 \mathbf{I}_L \in \mathbb{C}^{L \times L}.
\end{align}
Therefore, the likelihood of $\mathbf{Y}$, denoted by $p_{\textrm{\rm ML-K}}(\mathbf{Y};\boldsymbol\alpha)$, is:
\begin{align}
\label{fY}
	p_{\textrm{\rm ML-K}}(\mathbf{Y};\boldsymbol\alpha) \triangleq \frac{\exp \left(-\operatorname{tr}\left(\boldsymbol{\Sigma}_{\boldsymbol{\alpha}}^{-1} \mathbf{Y Y}^H\right)\right)}{\pi^{LM}\left| \boldsymbol{\Sigma}_{\boldsymbol{\alpha}} \right|^M}.
\end{align}
Let
\begin{align}
	\label{likelihood_Y}
	f_{\textrm{\rm ML-K}}\left(\boldsymbol \alpha\right) & \triangleq  -\frac{1}{M} \log{p_{\textrm{\rm ML-K}}(\mathbf{Y};\boldsymbol\alpha)} - L \log{\pi} \nonumber \\ &= \log \left| \boldsymbol{\Sigma}_{\boldsymbol{\alpha}} \right| + \operatorname{tr}\left(\boldsymbol{\Sigma}_{\boldsymbol{\alpha}}^{-1} \widehat{\boldsymbol{\Sigma}}_{\mathbf{Y}}\right),
\end{align} 
where $\widehat{\boldsymbol{\Sigma}}_{\mathbf{Y}} \triangleq \frac{1}{M} \mathbf{Y Y}^H \in \mathbb{C}^{L \times L}$ denotes the empirical covariance matrix, which is the sufficient statistics for estimating $\boldsymbol{\alpha}$.
The MLE problem of $\boldsymbol{\alpha}$, i.e., the maximization of $p_{\textrm{\rm ML-K}}(\mathbf{Y};\boldsymbol\alpha)$ or $\log p_{\textrm{\rm ML-K}}(\mathbf{Y};\boldsymbol\alpha)$ \red{w.r.t. $\boldsymbol{\alpha}$,} can be equivalently formulated below.\footnote{In Problem~\ref{Prob_ML} and Problem~\ref{Prob_MAP}, $\alpha_n \in\{0,1\}$ is relaxed to $\alpha_n \in [0,1]$, for all $n \in \mathcal{N}$, and the binary detection results are obtained by performing thresholding as in \cite{fengler2021nonbayesian,chen2019covariance,jiang2022ml,cui2021jointly}.}
\begin{Prob}[MLE of $\boldsymbol\alpha$ for Known $\mathbf{g}$]
	\label{Prob_ML}
	\begin{align}
		\min_{\boldsymbol\alpha} & \ f_{\textrm{\rm ML-K}}\left(\boldsymbol \alpha\right) \nonumber \\ \label{constraint}
		\text {s.t.} & \ \alpha_n \in [0,1], \ n \in \mathcal{N}.
	\end{align}
\end{Prob}

Problem~\ref{Prob_ML} is a challenging non-convex problem, which generally has no effective algorithms. Existing state-of-the-art iterative algorithms for obtaining stationary points of Problem~\ref{Prob_ML} include BCD-ML-K \cite{chen2019covariance,jiang2022ml,fengler2021nonbayesian,jiang2023statistical} and PG-ML-K \cite{wang2021efficient}. As illustrated in Section \ref{Introduction}, BCD-ML-K has a long computation time due to its sequential update and cannot leverage modern multi-core processors. PG-ML-K successfully reduces the computation time, leveraging the parallel update mechanism. However, PG-ML-K, which utilizes up to the first-order information of the objective function $f_{\textrm{\rm ML-K}}\left(\boldsymbol \alpha\right)$ of Problem \ref{Prob_ML}, still has unsatisfactory convergence speed and computation time.

\subsubsection{PSCA-ML-K}
\label{PSCA_ML_Algorithm}
In this part, we propose an algorithm, referred to as PSCA-ML-K, using PSCA. The main idea of PSCA-ML-K is to solve a sequence of successively refined approximate convex optimization problems, each of which can be equivalently separated into $N$ subproblems that are analytically solved in parallel. Specifically, at iteration $k$, we choose:
\begin{align}
	\label{Surrogate_ML}
	\widetilde{f}_{\textrm{\rm ML-K}}\left(\boldsymbol{\alpha};\boldsymbol{\alpha}^{(k)}\right) \triangleq \sum_{n \in \mathcal{N}} & \widetilde{f}_{\textrm{\rm ML-K}, n}\left(\alpha_n ;\boldsymbol{\alpha}^{(k)}\right) 
\end{align} as a strongly convex approximate function of $f_{\textrm{\rm ML-K}}\left( \boldsymbol{\alpha}\right)$ around $\boldsymbol{\alpha}^{(k)}$ \red{w.r.t. $\boldsymbol{\alpha}$}, where
\begin{align}
	\label{Surrogate_ML_separate}
	\widetilde{f}_{\textrm{\rm ML-K}, n}(\alpha_n;\boldsymbol{\alpha}^{(k)})  \triangleq & \partial_{n} f_{\textrm{\rm ML-K}}\left(\boldsymbol{\alpha}^{(k)}\right) \left(\alpha_n - \alpha_n^{(k)}\right) \nonumber \\  + & \frac{1}{2}\left(g_n \mathbf{s}_n^H \boldsymbol{\Sigma}_{\boldsymbol{\alpha}^{(k)}}^{-1} \mathbf{s}_n\right)^2 \left(\alpha_n - \alpha_n^{(k)}\right)^2
\end{align} is a strongly convex approximate function of $f_{\textrm{\rm ML-K}}\left(\boldsymbol{\alpha}\right)$ around $\boldsymbol{\alpha}^{(k)}$ w.r.t. $\alpha_n$, \red{and $\boldsymbol{\alpha}^{(k)}$ is obtained at iteration $k-1$}. Here,
\begin{align}
	\label{gradient_ML}
	&\partial_{n} f_{\textrm{\rm ML-K}}\left(\boldsymbol{\alpha}^{(k)}\right) \nonumber \\ &= g_n \left(\mathbf{s}_n^H \boldsymbol{\Sigma}_{\boldsymbol{\alpha}^{(k)}}^{-1} \mathbf{s}_n - \mathbf{s}_n^H \boldsymbol{\Sigma}_{\boldsymbol{\alpha}^{(k)}}^{-1} \widehat{\boldsymbol{\Sigma}}_{\mathbf{Y}} \boldsymbol{\Sigma}_{\boldsymbol{\alpha}^{(k)}}^{-1} \mathbf{s}_n \right).
\end{align} 
The quadratic coefficient $\left(g_n \mathbf{s}_n^H \boldsymbol{\Sigma}_{\boldsymbol{\alpha}^{(k)}}^{-1} \mathbf{s}_n\right)^2$ in (\ref{Surrogate_ML_separate}) converges to the $n$-th diagonal element of the Hessian matrix $\nabla^2 f_{\textrm{\rm ML-K}} (\boldsymbol{\alpha}^{(k)})$, i.e., $ g_n^2\left(\mathbf{s}_n^H \boldsymbol{\Sigma}_{\boldsymbol{\alpha}^{(k)}}^{-1} \mathbf{s}_n\right)\left(2\mathbf{s}_n^H \boldsymbol{\Sigma}_{\boldsymbol{\alpha}^{(k)}}^{-1} \widehat{\boldsymbol{\Sigma}}_{\mathbf{Y}} \boldsymbol{\Sigma}_{\boldsymbol{\alpha}^{(k)}}^{-1} \mathbf{s}_n - \mathbf{s}_n^H \boldsymbol{\Sigma}_{\boldsymbol{\alpha}^{(k)}}^{-1} \mathbf{s}_n\right)$, almost surely, as $M \to \infty$, since $\mathbf{s}_n^H \boldsymbol{\Sigma}_{\boldsymbol{\alpha}^{(k)}}^{-1} \widehat{\boldsymbol{\Sigma}}_{\mathbf{Y}} \boldsymbol{\Sigma}_{\boldsymbol{\alpha}^{(k)}}^{-1} \mathbf{s}_n$ converges to $\mathbf{s}_n^H \boldsymbol{\Sigma}_{\boldsymbol{\alpha}^{(k)}}^{-1} \mathbf{s}_n$, almost surely, as $M \to \infty$. Thus, the surrogate function $\widetilde{f}_{\textrm{\rm ML-K}}\left(\boldsymbol{\alpha};\boldsymbol{\alpha}^{(k)}\right)$ given in (\ref{Surrogate_ML})  approximately utilizes up to the second-order information of the objective function $f_{\textrm{\rm ML-K}}\left(\boldsymbol \alpha\right)$ of Problem~\ref{Prob_ML}. 
The approximate convex problem at iteration $k$ is given by:
\begin{align}
	\underset{\boldsymbol{\alpha} \in [0,1]^N}{\operatorname{min}}  \widetilde{f}_{\textrm{\rm ML-K}}\left(\boldsymbol{\alpha};\boldsymbol{\alpha}^{(k)}\right),
\end{align}
which can be equivalently separated into $N$ convex problems, one for each coordinate.
\begin{Prob}[Convex Approximate Problem of Problem~\ref{Prob_ML} for $n \in\mathcal{N}$ at Iteration $k$]
	\label{Prob_Surrogate_ML}\red{
\begin{align*}
    \widetilde{\alpha}_{\textrm{\rm ML-K}, n}^{(k)} \triangleq \underset{\alpha_n \in [0,1]}{\operatorname{argmin}} \  \widetilde{f}_{\textrm{\rm ML-K}, n}\left(\alpha_n;\boldsymbol{\alpha}^{(k)}\right).
\end{align*}}
\end{Prob}
\begin{figure*}[!ht]
\begin{align}
 \label{Frobenius_inversion_1}
\Re\left(\boldsymbol{\Sigma}_{\boldsymbol{\alpha}^{(k)}}^{-1}\right) = &\left(\Re \left(\boldsymbol{\Sigma}_{\boldsymbol{\alpha}^{\left(k\right)}}\right)+\Im \left(\boldsymbol{\Sigma}_{\boldsymbol{\alpha}^{\left(k\right)}}\right)\Re \left(\boldsymbol{\Sigma}_{\boldsymbol{\alpha}^{\left(k\right)}}\right)^{-1}\Im \left(\boldsymbol{\Sigma}_{\boldsymbol{\alpha}^{\left(k\right)}}\right)\right)^{-1}, \\ \label{Frobenius_inversion_2}\Im\left(\boldsymbol{\Sigma}_{\boldsymbol{\alpha}^{(k)}}^{-1}\right) = & \red{-}\Re \left(\boldsymbol{\Sigma}_{\boldsymbol{\alpha}^{\left(k\right)}}\right)^{-1}\Im \left(\boldsymbol{\Sigma}_{\boldsymbol{\alpha}^{\left(k\right)}}\right)\Big(\Re \left(\boldsymbol{\Sigma}_{\boldsymbol{\alpha}^{\left(k\right)}}\right) \Big.+ \left. \Im \left(\boldsymbol{\Sigma}_{\boldsymbol{\alpha}^{\left(k\right)}}\right)\Re \left(\boldsymbol{\Sigma}_{\boldsymbol{\alpha}^{\left(k\right)}}\right)^{-1}\Im \left(\boldsymbol{\Sigma}_{\boldsymbol{\alpha}^{\left(k\right)}}\right)\right)^{-1}.
	\end{align}
	\hrulefill
\end{figure*}
\begin{Lem}[Optimal Solution of Problem~\ref{Prob_Surrogate_ML}]
	\label{Lem_ML}
	The optimal point of Problem~\ref{Prob_Surrogate_ML} is given by:
	\begin{align}
		\label{Solution_ML}
		\widetilde{\alpha}_{\textrm{\rm ML-K}, n}^{(k)} = \min \left\{\max \left\{\alpha_n^{(k)} - \frac{\partial_{n} f_{\textrm{\rm ML-K}}\left(\boldsymbol{\alpha}^{(k)}\right)}{\left(g_n \mathbf{s}_n^H \boldsymbol{\Sigma}_{\boldsymbol{\alpha}^{(k)}}^{-1} \mathbf{s}_n \right)^2} ,0 \right \},1 \right \},
	\end{align}
	where $\partial_{n} f_{\textrm{\rm ML-K}}\left(\boldsymbol{\alpha}^{(k)}\right)$ is given by (\ref{gradient_ML}).
\end{Lem}
\begin{IEEEproof}
Please refer to Appendix~\ref{appendix:solution}.
\end{IEEEproof}

Note that $\widetilde{\alpha}_{\textrm{\rm ML-K}, n}^{(k)}$ in (\ref{Solution_ML}) approximately utilizes up to the second-order information of $f_{\textrm{\rm ML-K}}\left(\boldsymbol{\alpha}^{(k)}\right)$ via $\left(g_n \mathbf{s}_n^H \boldsymbol{\Sigma}_{\boldsymbol{\alpha}^{(k)}}^{-1} \mathbf{s}_n \right)^2$. Furthermore, PSCA-ML-K analytically solves $N$ convex approximate problems (Problem \ref{Prob_Surrogate_ML}) in parallel, whereas BCD-ML-K analytically and sequentially solves $N$ non-convex coordinate descent optimization problems. For any coordinate, the convex approximate problem, i.e., Problem \ref{Prob_Surrogate_ML} in PSCA-ML-K, and the non-convex coordinate descent problem in BCD-ML-K share the same optimal solution form, i.e., (\ref{Solution_ML}).

In each iteration $k$, we first compute the inverse of a complex positive definite matrix, i.e., $\boldsymbol{\Sigma}_{\boldsymbol{\alpha}^{(k)}}^{-1}$, by computing its real and imaginary parts according to (\ref{Frobenius_inversion_1}) and (\ref{Frobenius_inversion_2}), respectively, as shown at the top of the page, using Frobenius inversion.
Then, we are ready to compute $\widetilde{\alpha}_{\textrm{\rm ML-K}, n}^{(k)}, n \in \mathcal{N}$ in parallel according to (\ref{Solution_ML}). Finally, we update $\alpha_n^{(k)}, n \in \mathcal{N}$ in parallel according to:
\begin{align}
\label{Update_ML}
	\alpha_n^{(k+1)} = \left(1 - \rho^{(k)}\right)\alpha_n^{(k)} + \rho^{(k)} \widetilde{\alpha}_{\textrm{\rm ML-K}, n}^{(k)}, \ n \in \mathcal{N},
\end{align}
where $\rho^{(k)}$ is the step size of iteration $k$ satisfying:
\begin{align}
	\label{diminishing_rule}
	\rho^{(k)} \in (0,1], \ \sum_{k=0}^{\infty} \rho^{(k)}=+\infty, \ \sum_{k=0}^{\infty}\left(\rho^{(k)}\right)^2<+\infty.
\end{align}
The detailed procedure\footnote{All computations in each PSCA and PSCA-Net are implemented using vector and matrix operations wherever applicable to facilitate parallel computations.} is summarized in Algorithm \ref{Algorithm_PSCA_General} (PSCA-ML-K).\footnote{Each PSCA algorithm relies only on $\widehat{\boldsymbol{\Sigma}}_{\mathbf{Y}} = \frac{1}{M} \mathbf{Y Y}^H$ rather than $\mathbf{Y}$. $\widehat{\boldsymbol{\Sigma}}_{\mathbf{Y}}$ is computed only once, with the dominant term in the flop count \red{$8L^2M$}, and is treated as the input of each PSCA algorithm. The per-iteration computational complexity of each PSCA algorithm is irrelevant to $M$. Furthermore, $\boldsymbol\alpha^{(0)} = \mathbf{0}_N$ and $\boldsymbol\gamma^{(0)} = \mathbf{0}_N$ have been verified as good initial points for PSCAs due to the sparsity of device activities \cite{fengler2021nonbayesian,jiang2022ml,jiang2023statistical}.}

\begin{Thm}[Convergence of PSCA-ML-K]
	\label{Convergence_ML}
	Every limit point of the sequence $\left\{\boldsymbol{\alpha}^{(k)}\right\}_{k \in \mathbb{Z}_{+}}$ generated by PSCA-ML-K (at least one such point exists) is a stationary point of Problem~\ref{Prob_ML}.
\end{Thm}
\begin{IEEEproof}
	Please refer to Appendix~\ref{appendix:convergence}.
\end{IEEEproof}

In the following, we compare the main characteristics of PSCA-ML-K with BCD-ML-K \cite{fengler2021nonbayesian} and PG-ML-K \cite{wang2021efficient}. For ease of illustration, we present BCD-ML-K and PG-ML-K (without the active selection for simplicity) in Algorithm \ref{Algorithm_BCD} and Algorithm \ref{Algorithm_PG}, respectively. In Algorithm \ref{Algorithm_PG}, $d^{(k)}$ in Step 7 represents the final stepsize of iteration $k$ and is obtained by the non-monotone line search with $d_0^{(k)}$ as the initial value \cite{birgin2000nonmonotone}, implying that $d^{(k)} \leq d_0^{(k)}$ \cite{birgin2000nonmonotone}. \red{For fair comparison, in Algorithm \ref{Algorithm_PG}, we omit the active set selection \cite{wang2021efficient} and the non-monotone line search \cite{birgin2000nonmonotone}.}
(i) PSCA-ML-K and PG-ML-K are parallel algorithms, whereas BCD-ML-K is a sequential algorithm yielding a long computation time. 
(ii) In each iteration, the update of $\alpha_n^{(k)}, n \in \mathcal{N}$ in \red{Step 7} of PSCA-ML-K and that in Step 6 of BCD-ML-K are identical and approximately utilize up to the second-order information of $f_{\textrm{\rm ML-K}}\left(\boldsymbol \alpha\right)$ via $\left(g_n \mathbf{s}_n^H \boldsymbol{\Sigma}_{\boldsymbol{\alpha}^{(k)}}^{-1} \mathbf{s}_n\right)^2$, whereas that in \red{Step 9} of PG-ML-K only uses up to the first-order information of $f_{\textrm{\rm ML-K}}\left(\boldsymbol \alpha\right)$, yielding unsatisfactory convergence speed. Specifically, the coefficient of $\partial_{n} f_{\textrm{\rm ML-K}}\left(\boldsymbol{\alpha}^{(k)}\right)$ in \red{Step 7 of PSCA-ML-K} and that in Step 6 of BCD-ML-K is $\frac{1}{\left(g_n \mathbf{s}_n^H \boldsymbol{\Sigma}_{\boldsymbol{\alpha}^{(k)}}^{-1} \mathbf{s}_n\right)^2}$, whereas the coefficient of $\partial_{n} f_{\textrm{\rm ML-K}}\left(\boldsymbol{\alpha}^{(k)}\right)$ in \red{Step 9} of PG-ML-K is $d^{(k)}$.
(iii) In each iteration of PSCA-ML-K and PG-ML-K,  $\boldsymbol{\Sigma}_{\boldsymbol{\alpha}^{(k)}}^{-1}$ is updated only once, using Frobenius inversion according to (\ref{Frobenius_inversion_1}) and (\ref{Frobenius_inversion_2}), with high computational complexity. In contrast, in each iteration of BCD-ML-K, $\boldsymbol{\Sigma}_{\boldsymbol{\alpha}^{(k)}}^{-1}$ is updated $N$ times, each using the Woodbury matrix inverse lemma in Step 7 \red{of BCD-ML-K}, with low computational complexity. Overall, the update of $\boldsymbol{\Sigma}_{\boldsymbol{\alpha}^{(k)}}^{-1}$  in each iteration of PSCA-ML-K and PG-ML-K achieves a lower computational complexity than the $N$ updates of $\boldsymbol{\Sigma}_{\boldsymbol{\alpha}^{(k)}}^{-1}$  in each iteration of BCD-ML-K when $N \gg L$.
(iv) Regarding the information of $f_{\textrm{\rm ML-K}\left(\boldsymbol\alpha\right)}$, PSCA-ML-K and BCD-ML-K utilize up to second-order information, whereas PG-ML-K utilizes up to  first-order information. Therefore, PSCA-ML-K is expected to achieve a shorter computation time than BCD-ML-K and PG-ML-K. 

\begin{algorithm}[t]
	\caption{\red{PSCA Algorithms for Device Activity Detection: PSCA-ML-K, PSCA-MAP-K, PSCA-ML-UD, and PSCA-MAP-UR.}}
	\label{Algorithm_PSCA_General}
	\begin{algorithmic}[1]
		\STATE \textbf{Input:} $\widehat{\boldsymbol{\Sigma}}_{\mathbf{Y}}$.
		\STATE \textbf{Output:}\\
        \ \ Known $\mathbf{g}$: $\boldsymbol{\alpha}$.\\
        \ \ Unknown $\mathbf{g}$: $\boldsymbol{\gamma}$.
		\STATE \textbf{Initialize:}\\
        \ \ Known $\mathbf{g}$: $\boldsymbol{\alpha}^{(0)} = \mathbf{0}_N$ and $k = 0$.\\
        \ \ Unknown $\mathbf{g}$: $\boldsymbol{\gamma}^{(0)} = \mathbf{0}_N$ and $k = 0$.
		\REPEAT
        \STATE Known $\mathbf{g}$: Calculate $\boldsymbol{\Sigma}_{\boldsymbol \alpha^{(k)}}$ according to (\ref{cov_a}).\\
        Unknown $\mathbf{g}$: Calculate $\boldsymbol{\Sigma}_{\boldsymbol \gamma^{(k)}}$ according to (\ref{cov_b}).
        \STATE Known $\mathbf{g}$: Calculate $\boldsymbol{\Sigma}_{\boldsymbol \alpha^{(k)}}^{-1}$ according to (\ref{Frobenius_inversion_1}) and (\ref{Frobenius_inversion_2}).\\
        Unknown $\mathbf{g}$: Calculate $\boldsymbol{\Sigma}_{\boldsymbol \gamma^{(k)}}^{-1}$ according to (\ref{Frobenius_inversion_1}) and (\ref{Frobenius_inversion_2}).
        

		\STATE 
        Known $\mathbf{g}$: Calculate $\widetilde{\alpha}_{\textrm{\rm ML-K}, n}^{(k)}, n \in \mathcal{N}$  in parallel according to (\ref{Solution_ML}) for MLE or $\widetilde{\alpha}_{\textrm{\rm MAP-K}, n}^{(k)}, n \in \mathcal{N}$ in parallel according to (\ref{Solution_MAP}) for MAPE.\\
        Unknown $\mathbf{g}$: Calculate $\widetilde{\gamma}_{\textrm{\rm ML-UD}, n}^{(k)}, n \in \mathcal{N}$  in parallel according to (\ref{Solution_ML_U}) for MLE or $\widetilde{\gamma}_{\textrm{\rm MAP-UR}, n}^{(k)}, n \in \mathcal{N}$ in parallel according to (\ref{Solution_MAP_U}) for MAPE.
        
		\STATE Known $\mathbf{g}$: Calculate $\alpha_n^{(k)}, n \in \mathcal{N}$ in parallel according to (\ref{Update_ML}) for MLE or (\ref{Update_MAP}) for MAPE.\\
        Unknown $\mathbf{g}$: Calculate $\gamma_n^{(k)}, n \in \mathcal{N}$ in parallel according to (\ref{Update_ML_U}) for MLE or (\ref{Update_MAP_U}) for MAPE.
		\STATE Set $k = k+1$.\\
		\UNTIL $\boldsymbol\alpha^{(k)}$ or $\boldsymbol\gamma^{(k)}$ satisfies some stopping criterion.  
	\end{algorithmic}
\end{algorithm}

\begin{algorithm}[t]
	\caption{BCD-ML-K \cite{fengler2021nonbayesian}}
	\label{Algorithm_BCD}
	\begin{algorithmic}[1]
		\STATE \textbf{Input: }$\widehat{\boldsymbol{\Sigma}}_{\mathbf{Y}}$. 
		\STATE \textbf{Output: }$\boldsymbol\alpha$. 
		\STATE \textbf{Initialize: }$\boldsymbol\alpha^{(0)} = \mathbf{0}, \boldsymbol{\Sigma}_{\boldsymbol{\alpha}^{(0)}}^{-1} = \frac{1}{\sigma^2}\mathbf{I}_L, \textrm{and } k = 0$.
		\REPEAT
		\FOR{$n \in \mathcal{N}$}
		\STATE Calculate $\alpha_n^{(k+1)}$ according to (\ref{Solution_ML}) with $ \widetilde{\alpha}_n^{(k+1)}$ replaced by $\alpha_n^{(k+1)}$.\\
		\STATE Calculate $\boldsymbol{\Sigma}_{\boldsymbol \alpha^{(k)}}^{-1} = \boldsymbol{\Sigma}_{\boldsymbol \alpha^{(k)}}^{-1} - \frac{(\alpha_n^{(k+1)}-\alpha_n^{(k)})g_n\boldsymbol{\Sigma}_{\boldsymbol \alpha^{(k)}}^{-1}\mathbf{s}_n\mathbf{s}_n^H\boldsymbol{\Sigma}_{\boldsymbol \alpha^{(k)}}^{-1}}{1+(\alpha_n^{(k+1)}-\alpha_n^{(k)})\left(g_n \mathbf{s}_n^H \boldsymbol{\Sigma}_{\boldsymbol{\alpha}^{(k)}}^{-1} \mathbf{s}_n\right)}$.\\
		\ENDFOR
		\STATE Set $k = k+1$.
		\UNTIL $\boldsymbol\alpha^{(k)}$ satisfies some stopping criterion. 
	\end{algorithmic}
\end{algorithm}
\begin{algorithm}[t]
	\caption{PG-ML-K (without active set selection) \cite{wang2021efficient}}
	\label{Algorithm_PG}
	\begin{algorithmic}[1]
		\STATE \textbf{Input: }$\widehat{\boldsymbol{\Sigma}}_{\mathbf{Y}}$. 
		\STATE \textbf{Output: }$\boldsymbol\alpha$. 
		\STATE \textbf{Initialize: }$\boldsymbol\alpha^{(0)} = \mathbf{0}, \boldsymbol{\Sigma}_{\boldsymbol{\alpha}^{(0)}}^{-1} = \frac{1}{\sigma^2}\mathbf{I}_L, \textrm{and } k = 1$.
		\REPEAT
  \STATE Calculate $\boldsymbol{\Sigma}_{\boldsymbol \alpha^{(k)}}$ according to (\ref{cov_a}).
  \STATE Calculate $\boldsymbol{\Sigma}_{\boldsymbol \alpha^{(k)}}^{-1}$ according to (\ref{Frobenius_inversion_1}) and (\ref{Frobenius_inversion_2}).
  \STATE Calculate $\partial_{n} f_{\textrm{\rm ML-K}}\left(\boldsymbol{\alpha}^{(k)}\right)$ in parallel according to (\ref{gradient_ML})
  \STATE Calculate \\ 
  \small{\small{$\alpha_n^{(k+1)} =  \min \left\{\max \left\{\alpha_n^{(k)} - d^{(k)} \partial_{\alpha_n}f_{\textrm{\rm ML-}K}\left(\boldsymbol{\alpha}_n^{(k)}\right),0\right \},1\right \}$}},\\
   $n \in \mathcal{N}$ in parallel.
		\STATE Set $k = k+1$.\\
		\UNTIL $\boldsymbol\alpha^{(k)}$ satisfies some stopping criterion.   
	\end{algorithmic}
\end{algorithm}

\begin{table}[!ht]
	\caption{Per-Iteration Computational Complexities of PSCA-ML-K, BCD-ML-K \cite{fengler2021nonbayesian},  PG-ML-K \cite{wang2021efficient}, PSCA-MAP-K, \red{and BCD-MAP-K \cite{jiang2022ml}}.} 
 \label{table_known}
	\centering
	\begin{tabular}{|l|l|l|}
		\hline
		\textbf{Algorithms} & \textbf{Dominant Term} & \textbf{Order}\\ \hline
  \textbf{PSCA-ML-K} & \red{$40NL^2$} & $\mathcal{O}(NL^2)$\\ \hline
  \textbf{BCD-ML-K} & $56NL^2$ & $\mathcal{O}(NL^2)$\\ \hline 
  \textbf{PG-ML-K} & \red{$40NL^2$} & $\mathcal{O}(NL^2)$\\ \hline
  \textbf{PSCA-MAP-K (general)} & \red{$40NL^2+2^N$} & $\mathcal{O}(NL^2+2^N)$\\ \hline
    \textbf{PSCA-MAP-K (independent)} & \red{$40NL^2$} & $\mathcal{O}(NL^2)$\\ \hline
      \textbf{BCD-MAP-K (general)} & $56NL^2+2^N$ & $\mathcal{O}(NL^2+2^N)$\\ \hline 
  \textbf{BCD-MAP-K (independent)} & $56NL^2$ & $\mathcal{O}(NL^2)$\\ \hline 
	\end{tabular}
\end{table}

Finally, we analyze the computational complexities of PSCA-ML-K, BCD-ML-K, and PG-ML-K by characterizing the dominant terms and orders of the per-iteration flop counts under the assumption that $N \gg L$.
The results are summarized in Table \ref{table_known}, \red{and the proof is given in Appendix~\ref{appendix:computational}.} According to Table \ref{table_known}, PSCA-ML-K has the lowest per-iteration computational complexity among the three MLE-based device activity detection methods.\footnote{\red{Compared to PG-ML-K \cite{wang2021efficient} with the active set selection and the
non-monotone line search, PG-ML-K has a low per-iteration computational complexity in the dominant term and the same per-iteration computational complexity in order.}} 

\begin{figure*}[!ht]
	\begin{align}
		\label{Loss_alpha}
		\red{\textrm{Loss}\left( \left\{ \widetilde{\boldsymbol{\alpha}}^{[i]} \right\}_{i \in \mathcal{I}_{P}}, \left\{ \boldsymbol{\alpha}^{[i]} \right\}_{i \in \mathcal{I}_{P}}\right) = - \frac{1}{NI_{P}}\sum_{i \in \mathcal{I}_{P}}\sum_{n \in \mathcal{N}}\left( \alpha_n^{[i]} \log \left(\widetilde{\alpha}_n^{[i]} \right)+\left(1-\alpha_n^{[i]}\right) \log \left(1-\widetilde{\alpha}_n^{[i]}\right)\right).}
	\end{align}
	\hrulefill
\end{figure*}

\subsubsection{PSCA-ML-K-Net}
\label{PSCA_ML_Net}
\red{To reduce the overall computation time of PSCA-ML-K, it remains to optimize its step size sequence, which plays a crucial role in speeding up the convergence.} In this part, we propose a PSCA-ML-K-driven deep unrolling neural network approach, called PSCA-ML-K-Net, \red{to optimize PSCA-ML-K’s step size}. 
\begin{figure}[tp]
	\begin{center}
		{\resizebox{8cm}{!}{\includegraphics{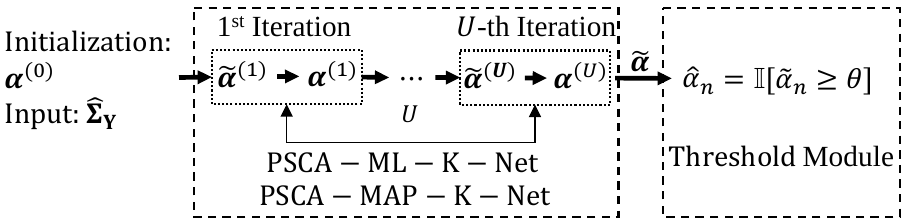}}}
	\end{center}
	\caption{Architectures of PSCA-ML-K-Net and PSCA-MAP-K-Net.}
 \label{NN}
\end{figure}

As illustrated in Fig. \ref{NN}, the neural network has $U$~($> 0$) building blocks, each implementing one iteration of PSCA-ML-K. The operations for complex numbers in PSCA-ML-K are readily implemented via the operations for real numbers using a standard neural network \cite{cui2021jointly}.\footnote{One building block of each PSCA-Net and one iteration of the underlying PSCA algorithm share the same computational complexity (in exact form and thus dominant term and order). Each PSCA-Net does not rely on any basic neural network module.} The step size of iteration $k$, i.e., $\rho^{(k)}$, is tunable and optimized by the neural network training \red{on the training samples to improve the convergence speed}. The choice of $U$ determines PSCA-ML-K's tradeoff between the error rate and computation time. Specifically, the accuracy increases at a decreasing rate with $U$, and the computation time increases linearly with $U$.

Now, we introduce the training, validation, and test procedures of PSCA-ML-K-Net. Consider $I_P$ training, $I_Q$ validation, and $I_R$ test data samples, i.e., \red{$\left\{\left( \mathbf{S}^{[i]}, \mathbf{g}^{[i]}, \widehat{\boldsymbol{\Sigma}}_{\mathbf{Y}}^{[i]}, \boldsymbol{\alpha}^{[i]} \right)\right\}_{i \in \mathcal{I}_j}$}, where $\mathcal{I}_j \triangleq \{1,2,...,I_j\}, j \in \{P,Q,R\}$. Let $\widetilde{\boldsymbol{\alpha}}^{[i]}$ represent the output \red{that corresponds} to the input $\left(\mathbf{S}^{[i]}, \mathbf{g}^{[i]}, \widehat{\boldsymbol{\Sigma}}_{\mathbf{Y}}^{[i]}\right)$. The binary cross-entropy (BCE) loss function given by (\ref{Loss_alpha}), as shown at the top of the next page, is used to measure the distance between $\boldsymbol{\alpha}^{[i]}$ and $\widetilde{\boldsymbol{\alpha}}^{[i]}, i \in \mathcal{I}_{P}$ \cite{cui2021jointly}. We train PSCA-ML-K-Net in an end-to-end manner using the adaptive moment estimation (ADAM) optimizer. \red{When the value of the loss function on the validation samples does not change, the training process is stopped, and the step size for PSCA-ML-K-Net is saved.} Finally, we convert $\widetilde{\boldsymbol{\alpha}}^{[i]}$ to the final output $\widehat{\boldsymbol{\alpha}}^{[i]}$ by a threshold module, i.e., $\widehat{\alpha}_n^{[i]} = \mathbb{I}\left[\widetilde{\alpha}_n^{[i]} \geq \theta\right],n\in\mathcal{N},i\in\mathcal{I}_{Q}\cup\mathcal{I}_{R}$, where $\theta\in (0,1)$ is a threshold that minimizes
the average error rate over the validation data $\frac{1}{NI_{Q}}\sum_{i \in \mathcal{I}_{Q}}\left\Vert \widehat{\boldsymbol{\alpha}}^{[i]} - \boldsymbol{\alpha}^{[i]}\right\Vert_1$ \cite{cui2021jointly}.
\subsection{MAPE-Based Device Activity Detection for Known Pathloss}
\label{MAPEK}
In this subsection, we view $\boldsymbol{\alpha}$ as realizations of Bernoulli random variables and study the MAPE of $\boldsymbol{\alpha}$ for known $\mathbf{g}$.
\subsubsection{MAPE Problem for Known Pathloss}
\label{MAPE_Problem}
In this part, we present the MAPE problem of $\boldsymbol{\alpha}$ for known $\mathbf{g}$.
Specifically, we adopt the general model \red{(or its approximation)} and the independent model for the p.m.f. of $\mathbf{A}$, i.e., the MVB model $p_{\mathbf{A}\textrm{\rm-G}}\left(\boldsymbol{\alpha}\right)$ in (\ref{MVB}) \cite{jiang2022ml,jiang2020mapbased} \red{(or (\ref{second}))} and the independent model in (\ref{independent}). The conditional p.m.f of $\mathbf{A}$ given $\mathbf{Y}$, denoted by $p_{\textrm{\rm MAP-K}}(\mathbf{Y},\boldsymbol\alpha)$, is expressed as:  
\begin{align}
  p_{\textrm{\rm MAP-K}}(\mathbf{Y},\boldsymbol\alpha) = p_{\textrm{\rm ML-K}}(\mathbf{Y};\boldsymbol\alpha)p_{\mathbf{A}\textrm{\rm-t}}\left(\boldsymbol{\alpha}\right),  
\end{align}
where $p_{\textrm{\rm ML-K}}(\mathbf{Y};\boldsymbol\alpha)$ is \red{given by (\ref{fY}), and $p_{\mathbf{A}\textrm{\rm-t}}\left(\boldsymbol{\alpha}\right), \textrm{t} \in \{\textrm{G}, \textrm{I}\}$} are given by (\ref{MVB}) and (\ref{independent}). For notation convenience, we do not differentiate $p_{\textrm{\rm MAP-K}}(\mathbf{Y},\boldsymbol\alpha)$ for the two distribution models for $\mathbf{A}$.
Let 
\begin{align}
\label{likelihood_Y_MAP}
	f_{\textrm{\rm MAP-K}}\left(\boldsymbol \alpha\right) & \triangleq -\frac{1}{M} \log p_{\textrm{\rm MAP-K}}(\mathbf{Y},\boldsymbol\alpha) - L \log{\pi} \nonumber \\
	& \red{\quad -\frac{1}{M}\log
    p_{\mathbf{A}\textrm{\rm-t}}\left(\mathbf{0}_N\right)} \nonumber \\
	& =  f_{\textrm{\rm ML-K}}\left(\boldsymbol \alpha\right) + f_{\textrm{t}}\left(\boldsymbol \alpha\right), \ \textrm{t} \in \{\textrm{G},\textrm{I}\}
\end{align}  
where \red{$f_{\textrm{\rm ML-K}}\left(\boldsymbol \alpha\right)$ is given by (\ref{likelihood_Y}), and $f_{\textrm{t}}\left(\boldsymbol \alpha\right)$} is given by:
\begin{align}
\label{ft}
    f_{\textrm{t}}\left(\boldsymbol \alpha\right) \triangleq 
    \begin{cases}
        -\frac{1}{M}\sum_{\omega\in\Psi}c_{\omega}\prod_{n\in\omega}\alpha_n, & \textrm{t} = \textrm{G}, \\
        -\frac{1}{M}\sum_{n \in \mathcal{N}}\log\left(\frac{p_n}{1-p_n}\right)\alpha_n, & \textrm{t} = \textrm{I}.
    \end{cases}
\end{align}
Note that $f_{\textrm{t}}\left(\boldsymbol \alpha\right)$ is from $p_{\mathbf{A}\textrm{\rm-t}}\left(\boldsymbol{\alpha}\right)$. Thus, the MAP problem of $\boldsymbol\alpha$, i.e., the maximization of $p_{\textrm{\rm MAP-K}}(\mathbf{Y},\boldsymbol\alpha)$ or $\log\left(p_{\textrm{\rm MAP-K}}(\mathbf{Y},\boldsymbol\alpha)\right)$ w.r.t. $\boldsymbol \alpha$, can be equivalently formulated as follows.
\begin{Prob}[MAPE of $\boldsymbol\alpha$ for Known $\mathbf{g}$]
	\label{Prob_MAP}
	\begin{align*}
		\min_{\boldsymbol\alpha} & \quad f_{\textrm{\rm MAP-K}}\left(\boldsymbol \alpha\right) \nonumber \\
		\text {s.t.} & \quad (\ref{constraint}).
	\end{align*}
\end{Prob}

Compared to Problem~\ref{Prob_ML}, Problem~\ref{Prob_MAP} additionally incorporates the prior distribution $p_{\mathbf{A}\textrm{\rm-t}}\left(\boldsymbol{\alpha}\right)$ in (\ref{independent}) or (\ref{MVB}) in the objective function. Since $f_{\textrm{\rm MAP-K}} - f_{\textrm{\rm ML-K}} = f_{\textrm{t}}(\boldsymbol{\alpha}) \to 0$, as $M \to \infty$, \red{the influence of} $p_{\mathbf{A}\textrm{\rm-t}}\left(\boldsymbol{\alpha}\right)$  becomes less critical, and the observation $\mathbf Y$ becomes more dominating, as $M$ increases \cite{jiang2022ml}.
Problem~\ref{Prob_MAP} is a more challenging non-convex problem than Problem~\ref{Prob_ML}, due to the additional polynomial term $f_{\textrm{\rm-t}}\left(\boldsymbol \alpha\right)$ \red{in (\ref{ft})}. Like BCD-ML-K, the existing state-of-the-art BCD-MAP-K \cite{jiang2020mapbased,jiang2023statistical,cui2021jointly} for MAPE-based device activity detection can guarantee convergence to \red{Problem~\ref{Prob_MAP}'s stationary points} but yields a longer computation time due to the sequential update mechanism.
\subsubsection{PSCA-MAP-K}
\label{PSCA_MAP_Algorithm}
In this part, we propose an algorithm, referred to as PSCA-MAP-K, using PSCA. Specifically, at iteration $k$, we choose:
\begin{align}
	\label{Surrogate_MAP}
	\widetilde{f}_{\textrm{\rm MAP-K}}\left(\boldsymbol{\alpha};\boldsymbol{\alpha}^{(k)}\right) \triangleq \sum_{n \in \mathcal{N}} & \widetilde{f}_{\textrm{\rm MAP-K}, n}\left(\alpha_n ;\boldsymbol{\alpha}^{(k)}\right) 
\end{align} as a strongly convex approximate function of $f_{\textrm{\rm MAP-K}}\left( \boldsymbol{\alpha}\right)$ around $\boldsymbol{\alpha}^{(k)}$ w.r.t. $\boldsymbol{\alpha}$, where
\begin{align}
	\label{Surrogate_MAP_separate}
	\widetilde{f}_{\textrm{\rm MAP-K}, n}(\alpha_n;\boldsymbol{\alpha}^{(k)})  \triangleq & \partial_{n} f_{\textrm{\rm MAP-K}}\left(\boldsymbol{\alpha}^{(k)}\right) \left(\alpha_n - \alpha_n^{(k)}\right) \nonumber \\  + & \frac{1}{2}\left(g_n \mathbf{s}_n^H \boldsymbol{\Sigma}_{\boldsymbol{\alpha}^{(k)}}^{-1} \mathbf{s}_n\right)^2 \left(\alpha_n - \alpha_n^{(k)}\right)^2
\end{align} is a strongly convex approximate function of \red{$f_{\textrm{\rm MAP-K}}\left(\boldsymbol{\alpha}\right)$} around $\boldsymbol{\alpha}^{(k)}$ w.r.t. $\alpha_n$. Here,
\begin{align}
	\label{gradient_MAP}
	\partial_{n} f_{\textrm{\rm MAP-K}}(\boldsymbol{\alpha}^{(k)}) = & \partial_{n} f_{\textrm{\rm ML-K}}\left(\boldsymbol{\alpha}^{(k)}\right) + \partial_{n}f_{\textrm{t}}(\boldsymbol{\alpha}^{(k)}),
\end{align} 
where
\begin{align}
\label{ft_gradient}
    \red{\partial_{n}f_{\textrm{t}}(\boldsymbol{\alpha}^{(k)})}  = 
    \begin{cases}
        -\frac{1}{M}\sum_{\omega\in\Psi}\left(c_{\omega}\prod_{n'\in\omega \backslash {\{n\}}}\alpha_{n'}^{(k)}\right), & \textrm{t} = \textrm{G}, \\
        -\log\left(\frac{p_n}{1-p_n}\right), & \textrm{t} = \textrm{I}.
    \end{cases}
\end{align}
For all $n \in \mathcal{N}$, by noting that $\widetilde{f}_{\textrm{\rm MAP-K}, n}(\alpha_n;\boldsymbol{\alpha}^{(k)})$ in (\ref{Surrogate_MAP_separate}) and $\widetilde{f}_{\textrm{\rm ML-K}, n}(\alpha_n;\boldsymbol{\alpha}^{(k)})$ in (\ref{Surrogate_ML_separate}) share the same quadratic coefficient, i.e., $\left(g_n \mathbf{s}_n^H \boldsymbol{\Sigma}_{\boldsymbol{\alpha}^{(k)}}^{-1} \mathbf{s}_n\right)^2$, \red{and the derivative of the additional term $\partial_{n}f_{\textrm{t}}(\boldsymbol{\alpha}^{(k)})$ w.r.t. $\alpha_n$ is 0}, we can conclude that $\left(g_n \mathbf{s}_n^H \boldsymbol{\Sigma}_{\boldsymbol{\alpha}^{(k)}}^{-1} \mathbf{s}_n\right)^2$ also converges to the $n$-th diagonal element of $\nabla^2 f_{\textrm{\rm MAP-K}} (\boldsymbol{\alpha}^{(k)})$, almost surely, as $M \to \infty$. 
Thus, the surrogate function $\widetilde{f}_{\textrm{\rm MAP-K}}\left(\boldsymbol{\alpha};\boldsymbol{\alpha}^{(k)}\right)$ given in (\ref{Surrogate_MAP}) approximately utilizes up to the second-order information of the objective function $f_{\textrm{\rm MAP-K}}\left(\boldsymbol \alpha\right)$ of Problem~\ref{Prob_MAP}.
Similarly, the approximate convex problem at iteration $k$ is given by:
\begin{align}
	\underset{\boldsymbol{\alpha} \in [0,1]^N}{\operatorname{min}}  \widetilde{f}_{\textrm{\rm MAP-K}}\left(\boldsymbol{\alpha};\boldsymbol{\alpha}^{(k)}\right),
\end{align}
which can be equivalently separated into $N$ convex problems, one for each coordinate.
\begin{Prob}[Convex Approximate Problem of Problem~\ref{Prob_MAP} for $n \in\mathcal{N}$ at Iteration $k$]
	\label{Prob_Surrogate_MAP}
	\begin{align*}
	\red{\widetilde{\alpha}_{\textrm{\rm MAP-K},n}^{(k)} \triangleq \underset{\alpha_n \in [0,1]}{\operatorname{argmin}} \  \widetilde{f}_{\textrm{\rm MAP-K}}\left(\alpha_n;\boldsymbol{\alpha}^{(k)}\right).}
    \end{align*}
\end{Prob}
\begin{Lem}[Optimal Solution of Problem~\ref{Prob_Surrogate_MAP}]
	\label{Lem_MAP}
	The optimal solution of Problem~\ref{Prob_Surrogate_MAP} is given by:
	\begin{align}
		\label{Solution_MAP}
        &\widetilde{\alpha}_{\textrm{\rm MAP-K},n}^{(k)} \nonumber \\ &= \min \left\{\max \left\{- \frac{\partial_{n} f_{\textrm{\rm MAP-K}}(\boldsymbol{\alpha}^{(k)})}{\left(g_n \mathbf{s}_n^H \boldsymbol{\Sigma}_{\boldsymbol{\alpha}^{(k)}}^{-1} \mathbf{s}_n \right)^2} + \alpha_n^{(k)} ,0 \right \},1 \right \},
	\end{align}
    \red{where $\partial_{n} f_{\textrm{\rm MAP-K}}\left(\boldsymbol{\alpha}^{(k)}\right)$ is given by (\ref{gradient_MAP}).}
\end{Lem}
\begin{IEEEproof}
Please refer to Appendix~\ref{appendix:solution}.
\end{IEEEproof}

Note that $\widetilde{\alpha}_{\textrm{\rm MAP-K},n}^{(k)}$ in (\ref{Solution_MAP}) shares a form similar to $\widetilde{\alpha}_{\textrm{\rm ML-K},n}^{(k)}$ in (\ref{Solution_ML}) \red{and approximately} utilizes the second-order information of $f_{\textrm{\rm MAP-K}}\left(\boldsymbol \alpha\right)$ via $\left(g_n \mathbf{s}_n^H \boldsymbol{\Sigma}_{\boldsymbol{\alpha}^{(k)}}^{-1} \mathbf{s}_n \right)^2$. By (\ref{gradient_MAP}) and $c_n^{(k)} \to 0$, we have $\widetilde{\alpha}_{\textrm{\rm MAP-K},n}^{(k)} \to \widetilde{\alpha}_{\textrm{\rm ML-K},n}^{(k)}$, as $M \to \infty$.
Thus, the influence of the prior distribution of $\boldsymbol{\alpha}$ decreases as $M$ increases \cite{jiang2022ml}.
Furthermore, PSCA-MAP-K analytically solves $N$ convex approximate problems (Problem \ref{Prob_Surrogate_MAP}) in parallel, whereas BCD-MAP-K analytically and sequentially solves $N$ non-convex coordinate descent optimization problems. In addition, the solution of each convex approximate problem \red{in PSCA-MAP-K} is much simpler than that of the corresponding non-convex coordinate descent problem in BCD-MAP-K (see Steps 8$\sim$13).

Finally, we update $\alpha_n^{(k)}, n \in \mathcal{N}$ in parallel according to:
\begin{align}
	\label{Update_MAP}
	\alpha_n^{(k+1)} = \left(1 - \rho^{(k)}\right)\alpha_n^{(k)} + \rho^{(k)} \widetilde{\alpha}_{\textrm{\rm MAP-K},n}^{(k)}, \  n \in \mathcal{N},
\end{align}
where $\rho^{(k)}$ is the step size satisfying (\ref{diminishing_rule}).
The detailed procedure is summarized in Algorithm \ref{Algorithm_PSCA_General} (PSCA-MAP-K).
\begin{Thm}[Convergence of PSCA-MAP-K]
	\label{Convergence_MAP}
	Every limit point of the sequence $\left\{\boldsymbol{\alpha}^{(k)}\right\}_{k \in \mathbb{Z}_{+}}$ generated by PSCA-MAP-K (at least one such point exists) is a stationary point of Problem~\ref{Prob_MAP}.
\end{Thm}
\begin{IEEEproof}
	Please refer to Appendix~\ref{appendix:convergence}.
\end{IEEEproof}

PSCA-MAP-K and PSCA-ML-K are almost the same except for the slight difference in \red{Step 7}.
In the following, we compare PSCA-MAP-K and BCD-MAP-K \cite{jiang2022ml}.
    (i) PSCA-MAP-K is a parallel algorithm, whereas BCD-MAP-K is a sequential algorithm.
    (ii) The update of $\alpha_n^{(k)}, n \in \mathcal{N}$ in each iteration of PSCA-MAP-K \red{(see Steps 7$\sim$8)} is much simpler than those in each iteration of BCD-MAP-K (see Steps 8$\sim$14).\footnote{Since the optimal point of the coordinate optimization problem \red{for MAPE} in \cite{jiang2022ml} is discontinuous, there does not exist a continuously differentiable surrogate function for the approximate problem \red{for MAPE} that has the same optimal point.}
    (iii) The single update of $\boldsymbol{\Sigma}_{\boldsymbol{\alpha}^{(k)}}^{-1}$ in each iteration of PSCA-MAP-K achieves a lower computational complexity than the $N$ updates of $\boldsymbol{\Sigma}_{\boldsymbol{\alpha}^{(k)}}^{-1}$ in each iteration of BCD-MAP-K when $N \gg L$.
    (iv) Regarding the information of $f_{\textrm{\rm MAP-K}\left(\boldsymbol\alpha\right)}$, PSCA-MAP-K and BCD-MAP-K both approximately utilize up to the second-order information via $\left(g_n \mathbf{s}_n^H \boldsymbol{\Sigma}_{\boldsymbol{\alpha}^{(k)}}^{-1} \mathbf{s}_n \right)^2$.
Based on the above observation, PSCA-MAP-K is expected to achieve a shorter computation time than BCD-MAP-K.

Finally, we analyze the dominant terms and orders of the per-iteration flop counts of PSCA-ML-K and BCD-MAP-K under the assumption that $N \gg L$ to compare their per-iteration computational complexities.
The results are summarized in Table \ref{table_known}, and the proof is given in Appendix~\ref{appendix:computational}. According to Table \ref{table_known}, PSCA-MAP-K has a lower per-iteration computational complexity than BCD-MAP-K and a higher computational complexity than PSCA-ML-K.

\begin{algorithm}[t]
	\caption{BCD-MAP-K \cite{jiang2022ml}}
 	\label{Algorithm_BCD_MAP}
		\begin{algorithmic}[1]
			\STATE \textbf{Input: }empirical covariance matrix $\widehat{\boldsymbol{\Sigma}}_{\mathbf{Y}}$. 
			\STATE \textbf{Output: }activities of devices $\boldsymbol\alpha$. 
			\STATE \textbf{Initialize: }$\boldsymbol\alpha^{(0)} = \mathbf{0}, \boldsymbol{\Sigma}_{\boldsymbol{\alpha}^{(0)}}^{-1} = \frac{1}{\sigma^2}\mathbf{I}_L, k = 0$.
			\REPEAT
   \FOR{$n \in \mathcal{N}$}
			\STATE Calculate $f'_{t}\left(\boldsymbol \alpha^{(k)}\right)$ in (\ref{gradient_MAP}).
   \IF{$f'_{t}\left(\boldsymbol \alpha^{(k)}\right) \leq 0$}
   \STATE Calculate $\alpha_n^{(k+1)}$ according to (\ref{BCDMAP}) with $\overline{\alpha}_n^{(k+1)}$ replaced by $\alpha_n^{(k+1)}$.\
   \ELSIF{$0<f'_{t}\left(\boldsymbol \alpha^{(k)}\right)<\frac{g_n\left(\mathbf{s}_n^H \boldsymbol{\Sigma}_{\boldsymbol{\alpha}^{(k)}}^{-1} \mathbf{s}_n\right)^2}{4\mathbf{s}_n^H \boldsymbol{\Sigma}_{\boldsymbol{\alpha}^{(k)}}^{-1} \widehat{\boldsymbol{\Sigma}}_{\mathbf{Y}} \boldsymbol{\Sigma}_{\boldsymbol{\alpha}^{(k)}}^{-1} \mathbf{s}_n}$}
   \STATE Solve $\alpha_n^{(k+1)} = \underset{\alpha_n  \in \{\overline{\alpha}_n^{(k+1)},0,1\}} {\operatorname{argmin}}f_{\textrm{\rm MAP-}K}\left(\alpha_n;\boldsymbol \alpha^{(k)}\right)$.
   \ELSIF{$f'_{t}\left(\boldsymbol \alpha^{(k)}\right) \geq \frac{g_n\left(\mathbf{s}_n^H \boldsymbol{\Sigma}_{\boldsymbol{\alpha}^{(k)}}^{-1} \mathbf{s}_n\right)^2}{4\mathbf{s}_n^H \boldsymbol{\Sigma}_{\boldsymbol{\alpha}^{(k)}}^{-1} \widehat{\boldsymbol{\Sigma}}_{\mathbf{Y}} \boldsymbol{\Sigma}_{\boldsymbol{\alpha}^{(k)}}^{-1} \mathbf{s}_n}$.}
   \STATE Set $\alpha_n^{(k+1)} = 1$.
   \ENDIF
   \STATE Calculate $\boldsymbol{\Sigma}_{\boldsymbol \alpha^{(k)}}^{-1} = \boldsymbol{\Sigma}_{\boldsymbol \alpha^{(k)}}^{-1} - \frac{(\alpha_n^{(k+1)}-\alpha_n^{(k)})g_n\boldsymbol{\Sigma}_{\boldsymbol \alpha^{(k)}}^{-1}\mathbf{s}_n\mathbf{s}_n^H\boldsymbol{\Sigma}_{\boldsymbol \alpha^{(k)}}^{-1}}{1+(\alpha_n^{(k+1)}-\alpha_n^{(k)})\left(g_n \mathbf{s}_n^H \boldsymbol{\Sigma}_{\boldsymbol{\alpha}^{(k)}}^{-1} \mathbf{s}_n\right)}$.\\
\ENDFOR
			\STATE Set $k = k+1$.\\
			\UNTIL $\boldsymbol\alpha^{(k)}$ satisfies some stopping criterion.
			
		\end{algorithmic}
\end{algorithm}
\begin{figure*}[!hbp]
	\hrulefill
	\begin{align}
		\label{BCDMAP}
		\overline{\alpha}_n^{(k+1)} = \min \left\{\max \left\{\alpha_n^{(k)} + \frac{1}{2f'_{\textrm{\rm t}}\left(\boldsymbol \alpha\right)}\left(1- \sqrt{1-\frac{4f'_{\textrm{\rm t}}\left(\boldsymbol \alpha^{(k)}\right)\mathbf{s}_n^H \boldsymbol{\Sigma}_{\boldsymbol{\alpha}^{(k)}}^{-1} \widehat{\boldsymbol{\Sigma}}_{\mathbf{Y}} \boldsymbol{\Sigma}_{\boldsymbol{\alpha}^{(k)}}^{-1} \mathbf{s}_n}{g_n\left(\mathbf{s}_n^H \boldsymbol{\Sigma}_{\boldsymbol{\alpha}^{(k)}}^{-1} \mathbf{s}_n\right)^2}}\right)-\frac{1}{g_n\mathbf{s}_n^H \boldsymbol{\Sigma}_{\boldsymbol{\alpha}^{(k)}}^{-1} \mathbf{s}_n},0 \right \},1 \right \}
	\end{align}
\end{figure*}

\subsubsection{PSCA-MAP-K-Net}
\label{PSCA_MAP_Net}
In this part, we propose a PSCA-MAP-K-driven deep unrolling neural network, called PSCA-MAP-K-Net, \red{to optimize PSCA-MAP-K’s step size and} reduce the overall computation time. As illustrated in Fig. \ref{NN}, similar to PSCA-ML-K-Net, PSCA-MAP-K-Net has $U$~($> 0$) building blocks, each implementing one iteration of PSCA-MAP-K, and sets the step size of each iteration as a tunable parameter. 
\red{In PSCA-MAP-K-Net, the parameters in the independent model in (\ref{independent}) and the general model in (\ref{MVB}) (or its first-order approximation in (\ref{independent}) or second-order approximation in (\ref{second})) are set as tunable parameters if unknown and fixed parameters if known.}
All tunable parameters of PSCA-MAP-K are optimized by the neural network training. The BCE loss function and the training, validation, and test procedures of PSCA-MAP-K-Net and the thresholding module are the same as those of PSCA-ML-K-Net in Section~\ref{MLEK}.

%% file: unknown.tex
\section{MLE and MAPE-Based Device Activity Detection for Unknown Pathloss}
\label{Unknown}
In this section, we investigate MLE and MAPE-based device activity detection in the unknown pathloss case.\footnote{In the unknown pathloss case, the \red{MLE and MAPE-based} effective pathloss estimates converge to the real effective pathloss in probability, as $M$ goes to infinity \cite{chen2019covariance,kay1993fundamentalsOS}.}

\subsection{MLE-Based Device Activity Detection for Unknown Pathloss}
\label{MLEU}
In this subsection, we view $\boldsymbol{\gamma}$ as deterministic but unknown constants and investigate the MLE of $\boldsymbol{\gamma}$ for unknown $\mathbf{g}$.
\subsubsection{MLE Problem for Unknown Pathloss}
\label{MLE_Problem_U}
In this part, we present the MLE problem of $\boldsymbol{\gamma}$ for unknown $\mathbf{g}$. Given $\boldsymbol\gamma, \mathbf{Y}_{:, m}, m \in \mathcal{M}$ follow i.i.d. $C \mathcal{N}\left(0, \boldsymbol{\Sigma}_{\boldsymbol{\gamma}}\right)$ \red{with the covariance matrix} \cite{fengler2021nonbayesian,jiang2022ml}:
\begin{align}
    \label{cov_b}	\boldsymbol{\Sigma}_{\boldsymbol{\gamma}} \triangleq \mathbf{S}\boldsymbol{\Lambda}\mathbf{S}^H+\sigma^2 \mathbf{I}_L \in \mathbb{C}^{L \times L}.
\end{align}
The likelihood of $\mathbf{Y}$, denoted by $p_{\textrm{\rm ML-UD}}(\mathbf Y; \boldsymbol\gamma )$, is:
\begin{align}
\label{fY_U}
	p_{\textrm{\rm ML-UD}}(\mathbf Y; \boldsymbol\gamma ) \triangleq \frac{\exp \left(-\operatorname{tr}\left(\boldsymbol{\Sigma}_{\boldsymbol{\gamma}}^{-1} \mathbf{Y Y}^H\right)\right)}{\pi^{LM}\left| \boldsymbol{\Sigma}_{\boldsymbol{\gamma}} \right|^M}.
\end{align}
Let
\begin{align}
	\label{likelihood_Y_U}
	f_{\textrm{\rm ML-UD}}\left(\boldsymbol \gamma\right) & \triangleq  -\frac{1}{M} \log{p_{\textrm{\rm ML-UD}}(\mathbf Y; \boldsymbol\gamma )} - L \log{\pi} \nonumber \\ &= \log \left| \boldsymbol{\Sigma}_{\boldsymbol{\gamma}} \right| + \operatorname{tr}\left(\boldsymbol{\Sigma}_{\boldsymbol{\gamma}}^{-1} \widehat{\boldsymbol{\Sigma}}_{\mathbf{Y}}\right).
\end{align} 
Notice that $p_{\textrm{\rm ML-UD}}(\mathbf Y; \boldsymbol\gamma )$ in (\ref{fY_U}) and $f_{\textrm{\rm ML-UD}}\left(\boldsymbol \gamma\right)$ in (\ref{likelihood_Y_U}) share \red{similar function forms to}  $p_{\textrm{\rm ML-K}}(\mathbf Y; \boldsymbol\alpha )$ in (\ref{fY}) and $f_{\textrm{\rm ML-K}}\left(\boldsymbol \alpha\right)$ in (\ref{likelihood_Y}), respectively. The MLE problem of $\boldsymbol{\gamma}$, i.e., the maximization of $p_{\textrm{\rm ML-UD}}(\mathbf Y; \boldsymbol\gamma )$ or $\log p_{\textrm{\rm ML-UD}}(\mathbf Y; \boldsymbol\gamma )$ w.r.t. $\boldsymbol{\gamma}$, can be equivalently formulated below.
\begin{Prob}[MLE of $\boldsymbol\gamma$ for Unknown $\mathbf{g}$]
	\label{Prob_ML_U}
	\begin{align}
		\min_{\boldsymbol\gamma} & \ f_{\textrm{\rm ML-UD}}\left(\boldsymbol \gamma\right) \nonumber \\
  \label{constraint_U}
		\text {s.t.} & \ \gamma_n \ge 0, \ n \in \mathcal{N}.
	\end{align}
\end{Prob}

Problem~\ref{Prob_ML_U} differs from Problem~\ref{Prob_ML} \red{mainly} in terms of variables and constraints. The existing state-of-the-art BCD-ML-UD \cite{fengler2021nonbayesian} and PG-ML-UD \cite{wang2021efficient} for Problem~\ref{Prob_ML_U}, which resemble BCD-ML-K \cite{fengler2021nonbayesian} and PG-ML-K \cite{wang2021efficient} for Problem~\ref{Prob_ML}, respectively, can guarantee convergence to Problem~\ref{Prob_ML_U}'s stationary points but yield longer computation times.
\subsubsection{PSCA-ML-UD}
\label{PSCA_ML_Algorithm_U}

In this part, we propose an algorithm, referred to as PSCA-ML-UD, using PSCA. Specifically, at iteration $k$, we choose:
\begin{align}
	\label{Surrogate_ML_U}
	\widetilde{f}_{\textrm{\rm ML-UD}}\left(\boldsymbol{\gamma};\boldsymbol{\gamma}^{(k)}\right) \triangleq \sum_{n \in \mathcal{N}} & \widetilde{f}_{\textrm{\rm ML-UD}, n}\left(\gamma_n ;\boldsymbol{\gamma}^{(k)}\right) 
\end{align} as a strongly convex approximate function of $f_{\textrm{\rm ML-UD}}\left( \boldsymbol{\gamma}\right)$ around $\boldsymbol{\gamma}^{(k)}$ \red{w.r.t. $\boldsymbol{\gamma}$}, where
\begin{align}
	\label{Surrogate_ML_U_separate}
	\widetilde{f}_{\textrm{\rm ML-UD}, n}(\gamma_n;\boldsymbol{\gamma}^{(k)})  \triangleq & \partial_{n} f_{\textrm{\rm ML-UD}}\left(\boldsymbol{\gamma}^{(k)}\right) \left(\gamma_n - \gamma_n^{(k)}\right) \nonumber \\  + & \frac{1}{2}\left(\mathbf{s}_n^H \boldsymbol{\Sigma}_{\boldsymbol{\gamma}^{(k)}}^{-1} \mathbf{s}_n\right)^2 \left(\gamma_n - \gamma_n^{(k)}\right)^2
\end{align} is a strongly convex approximate function of \red{$f_{\textrm{\rm ML-UD}}\left(\boldsymbol{\gamma}\right)$} around $\boldsymbol{\gamma}^{(k)}$ w.r.t. $\gamma_n$, \red{and $\boldsymbol{\gamma}^{(k)}$ is obtained at iteration $k-1$}. Here,
\begin{align}
	\label{gradient_ML_U}
	\partial_{n} f_{\textrm{\rm ML-UD}}\left(\boldsymbol{\gamma}^{(k)}\right) =  \mathbf{s}_n^H \boldsymbol{\Sigma}_{\boldsymbol{\gamma}^{(k)}}^{-1} \mathbf{s}_n - \mathbf{s}_n^H \boldsymbol{\Sigma}_{\boldsymbol{\gamma}^{(k)}}^{-1} \widehat{\boldsymbol{\Sigma}}_{\mathbf{Y}} \boldsymbol{\Sigma}_{\boldsymbol{\gamma}^{(k)}}^{-1} \mathbf{s}_n.
\end{align} 
Similarly, $\left( \mathbf{s}_n^H \boldsymbol{\Sigma}_{\boldsymbol{\gamma}^{(k)}}^{-1} \mathbf{s}_n\right)^2$ converges to the $n$-th diagonal element of $\nabla^2 f_{\textrm{\rm ML-K}} (\boldsymbol{\gamma}^{(k)})$, almost surely, as $M \to \infty$. The approximate convex problem at iteration $k$ is given by:
\begin{align}
    \underset{\gamma_n \ge 0, n \in \mathcal{N}}{\operatorname{min}}  \widetilde{f}_{\textrm{\rm ML-UD}}\left(\boldsymbol{\gamma};\boldsymbol{\gamma}^{(k)}\right),
\end{align}
which can be equivalently separated into $N$ convex problems, one for each coordinate.
\begin{Prob}[Convex Approximate Problem of Problem~\ref{Prob_ML_U} for $n \in\mathcal{N}$ at Iteration $k$]
	\label{Prob_Surrogate_ML_U}
	\begin{align*}
	\red{\widetilde{\gamma}_{\textrm{\rm ML-UD}, n}^{(k)} \triangleq \underset{\gamma_n \ge 0}{\operatorname{argmin}} \  \widetilde{f}_{\textrm{\rm ML-UD}}\left(\gamma_n;\boldsymbol{\gamma}^{(k)}\right).}
	\end{align*}
\end{Prob}
\begin{Lem}[Optimal Solution of Problem \ref{Prob_Surrogate_ML_U}]
	\label{Lem_ML_U}
	The optimal point of Problem \ref{Prob_Surrogate_ML_U} is given by:
	\begin{align}
		\label{Solution_ML_U}
		\widetilde{\gamma}_{n}^{(k)} = \max \left\{\gamma_n^{(k)} - \frac{\partial_{n} f_{\textrm{\rm ML-UD}}\left(\boldsymbol{\gamma}^{(k)}\right)}{\left( \mathbf{s}_n^H \boldsymbol{\Sigma}_{\boldsymbol{\gamma}^{(k)}}^{-1} \mathbf{s}_n \right)^2},0 \right \},
	\end{align}
	where $\partial_{n} f_{\textrm{\rm ML-UD}}\left(\boldsymbol{\gamma}^{(k)}\right)$ is given by (\ref{gradient_ML_U}).
\end{Lem}
\begin{IEEEproof}
\blue{Please refer to Appendix~\ref{appendix:solution}.}
\end{IEEEproof}

Note that the slight difference between $\widetilde{\gamma}_{\textrm{\rm ML-UD}, n}^{(k)}$ in (\ref{Solution_ML_U}) and $\widetilde{\alpha}_{\textrm{\rm ML-K}, n}^{(k)}$ in (\ref{Solution_ML}) is caused by the difference between the constraints in (\ref{constraint_U}) and those in (\ref{constraint}).

Finally, we update $\gamma_n^{(k)}, n \in \mathcal{N}$ in parallel according to:
\begin{align}
	\label{Update_ML_U}
	\gamma_n^{(k+1)} = \left(1 - \rho^{(k)}\right)\gamma_n^{(k)} + \rho^{(k)} \widetilde{\gamma}_{\textrm{\rm ML-UD},n}^{(k)}, \ n \in \mathcal{N},
\end{align}
where $\rho^{(k)}$ is the step size satisfying (\ref{diminishing_rule}).
The detailed procedure is summarized in Algorithm~\ref{Algorithm_PSCA_General} (PSCA-ML-UD).
\begin{Thm}[Convergence of PSCA-ML-UD]
	\label{Convergence_ML_U}
	Every limit point of the sequence $\left\{\boldsymbol{\gamma}^{(k)}\right\}_{k \in \mathbb{Z}_{+}}$ generated by PSCA-ML-UD (at least one such point exists) is a stationary point of Problem~\ref{Prob_ML_U}.
\end{Thm}
\begin{IEEEproof}
	The proof follows from that of Theorem~\ref{Convergence_ML}.
\end{IEEEproof}

Since PSCA-ML-UD, BCD-ML-UD, and PG-ML-UD resemble PSCA-ML-K, BCD-ML-K, and PG-ML-K, respectively, their comparisons in characteristics follow those for PSCA-ML-K, BCD-ML-K, and PG-ML-K in Section \ref{PSCA_ML_Algorithm}.

Finally, we analyze the computational complexities of PSCA-ML-UD, BCD-ML-UD \cite{fengler2021nonbayesian}, and PG-ML-UD \cite{wang2021efficient} by characterizing the dominant terms and orders of the per-iteration flop counts under the assumption that $N \gg L$. The results are presented in Table~\ref{table_unknown}, and the proof is given in Appendix~\ref{appendix:computational}. \red{For fair comparison, in PG-ML-UD  \cite{wang2021efficient}, we omit the active set selection \cite{wang2021efficient} and the non-monotone line search \cite{birgin2000nonmonotone}.} According to Table~\ref{table_unknown}, PSCA-ML-UD has the lowest per-iteration computational complexity among the three MLE-based device activity detection methods and the same computational complexity as PSCA-ML-K.\footnote{\red{Compared to PG-ML-UD \cite{wang2021efficient} with the active set selection and the
non-monotone line search, PG-ML-UD has a low per-iteration computational complexity in the dominant term and the same per-iteration computational complexity in order.}}  Since (\ref{Solution_ML_U}) has a slightly simpler form than (\ref{Solution_ML}), PSCA-ML-UD is expected to have a marginally shorter computation time than PSCA-ML-K. 
\begin{table}[!ht]
	\caption{Per-Iteration Computational Complexities of PSCA-ML-UD,  BCD-ML-UD \cite{fengler2021nonbayesian}, PG-ML-UD \cite{wang2021efficient}, and PSCA-MAP-UR.}  

 \label{table_unknown}
	\centering
	\begin{tabular}{|l|l|l|}
		\hline
		\textbf{Algorithms} & \textbf{Dominant Term} & \textbf{Order}\\ \hline
  \textbf{PSCA-ML-UD} & \red{$40NL^2$} & $\mathcal{O}(NL^2)$\\ \hline
		\textbf{BCD-ML-UD} & $56NL^2$ & $\mathcal{O}(NL^2)$\\ \hline 
		\textbf{PG-ML-UD} & \red{$40NL^2$} & $\mathcal{O}(NL^2)$\\ \hline
   \textbf{PSCA-MAP-UR} & \red{$40NL^2$} & $\mathcal{O}(NL^2)$\\ \hline
	\end{tabular}
\end{table}
\subsubsection{PSCA-ML-UD-Net}
\label{PSCA_ML_Net_U}
\begin{figure}[tp]
	\begin{center}
		{\resizebox{8cm}{!}{\includegraphics{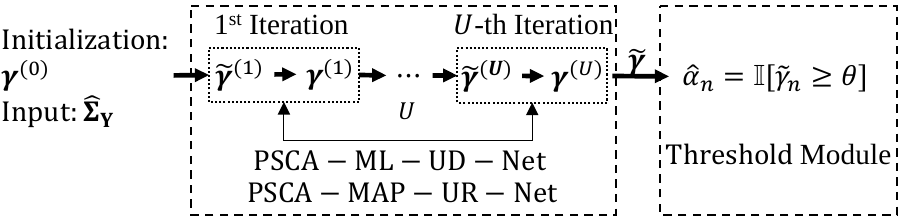}}}
	\end{center}
	\caption{Architectures of PSCA-ML-UD-Net and PSCA-MAP-UR-Net.}
 \label{NN_U}
\end{figure}

In this part, we propose a PSCA-ML-UD-driven deep unrolling neural network, called PSCA-ML-UD-Net, \red{to optimize PSCA-ML-UD’s step size and} reduce the overall computation time. As illustrated in Fig. \ref{NN_U}, \blue{similar to PSCA-ML-K-Net,} PSCA-ML-UD-Net has $U$~($> 0$) building blocks, each \blue{implementing} one iteration of PSCA-ML-UD and sets the step size of each iteration as a tunable parameter. Consider $I_P$ training, $I_Q$ validation, and $I_R$ test data samples, i.e., \red{$\left\{ \left(\mathbf{S}^{[i]}, \mathbf{g}^{[i]}, \widehat{\boldsymbol{\Sigma}}_{\mathbf{Y}}^{[i]}, \boldsymbol{\alpha}^{[i]} \right)\right\}_{i \in \mathcal{I}_j}, j \in \{P,Q,R\}$.}
Let $\widetilde{\boldsymbol{\gamma}}^{[i]}$ represent the output that corresponds to the input $\left(\mathbf{S}^{[i]}, \widehat{\boldsymbol{\Sigma}}_{\mathbf{Y}}^{[i]}\right)$. The BCE loss function given by (\ref{Loss_gamma}), as shown at the top of the next page, is used to measure the distance between $\widetilde{\boldsymbol{\gamma}}^{[i]}{/}\mathbf{g}^{[i]}$ and $\boldsymbol{\alpha}^{[i]}$, \red{for all} $i \in \mathcal{I}_{P}$  \cite{cui2021jointly}. We train PSCA-ML-UD-Net using the ADAM optimizer. Finally, we convert $\widetilde{\boldsymbol{\gamma}}^{[i]}$ to the final output $\widehat{\boldsymbol{\alpha}}^{[i]}$ by a threshold module, i.e., $\widehat{\alpha}_n^{[i]} = \mathbb{I}\left[\widetilde{\gamma}_n^{[i]} \geq \theta\right],n\in\mathcal{N},i\in\mathcal{I}_{Q}\cup\mathcal{I}_{R}$, where $\theta\in (0,1)$ is a threshold that minimizes
\red{the average error rate over the validation data $\frac{1}{NI_{Q}}\sum_{i \in \mathcal{I}_{Q}}\left\Vert \widehat{\boldsymbol{\alpha}}^{[i]} - \boldsymbol{\alpha}^{[i]}\right\Vert_1$ \cite{cui2021jointly}.}
\begin{figure*}[!ht]
    \begin{align}
        \label{Loss_gamma}
        \blue{
        \text{Loss}\left( \left\{ \widetilde{\boldsymbol{\gamma}}^{[i]}{/}\mathbf{g}^{[i]} \right\}_{i \in \mathcal{I}_P}, \left\{ \boldsymbol{\alpha}^{[i]} \right\}_{i \in \mathcal{I}_{P}}\right) = - \frac{1}{NI_{P}}\sum_{i \in \mathcal{I}_{P}}\sum_{n \in \mathcal{N}}\left( \left(\alpha_n\right)^{[i]} \log \left(\widetilde{\gamma}_n^{[i]}{/}g_n^{[i]} \right)+\left(1-\left(\alpha_n\right)^{[i]}\right) \log \left(1-(\widetilde{\gamma}_n^{[i]}{/}g_n^{[i]}\right)\right).
        }
    \end{align}
    \hrulefill
\end{figure*}
\subsection{MAPE-Based Device Activity Detection for Unknown Pathloss}
\label{MAPEU}

In this subsection, we view $\boldsymbol{\gamma}$ as realizations of independent random variables and investigate the MAPE of $\boldsymbol{\gamma}$ for unknown $\mathbf{g}$.
\subsubsection{MAPE Problem for Unknown Pathloss}
\label{MAPE_Problem_U}
In this part, we present the MAPE problem of $\boldsymbol{\gamma}$ for unknown $\mathbf{g}$.
Specifically, we consider the approximate p.d.f. of $\boldsymbol{\gamma}$, i.e., $p^{\epsilon }_{\mathbf{\Upsilon}}\left(\boldsymbol{\gamma}\right)$ in (\ref{pGamma}). The conditional p.d.f of $\mathbf{\Upsilon}$ given $\mathbf Y$, denoted by $p_{\textrm{\rm MAP-UR}}(\mathbf Y, \boldsymbol\gamma )$, is expressed as:
\begin{align}
    p_{\textrm{\rm MAP-UR}}(\mathbf Y, \boldsymbol\gamma ) = p_{\textrm{\rm ML-UD}}(\mathbf Y; \boldsymbol\gamma )p^{\epsilon }_{\mathbf{\Upsilon}}\left(\boldsymbol{\gamma}\right),
\end{align}
 where $p_{\textrm{\rm ML-UD}}(\mathbf Y; \boldsymbol\gamma )$ and $ p^{\epsilon }_{\mathbf{\Upsilon}}\left(\boldsymbol{\gamma}\right)$ are given by (\ref{fY_U}) and (\ref{pGamma}), respectively.
Let 
\begin{align}
\label{likelihood_Y_MAP_U}
	f_{\textrm{\rm MAP-UR}}\left(\boldsymbol \gamma\right) & \triangleq -\frac{1}{M} \log p_{\textrm{\rm MAP-UR}}(\mathbf Y, \boldsymbol\gamma ) - L \log{\pi} \nonumber \\
	& =  f_{\textrm{\rm ML-UD}}\left(\boldsymbol \gamma\right) + f(\boldsymbol \gamma),
\end{align}  
where $f_{\textrm{\rm ML-UD}}\left(\boldsymbol \alpha\right)$ is given by (\ref{likelihood_Y_MAP_U}), and $f(\boldsymbol \gamma)$ is given by:
\begin{align}
\label{feps}
    f(\boldsymbol \gamma) = - \frac{1}{M}\sum_{n \in \mathcal{N}}\log\left(p^{\epsilon }_{\Upsilon_n}\left(\gamma_n\right)\right).
\end{align}
Note that $f(\boldsymbol \gamma)$ in (\ref{feps}) comes from $p^{\epsilon }_{\mathbf{\Upsilon}}\left(\boldsymbol{\gamma}\right)$ in (\ref{pGamma}). Thus, the MAP problem of $\boldsymbol\gamma$, i.e., the maximization of $p_{\textrm{\rm MAP-UR}}(\mathbf Y, \boldsymbol\gamma )$ or $\log\left(p_{\textrm{\rm MAP-UR}}(\mathbf Y, \boldsymbol\gamma )\right)$ w.r.t. $\boldsymbol \gamma$, can be equivalently formulated below.

\begin{Prob}[MAPE of $\boldsymbol\gamma$ for Unknown $\mathbf{g}$]
	\label{Prob_MAP_U}
	\begin{align}
		\min_{\boldsymbol\gamma} & \quad f_{\textrm{\rm MAP-UR}}\left(\boldsymbol \gamma\right) \nonumber \\
		\text {s.t.} & \quad \gamma_n \in \left[0, g_{u,n}\right], \ n \in \mathcal{N}.
	\end{align}
\end{Prob}

Note that Problem \ref{Prob_MAP_U} has never been investigated before. Compared to Problem \ref{Prob_ML_U}, Problem \ref{Prob_MAP_U} additionally incorporates the prior distribution $p_{\mathbf{\Upsilon}}\left(\boldsymbol{\gamma}\right)$ in (\ref{pGamma}) in the objective function. Similarly, since $f_{\textrm{\rm MAP-UR}}\left(\boldsymbol \gamma\right) - f_{\textrm{\rm ML-UD}}\left(\boldsymbol \gamma\right) = f(\boldsymbol \gamma) \to 0$, as $M\to0$, the influence of $p^{\epsilon }_{\mathbf{\Upsilon}}\left(\boldsymbol{\gamma}\right)$ becomes less critical, and the observation of $\mathbf Y$ becomes more dominating, as $M$ increases. Furthermore, Problem~\ref{Prob_MAP_U} is a more challenging non-convex problem than Problem~\ref{Prob_ML_U}, due to the additional term $f(\boldsymbol \gamma)$ in (\ref{feps}). Besides, Problem \ref{Prob_MAP_U} is more challenging than Problem \ref{Prob_MAP}, since $p^{\epsilon }_{\mathbf{\Upsilon}}\left(\boldsymbol{\gamma}\right)$ in (\ref{pGamma}) is more complex than $p_{\mathbf{A}\textrm{\rm-G}}\left(\boldsymbol{\alpha}\right)$ in (\ref{MVB}) and $p_{\mathbf{A}\textrm{\rm-I}}\left(\boldsymbol{\alpha}\right)$ in (\ref{independent}).

\subsubsection{PSCA-MAP-UR}
\label{PSCA_MAP_Algorithm_U}
In this part, we propose an algorithm, referred to as PSCA-MAP-UR, using PSCA. Specifically, at iteration $k$, we choose:
\begin{align}
	\label{Surrogate_MAP_U}
	\widetilde{f}_{\textrm{\rm MAP-UR}}\left(\boldsymbol{\gamma};\boldsymbol{\gamma}^{\left(k\right)}\right) \triangleq \sum_{n \in \mathcal{N}} & \widetilde{f}_{\textrm{\rm MAP-UR},n}\left(\gamma_n ;\boldsymbol{\gamma}^{\left(k\right)}\right) 
\end{align} as a strongly convex approximate function of \red{$f_{\textrm{\rm MAP-UR}}\left( \boldsymbol{\gamma}\right)$} around $\boldsymbol{\gamma}^{\left(k\right)}$ \red{w.r.t. $\boldsymbol{\gamma}$}, where
\begin{align}
	\label{Surrogate_MAP_U_separate}
	\widetilde{f}_{\textrm{\rm MAP-UR}, n}\left(\gamma_n;\boldsymbol{\gamma}^{\left(k\right)}\right)  \triangleq & \partial_{n} f_{\textrm{\rm MAP-UR}}\left(\boldsymbol{\gamma}^{\left(k\right)}\right) \left(\gamma_n - \gamma_n^{\left(k\right)}\right) \nonumber \\   + & \frac{1}{2}\left( \mathbf{s}_n^H \boldsymbol{\Sigma}_{\boldsymbol{\gamma}^{\left(k\right)}}^{-1} \mathbf{s}_n\right)^2 \left(\gamma_n - \gamma_n^{\left(k\right)}\right)^2
\end{align} is a strongly convex approximate function of \red{$f_{\textrm{\rm MAP-UR}}\left( \boldsymbol{\gamma}\right)$} around $\boldsymbol{\gamma}^{\left(k\right)}$ w.r.t. $\gamma_n$. Here,
\begin{align}
	\label{gradient_MAP_U}
	\partial_{n} f_{\textrm{\rm MAP-UR}}\left(\boldsymbol \gamma^{\left(k\right)}\right) = \partial_{n} f_{\textrm{\rm ML-UD}}\left(\boldsymbol \gamma^{\left(k\right)}\right) - \frac{1}{M}\frac{{p^{\epsilon }_{\Upsilon_n}}'\left(\gamma_n\right)}{p^{\epsilon }_{\Upsilon_n}\left(\gamma_n\right)}, 
\end{align}  
where ${p^{\epsilon }_{\Upsilon_n}}'\left(\gamma_n\right)$ is given by (\ref{pepsd})\red{, as shown at the bottom of the next page.} 
\begin{figure*}[!hbp]
	\hrulefill
	\begin{align}
 \label{pepsd}
{p^{\epsilon }_{\Upsilon_n}}'\left(\gamma_n\right) &=
 \begin{cases}
     \frac{6\left(1-p_n\right)}{\epsilon^3}\gamma_n\left(\gamma_n-\epsilon\right), & 0, \leq \gamma_n < \epsilon \\ 
     0, & \epsilon \le \gamma_n < g_{l,n} - \epsilon, \\
     \frac{12p_n\left(\phi^{-1}\right)'\left(g_{l,n}\right)d_{u,n}}{\epsilon^3\left(d_{u, n}^2-d_{l,n}^2\right)}\left(\gamma_n-g_{l,n}\right)\left(\gamma_n-g_{l,n}+\epsilon\right) \\ - \frac{12p_n\left(\left(\left(\phi^{-1}\right)'\left(g_{l,n}\right)\right)^2 +\left(\phi^{-1}\right)''\left(g_{l,n}\right)d_{u,n}\right)}{\epsilon^2\left(d_{u, n}^2-d_{l,n}^2\right)}\left(\gamma_n-g_{l,n}+\epsilon\right)\left(3\gamma_n-3g_{l,n}+\epsilon\right), & g_{l,n} - \epsilon \le \gamma_n < g_{l,n}, \\
         -\frac{2p_n}{d_{u,n}^2 - d_{l,n}^2}\left(\left(\left(\phi^{-1}\right)'\left(\gamma_n\right)\right)^2+ \phi^{-1}\left(\gamma_n\right)\left(\phi^{-1}\right)''\left(\gamma_n\right)\right), & g_{l,n} \le \gamma_n \le g_{u,n}.
 \end{cases}
	\end{align}
\end{figure*}
Similarly, we can conclude that $\left( \mathbf{s}_n^H \boldsymbol{\Sigma}_{\boldsymbol{\gamma}^{(k)}}^{-1} \mathbf{s}_n\right)^2$ also converges to the $n$-th diagonal element of $\nabla^2f_{\textrm{\rm MAP-UR}} (\boldsymbol{\gamma}^{(k)})$, almost surely, as $M \to \infty$. 
The approximate convex problem at iteration $k$ is given by:
\begin{align}
	\underset{\gamma_n \in [0,g_{u,n}],n \in \mathcal{N}}{\operatorname{min}}  \widetilde{f}_{\textrm{\rm MAP-UR}}\left(\boldsymbol{\gamma};\boldsymbol{\gamma}^{\left(k\right)}\right),
\end{align}
which can be equivalently separated into $N$ convex problems, one for each coordinate.
\begin{Prob}[Convex Approximate Problem of Problem~\ref{Prob_Surrogate_MAP_U} for $n \in\mathcal{N}$ at Iteration $k$]
	\label{Prob_Surrogate_MAP_U}
	\begin{align*}
		\red{\widetilde{\gamma}_{\textrm{\rm MAP-UR}, n}^{\left(k\right)} \triangleq \underset{\gamma_n \in [0,g_{u,n}]}{\operatorname{argmin}} \  \widetilde{f}_{\textrm{\rm MAP-UR}, n}\left(\gamma_n;\boldsymbol{\gamma}^{\left(k\right)}\right).}
	\end{align*}
\end{Prob}
\begin{Lem}[Optimal Solution of Problem \ref{Prob_Surrogate_MAP_U}]
	\label{Lem_MAP_U}
	The optimal point of Problem \ref{Prob_Surrogate_MAP_U} is given by: 
 \begin{align}
     \label{Solution_MAP_U}
&\widetilde{\gamma}_{\textrm{\rm MAP-UR}, n}^{\left(k\right)} \nonumber \\ &= \min \left\{\max \left\{\gamma_n^{\left(k\right)} - \frac{\partial_{n} f_{\textrm{\rm MAP-UR}}\left(\boldsymbol{\gamma}^{\left(k\right)}\right)}{\left(\mathbf{s}_n^H \boldsymbol{\Sigma}_{\boldsymbol{\gamma}^{\left(k\right)}}^{-1} \mathbf{s}_n \right)^2} ,0 \right \},g_{u,n} \right \},
 \end{align}
where $\partial_{n} f_{\textrm{\rm MAP-UR}}\left(\boldsymbol{\gamma}^{\left(k\right)}\right)$ is given by (\ref{gradient_MAP_U}).
\end{Lem}

Note that $\widetilde{\gamma}_{\textrm{\rm MAP-UR},n}^{(k)}$ in (\ref{Solution_MAP_U}) shares a form similar to $\widetilde{\gamma}_{\textrm{\rm ML-UD},n}^{(k)}$ in (\ref{Solution_ML_U}) and approximately utilizes the second-order information of $f_{\textrm{\rm MAP-UR}}\left(\boldsymbol \gamma\right)$ via $\left( \mathbf{s}_n^H \boldsymbol{\Sigma}_{\boldsymbol{\gamma}^{(k)}}^{-1} \mathbf{s}_n \right)^2$. By (\ref{gradient_MAP_U}) and $\frac{1}{M}\frac{{p^{\epsilon }_{\Upsilon_n}}'\left(\gamma_n\right)}{p^{\epsilon }_{\Upsilon_n}\left(\gamma_n\right)} \to 0$, as $M \to \infty$, we have $\widetilde{\gamma}_{\textrm{\rm MAP-UR},n}^{(k)} \to \widetilde{\gamma}_{\textrm{\rm ML-UD},n}^{(k)}$, as $M \to \infty$.
Thus, the influence of the prior distribution of $\boldsymbol{\gamma}$ decreases as $M$ increases.

Finally, we update $\gamma_n^{(k)}, n \in \mathcal{N}$ in parallel according to:
\begin{align}
	\label{Update_MAP_U}
	\gamma_n^{(k+1)} = \left(1 - \rho^{(k)}\right)\gamma_n^{(k)} + \rho^{(k)} \widetilde{\gamma}_{\textrm{\rm MAP-UR},n}^{(k)}, \ n \in \mathcal{N},
\end{align}
where $\rho^{(k)}$ is the step size satisfying (\ref{diminishing_rule}).
The detailed procedure is summarized in Algorithm~\ref{Algorithm_PSCA_General} (PSCA-MAP-UR).
\begin{Thm}[Convergence of PSCA-MAP-UR]
	\label{Convergence_MAP_U}
	Every limit point of the sequence $\left\{\boldsymbol{\gamma}^{(k)}\right\}_{k \in \mathbb{Z}_{+}}$ generated by PSCA-MAP-UR (at least one such point exists) is a stationary point of Problem \ref{Prob_MAP_U}.
\end{Thm}
\begin{IEEEproof}
	The proof follows from that of Theorem~\ref{Convergence_MAP}.
\end{IEEEproof}

Finally, we analyze the per-iteration computational complexity of the per-iteration flop counts of PSCA-MAP-UR under the assumption that $N \gg L$. The results are illustrated in Table \ref{table_unknown}, \red{and the proof is given in Appendix~\ref{appendix:computational}.} Furthermore, according to Table \ref{table_unknown}, PSCA-MAP-UR has the same computational complexity in the dominant term \red{of flop counts} as PSCA-ML-UD. However, since (\ref{Solution_MAP_U}) is more complex than (\ref{Solution_ML_U}), PSCA-MAP-UR has a higher computational complexity in \red{absolute} flop counts than PSCA-ML-UD.

\subsubsection{PSCA-MAP-UR-Net}
\label{PSCA_MAP_Net_U}
In this part, we propose a PSCA-MAP-UR-driven deep unrolling neural network, called PSCA-MAP-UR-Net, \red{to optimize PSCA-MAP-UR’s step size and} reduce the overall computation time. As illustrated in Fig. \ref{NN_U}, similar to PSCA-ML-UD-Net, PSCA-MAP-UR-Net has $U$~($> 0$) building blocks, each implementing one iteration of PSCA-MAP-UR and sets the step size of each iteration as a tunable parameter. All tunable parameters are optimized by the neural network training. The BCE loss function and the training, validation, and test procedures of PSCA-MAP-UR-Net and the thresholding module are the same as those of PSCA-ML-UD-Net in Section~\ref{MLEU}.

%% file: numerical.tex
\section{Numerical Results}
\label{Numerical_Results}

\begin{figure}[tp]
\centering
\begin{minipage}{0.24\textwidth}
\centering
{\resizebox{4.34cm}{!}{\includegraphics{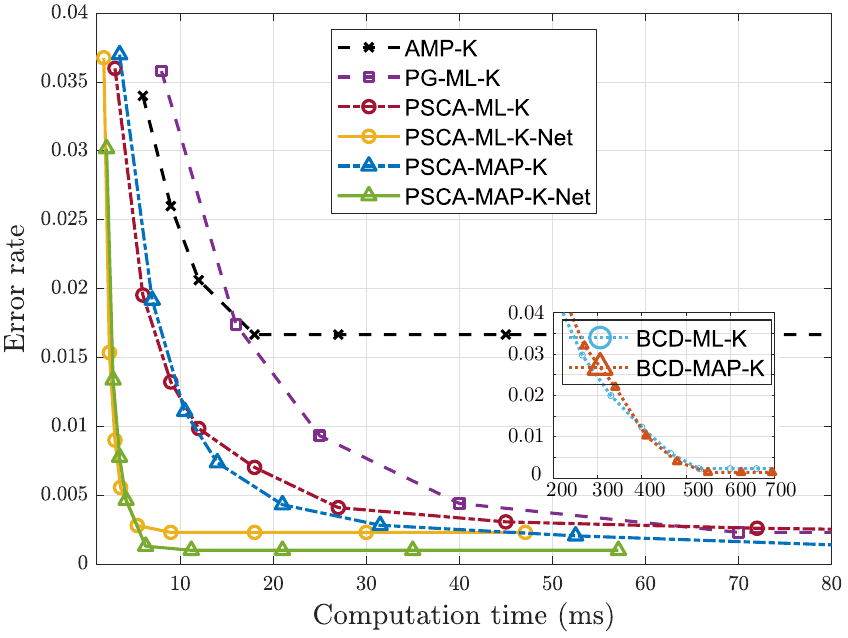}}}
\caption{\small{Error rate versus computation time at $N=1000, L=40, M=256$, and $P=23$dBm for known $\mathbf{g}$ under the i.i.d. activity model.}}
\label{Fig_Tradeoff}
\end{minipage}
\begin{minipage}
{0.24\textwidth}
\vspace{-0.25cm}
\centering
{\resizebox{4.34cm}{!}{\includegraphics{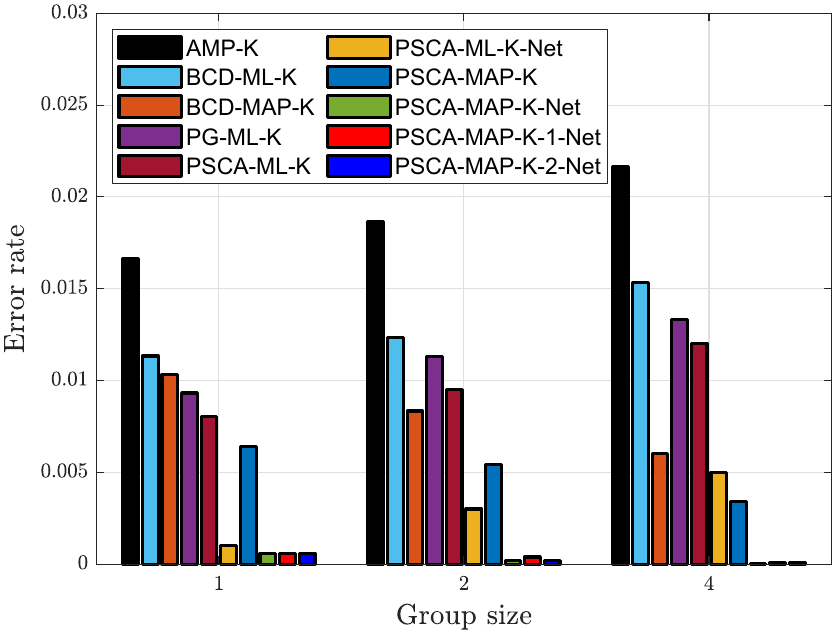}}}
\caption{\small{Error rate versus group size at $N=1000, L=40, M=256$, and $P=23$dBm  for known $\mathbf{g}$ under the group activity model.}}
\label{Fig_group}
\end{minipage}
\end{figure}

In this section, we evaluate the error rates and computation \red{times} of the proposed PSCA-based algorithms and PSCA-Nets for known and unknown $\mathbf{g}$.\footnote{Due to the page limit, we no longer compare PSCA-Nets with the BCD-driven and AMP-driven deep unrolling networks, namely BCD-Nets and AMP-Nets. Since each PSCA-based algorithm outperforms the corresponding BCD-based and AMP-based algorithms, each PSCA-Net naturally outperforms the corresponding BCD-Net and AMP-Net. We do not present the error rate versus the computation time for unknown $\mathbf{g}$, which resembles Fig.~\ref{Fig_Tradeoff}.} For known $\mathbf{g}$, we consider five state-of-the-art  algorithms, i.e., AMP-K \cite{liu2018massive}, BCD-ML-K \cite{fengler2021nonbayesian}, BCD-MAP-K (general) \cite{jiang2020mapbased,jiang2022ml,jiang2023statistical}, BCD-MAP-K (independent) \cite{jiang2020mapbased, jiang2022ml}, and PG-ML-K \cite{wang2021efficient}, as the baseline schemes. For unknown $\mathbf{g}$, we take three state-of-the-art algorithms, i.e., EM-AMP-UD \cite{hara2022blinda}, BCD-ML-UD \cite{fengler2021nonbayesian}, and PG-ML-UD \cite{wang2021efficient}, as the baseline schemes. The computational complexities of AMP-K and EM-AMP-UD are $\mathcal{O}(NLM)$. The computational complexities of the other algorithms have been illustrated in Table \ref{table_known} and Table \ref{table_unknown}. 

\red{In the simulation, we set $d_{l, n}=20$m, $d_{u, n}=200$m, and $p_n = 0.05$, for all $n \in \mathcal{N}$.} The path loss model, i.e., $\phi(\cdot)$, is given by $g_n = 10\eta\log_{10}(\frac{4\pi d_n}{\lambda})$ in dB, for all $n \in \mathcal{N}$, where the pathloss exponent $\eta$ is 2.5, and the wavelength $\lambda$ is 0.086 m. The bandwidth is 1 MHz. The noise power spectral density is -174 dBm/Hz. Thus, the noise power is -114 dBm, i.e., $\sigma^2 = 10^{-11.4}$. For all $n \in \mathcal{N}$, the elements of $\mathbf{s}_n$ are generated according to i.i.d. $\mathcal{CN}(0, 1)$ and normalized with the norm being $\sqrt{L}$ \cite{cui2021jointly}. We calculate the average transmit signal-to-noise ratio (SNR) by $10\log_{10}\left(\frac{P}{\sigma^{2}}\right)$.  We consider the i.i.d. activity model, i.e., the independent model in (\ref{independent}), and the group activity model in \cite{jiang2022ml}. \red{Unless otherwise stated, we set $N=1000, L=40, M=256$, and $P=23$dBm, and set the numbers of iterations of the proposed PSCA-based algorithms, their implementations in PSCA-Nets (i.e., $U$), PG-ML-K (UD), BCD-ML-K (UD), BCD-MAP-K, AMP-K, and EM-AMP-UD  to be 30, 15, 5, 5, 5, 20, and 20, respectively. The number of iterations for each method is chosen to achieve a reasonable error rate and computation time tradeoff.\footnote{\red{In the simulation, for PG-ML-K and PG-ML-UD \cite{wang2021efficient}, we conduct the non-monotone line search to obtain the step sizes \cite{birgin2000nonmonotone}.}}} 

For each method, we choose $\mathbf{0}$ as the initial point for the estimate. The step size of each PSCA and the initial step size of each PSCA-Net are selected according to $\rho^{(k+1)} = \rho^{(k)}(1-0.5\rho^{(k)})$ with $\rho^{(0)} = 0.5$. In the simulation, we consider $I_P = 6000$ training, $I_Q = 2000$ validation, and $I_R = 2000$ test data samples, which are the realizations of the random variables.
The maximization epochs, learning rate, and batch size in the training process are 100, 0.001, and 64, respectively, \cite{cui2021jointly}.  
For each method except for AMP-K and EM-AMP-UD, the optimal threshold $\theta$ is set as the one that minimizes the average error rate over the validation samples. For AMP-K and EM-AMP-UD, the threshold is given by \cite{liu2018massive} and \cite{hara2022blinda}, respectively.
The average error rate of each method is obtained by averaging over the test samples.

\begin{figure}[htbp]
    \centering
   \begin{minipage}{0.49\textwidth}
\centering
\subfigure[\scriptsize{\red{$\left\Vert\boldsymbol\alpha^{(t)} - \boldsymbol\alpha^{(t-1)}\right\Vert_2$ versus $t$ for known $\mathbf{g}$.}}]
		{\resizebox{4.28cm}{!}{\includegraphics{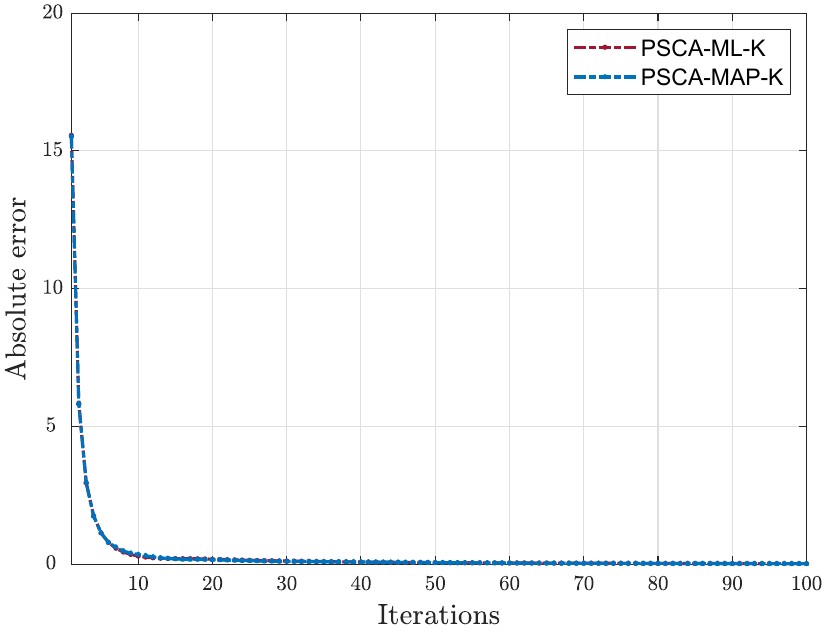}}}	
\subfigure[\scriptsize{\red{$\left\Vert\boldsymbol\gamma^{(t)} - \boldsymbol\gamma^{(t-1)}\right\Vert_2$ versus $t$ for unknown $\mathbf{g}$.}}]
		{\resizebox{4.36cm}{!}{\includegraphics{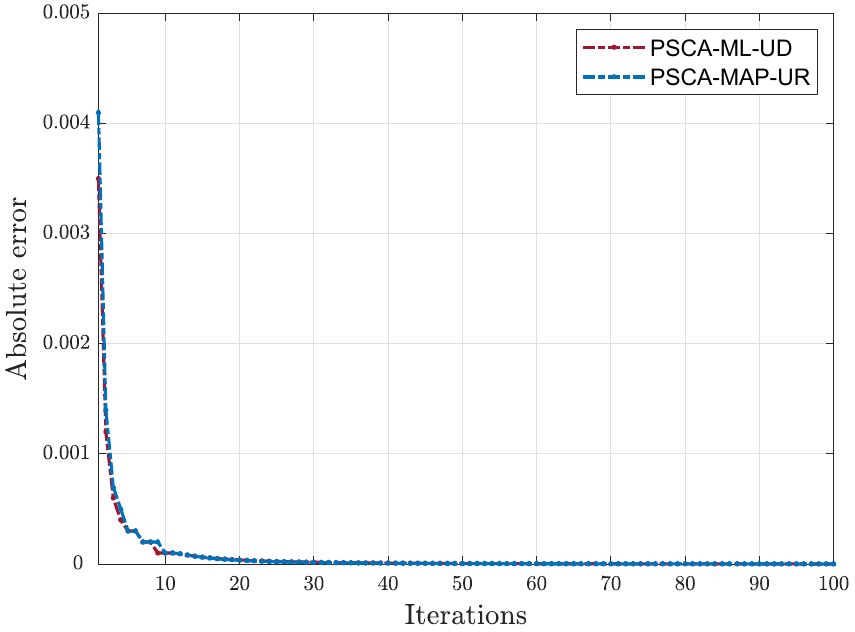}}}	
  \caption{\small{Convergence at $N=1000, L=40, M=256$, and $P=23$dBm under the i.i.d. activity model.}}
	\label{conver}
\end{minipage}
\end{figure}

\begin{figure}[htbp]
    \centering
    \begin{minipage}{0.49\textwidth}
\centering
        \subfigure[\scriptsize{\red{Generalization ability w.r.t. $N$ at $L=40$ and $P=23$dBm.}}]
		{\resizebox{4.32cm}{!}{\includegraphics{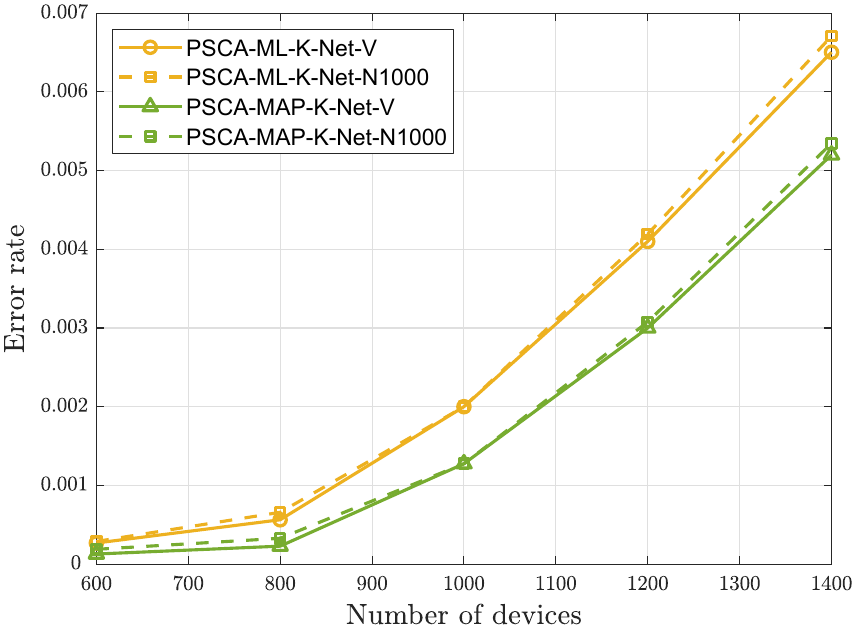}}}	\subfigure[\scriptsize{\red{Generalization ability w.r.t. $L$ at $N=1000$ and $P=23$dBm.}}]
		{\resizebox{4.28cm}{!}{\includegraphics{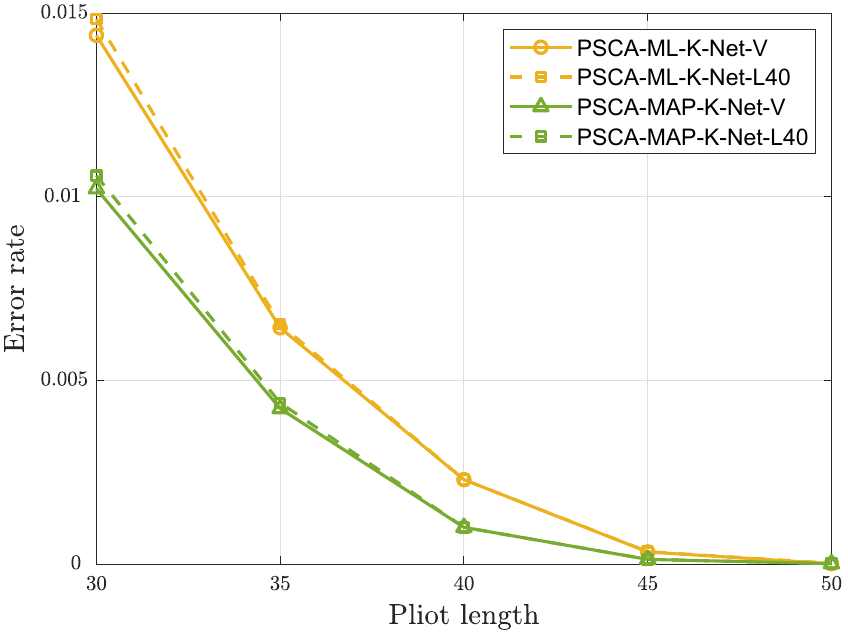}}}
        \subfigure[\scriptsize{\red{Generalization ability w.r.t. $P$ at $N=1000$ and $L=40$.}}]
		{\resizebox{4.32cm}{!}{\includegraphics{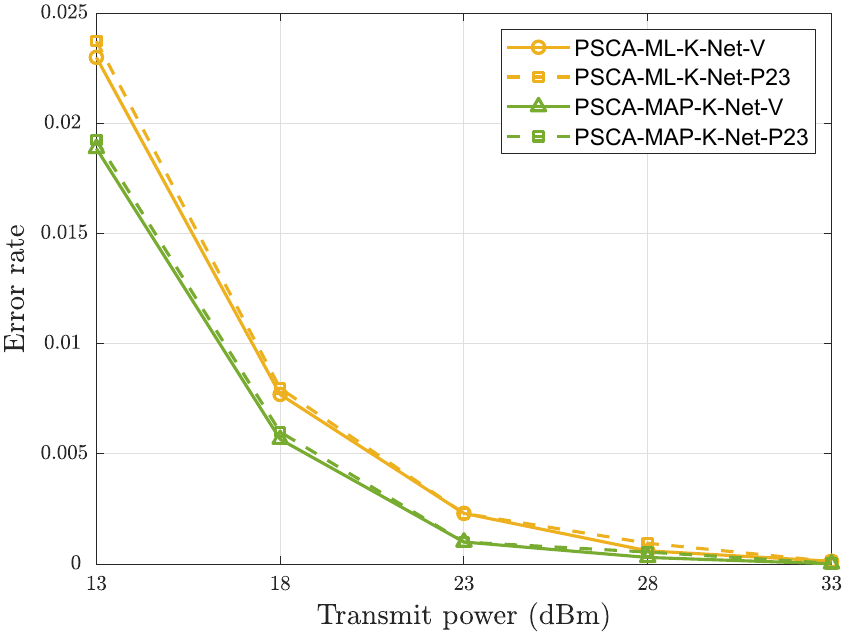}}}
  \caption{\small{\red{Generalization ability w.r.t. $N$, $L$, and $P$ at $M=256$ for known $\mathbf{g}$ under the i.i.d. activity model.}}}
  \label{add_11}
\end{minipage}
\end{figure}
\begin{figure}[htbp]
    \centering
    \begin{minipage}{0.49\textwidth}
\centering
\subfigure[\scriptsize{\red{Generalization ability w.r.t. $N$ at $L=40$ and $P=23$dBm.}}]
		{\resizebox{4.32cm}{!}{\includegraphics{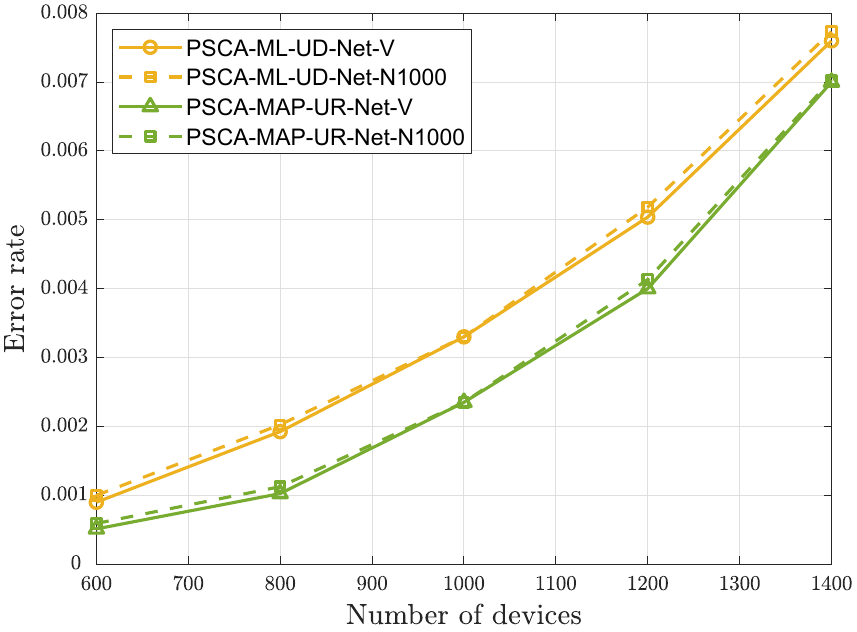}}}	\subfigure[\scriptsize{\red{Generalization ability w.r.t. $L$  at $N=1000$ and $P=23$dBm.}}]
		{\resizebox{4.28cm}{!}{\includegraphics{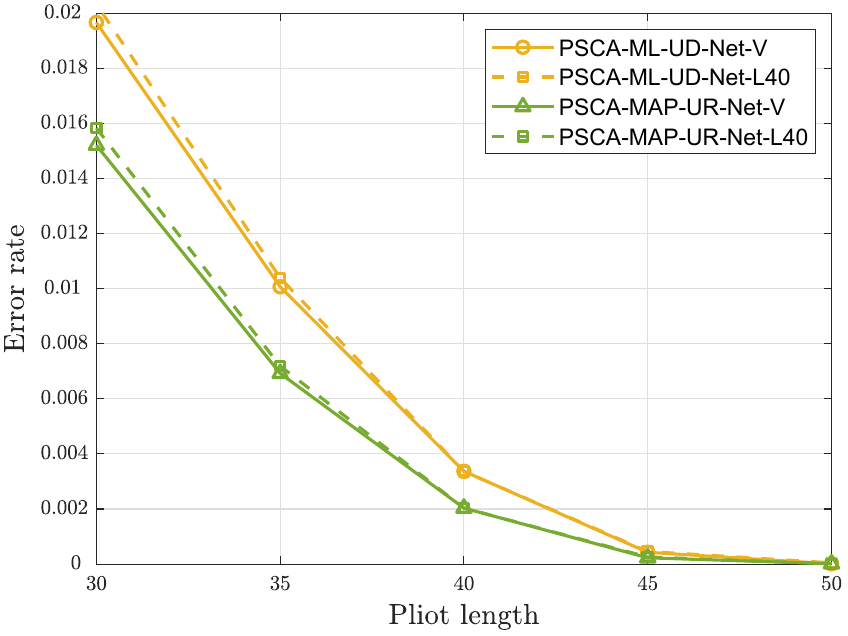}}}
        \subfigure[\scriptsize{\red{Generalization ability w.r.t. $P$  at $N=1000$ and $L=40$.}}]
		{\resizebox{4.32cm}{!}{\includegraphics{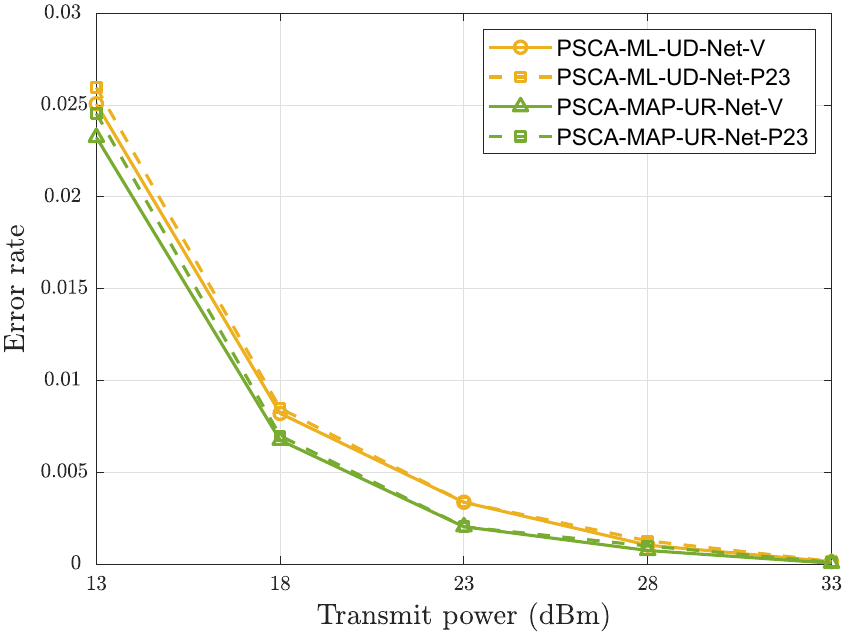}}}
  \caption{\small{\red{Generalization ability w.r.t. $N$, $L$, and $P$ at $M=256$ for unknown $\mathbf{g}$ under the i.i.d. activity model.}}}
  \label{add_21}
\end{minipage}
\end{figure}

\begin{figure*}[tp]
\centering
\begin{minipage}{0.49\textwidth}
\centering
\subfigure[\scriptsize{Error rate versus $L$.}]
		{\resizebox{4.34cm}{!}{\includegraphics{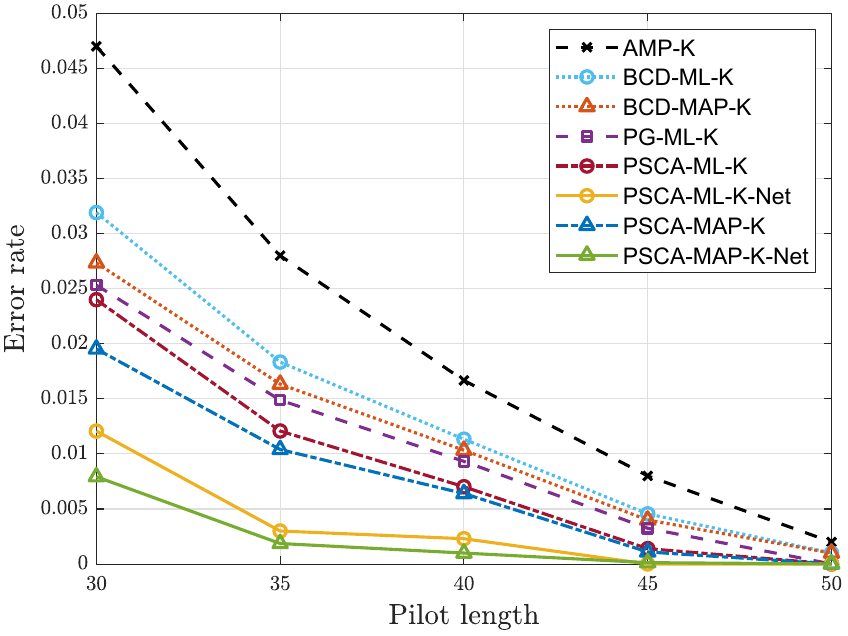}}}	\subfigure[\scriptsize{Computation time versus $L$.}]
		{\resizebox{4.24cm}{!}{\includegraphics{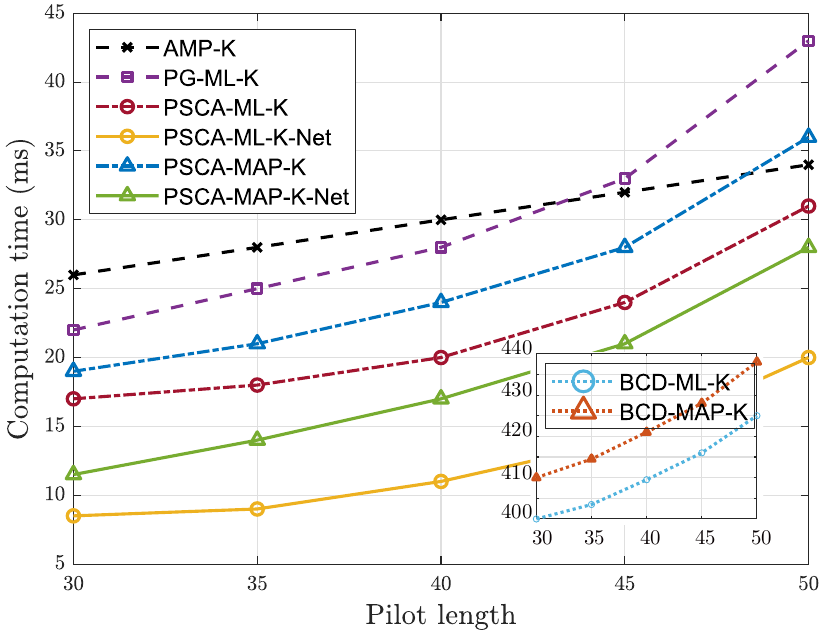}}}
  \caption{\small{Error rate and computation time versus $L$ at $N=1000, M=256$, and $P=23$dBm for known $\mathbf{g}$ under the i.i.d. activity model.}}
	\label{Fig_L_K}
\end{minipage}
\begin{minipage}{0.49\textwidth}
\centering
\subfigure[\scriptsize{Error rate versus $M$.}]
		{\resizebox{4.34cm}{!}{\includegraphics{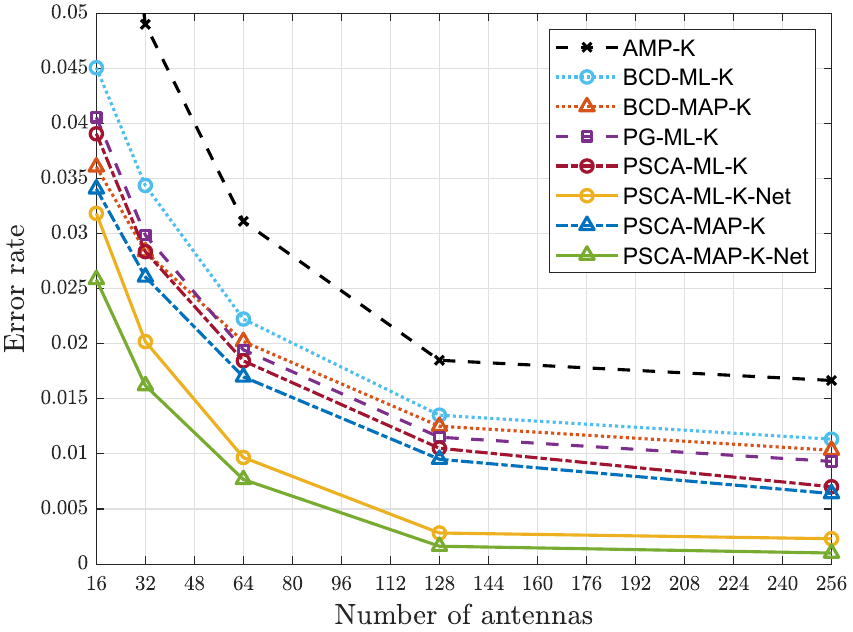}}}
\subfigure[\scriptsize{Computation time versus $M$.}]
		{\resizebox{4.24cm}{!}{\includegraphics{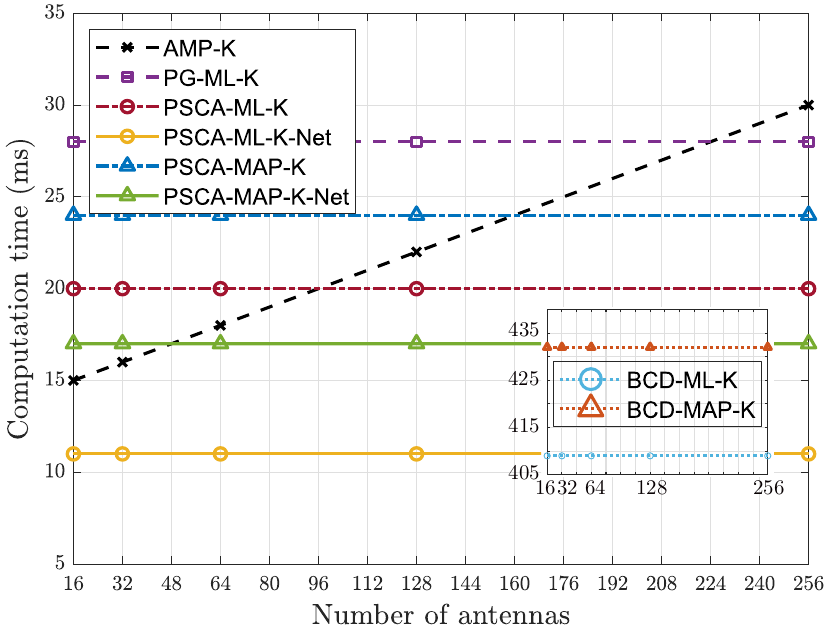}}}
\caption{\small{Error rate and computation time versus $M$ at $N=1000, L=40$, and $P=23$dBm for known $\mathbf{g}$ under the i.i.d. activity model.}}
\label{Fig_M_K}
\end{minipage}
\end{figure*}

\begin{figure*}[tp]
\centering
\begin{minipage}{0.49\textwidth}
\vspace{-0.29cm}
\centering
\subfigure[\scriptsize{Error rate versus $L$.}]
		{\resizebox{4.34cm}{!}{\includegraphics{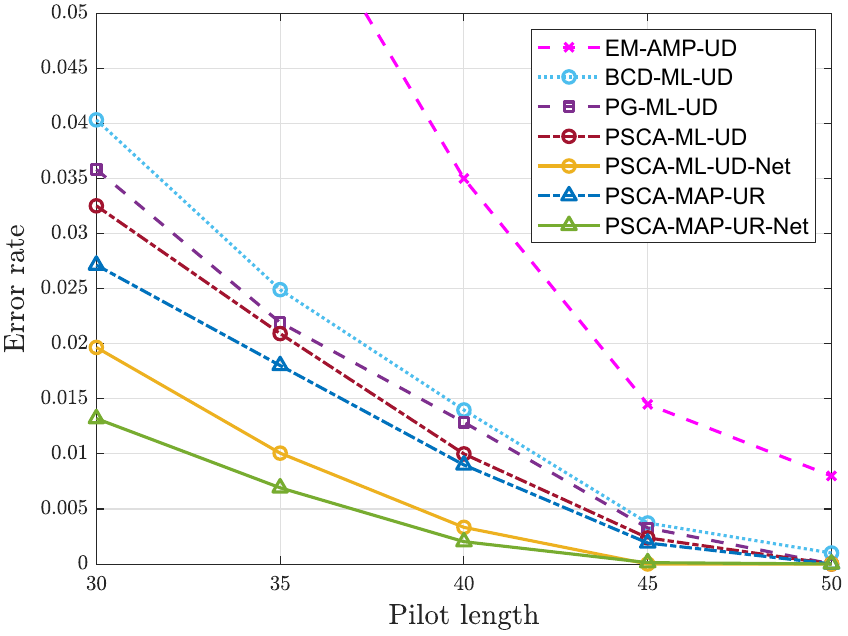}}}	\subfigure[\scriptsize{Computation time versus $L$.}]
		{\resizebox{4.24cm}{!}{\includegraphics{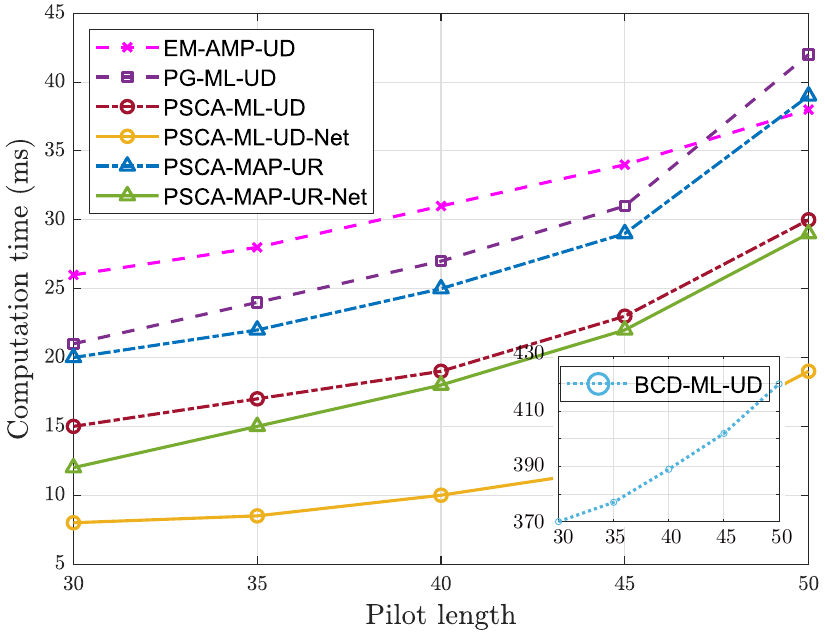}}}
  \caption{\small{Error rate and computation time versus $L$ at $N=1000, M=256$, and $P=23$dBm for unknown $\mathbf{g}$ under the i.i.d. activity model.}}
	\label{Fig_L_U}
\end{minipage}
\begin{minipage}{0.49\textwidth}
\centering
\subfigure[\scriptsize{Error rate versus $M$.}]
		{\resizebox{4.34cm}{!}{\includegraphics{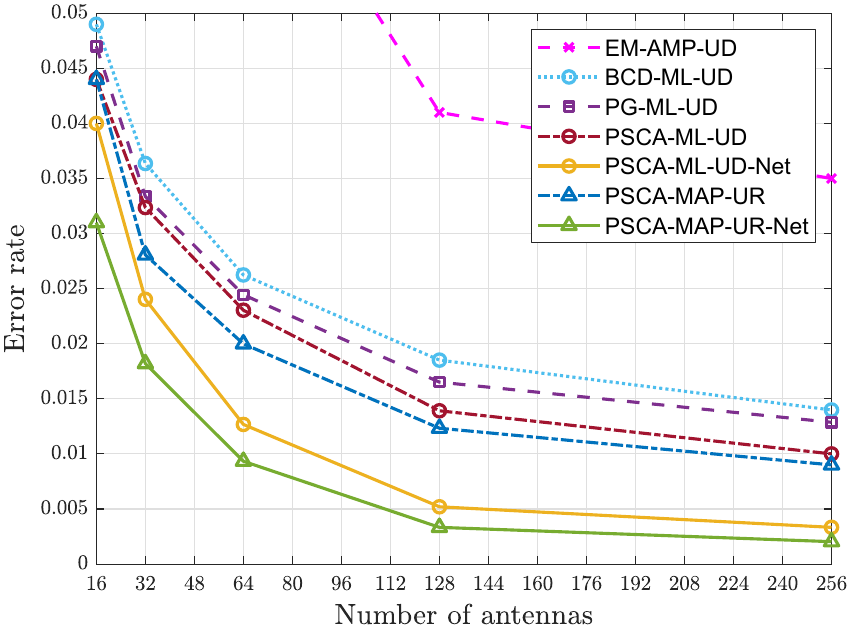}}}
\subfigure[\scriptsize{Computation time versus $M$.}]
		{\resizebox{4.24cm}{!}{\includegraphics{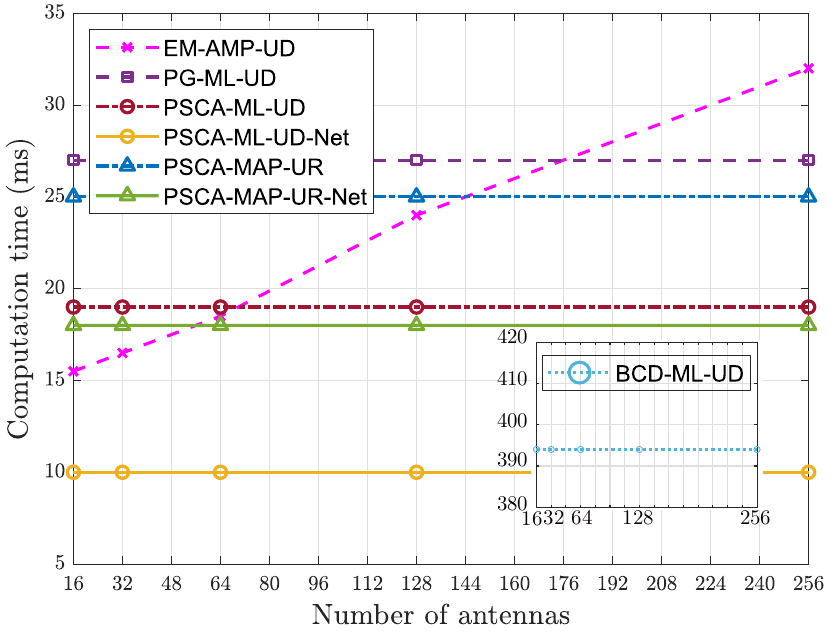}}}
\caption{\small{Error rate and computation time versus $M$ at $N=1000, L=40$, and $P=23$dBm for unknown $\mathbf{g}$ under the i.i.d. activity model.}}
\label{Fig_M_U}
\end{minipage}
\end{figure*}

Fig.~\ref{conver} illustrates the convergence results of the proposed PSCA's convergence for known and unknown $\mathbf{g}$, verifying Theorems 2-5.
Fig.~\ref{add_11} and Fig.~\ref{add_21} plot the error rates of these PSCA-Nets versus the number of devices $N$, pilot length $L$, and transmit power $P$, respectively, for known and unknown $\mathbf{g}$, respectively. Consider any \text{(i, j)} $\in$ \{(ML, K), (ML, UD), (MAP, K), (MAP, UR)\}. PSCA-i-j-Net-N1000, PSCA-i-j-Net-L40, and PSCA-i-j-Net-P23 represent PSCA-i-j-Net trained at $N=1000$, $L=40$, and $P=23$, respectively.  PSCA-i-j-Net-N-V, PSCA-i-j-Net-L-V, and PSCA-i-j-Net-P-V represent PSCA-i-j-Net retrained for each $N$, $L$, and $P$, respectively. From Fig.~\ref{add_11} and Fig.~\ref{add_21}, we can see that the gap between PSCA-i-j-Net-k-l and PSCA-i-j-Net-k-V is small for all $\text{(k, l)} \in \{\text{(N, 1000), (L, 40), (P, 23)}\}$ and $\text{k} \in \{\text{N, L, P}\}$, demonstrating the generalization abilities of the proposed PSCA-Nets. Therefore, we do not need to retrain PSCA-Nets unless the parameters $N$, $L$, and $P$ change dramatically.

Fig. \ref{Fig_L_K}, Fig.~\ref{Fig_M_K}, Fig.~\ref{Fig_L_U}, and Fig.~\ref{Fig_M_U} plot the error rate and computation time of each method versus the pilot length $L$ and the number of antennas $M$ for known and unknown $\mathbf{g}$ under the i.i.d. activity model. We can observe from Fig.~\ref{Fig_L_K} (a), Fig.~\ref{Fig_M_K} (a), Fig.~\ref{Fig_L_U} (a), and Fig.~\ref{Fig_M_U} (a) that the error rate of each method decreases with $L$ and $M$. We can observe from \red{Fig.~\ref{Fig_L_K} (b) and Fig.~\ref{Fig_L_U} (b)} that the computation time of each ML or MAPE-based algorithm increases with $L$ quadratically, whereas the computation times of AMP-K and EM-AMP-UD increase with $L$ linearly. Besides, we can observe from Fig.~\ref{Fig_M_K} (b) and Fig.~\ref{Fig_M_U} (b) that the computation time of each MLE or MAPE-based algorithm does not change with $M$, whereas the computation time of each \red{MMSE-based} algorithm increases with $M$ linearly. These observations are in accordance with the computational complexity analysis in Table~\ref{table_known} and Table~\ref{table_unknown}.

Furthermore, we can make the following observations from Fig.~\ref{Fig_L_K}, Fig.~\ref{Fig_M_K}, Fig.~\ref{Fig_L_U}, and Fig.~\ref{Fig_M_U}. 
For classical estimation, PSCA-ML-K (UD) reduces the error rate (by up to 69.2\%) and computation time (by up to 96.1\%) at all $L$ and $M$ compared to BCD-ML-K (UD) and PG-ML-K (UD), thanks to the parallel optimization algorithm design and the use of higher-order information. 
For Bayesian estimation, PSCA-MAP-K reduces the error rate (by up to 86.6\%) and computation time (by up to 95.1\%) compared to AMP-K and BCD-MAP-K. 
PSCA-MAP-K (UR) reduces the error rate (by up to 21.4\%) using the prior distribution of $\mathbf{A}$ at the sacrifice of additional computation time (by up to 33.3\%) compared to PSCA-ML-K (UD), and the error rate reduction decreases with $L$ and $M$. 
Each PSCA-Net further reduces the error rate (by up to 90.8\%) and computation time (by up to 50.1\%) compared to the underlying PSCA-based algorithm at all considered values of $L$ and $M$, \red{owing} to the optimization of the step sizes by the neural network training. 
Each method for unknown $\mathbf{g}$ performs worse than the counterpart for known $\mathbf{g}$, \red{as expected}.

Fig. \ref{Fig_Tradeoff} shows the error rate versus the computation time for known $\mathbf{g}$ under the i.i.d. activity model. In Fig.~\ref{Fig_Tradeoff}, the points on each curve are obtained by changing the number of iterations of the corresponding algorithm. 
The error rates of all MLE-based (MAPE-based) algorithms almost reduce to the same value, as the computation time incereases, since they can be shown to converge to the stationary points of Problem~\ref{Prob_ML} (Problem~\ref{Prob_MAP}). 
Besides, the proposed PSCA-based algorithms and PSCA-Nets achieve much better \red{error rate and computation time tradeoffs} than the baseline algorithms.

Fig. \ref{Fig_group} plots the error rate versus the group size for known $\mathbf{g}$ under the i.i.d. activity model. In Fig. \ref{Fig_group}, PSCA-MAP-K-1-Net and PSCA-MAP-K-2-Net represent the PSCA-MAP-K-Net under the first-order approximation given by (\ref{independent}) and the second-order approximation \red{given by} (\ref{second}), respectively. From Fig.~\ref{Fig_group}, we can see that the error rates of the MLE-based and MMSE-based algorithms increase with the group size, as the variance of the number of active devices increases, and the correlation among devices is not leveraged \cite{jiang2022ml}. In contrast, the error rates of the MAPE-based algorithms decrease with the group size as the exploitation of correlation successfully narrows down the set of possible activity states \cite{jiang2022ml,cui2021jointly}. \red{Besides, PSCA-MAP-K-Net can achieve a lower error rate than PSCA-MAP-K-2-Net, and PSCA-MAP-K-2-Net can achieve a lower error rate than PSCA-MAP-K-1-Net, as more prior information of device activities is utilized.}

\section{Conclusion}
\label{conclusion}
This paper obtained fast and accurate device activity detection methods for massive grant-free access by efficiently solving four non-convex MLE and MAPE-based device activity detection problems in the known and unknown pathloss cases, using optimization and learning techniques. First, for each non-convex estimation problem, we proposed an innovative PSCA-based iterative algorithm, which allows parallel computations, uses up to the objective function’s second-order information, converges to the problem’s stationary points, and has a low per-iteration computational complexity. Next, for each PSCA-based iterative algorithm, we proposed a deep unrolling neural network implementation, which elegantly marries the underlying PSCA-based algorithm’s parallel computation mechanism with the parallelizable neural network architecture and effectively optimizes its step sizes based on vast data samples. Analytical and numerical results validated the substantial gains of the proposed methods in the error rate and computation time over the state-of-the-art methods, revealing their significant values for mMTC and URLLC in 5G and beyond.

%% file: appendix.tex
\begin{appendices}
\section{Derivations of $p_{G_n}(g_n)$ and $p_{\Upsilon_n}\left(\gamma_n\right)$}
\label{appendix:pdf}

Let $F_{G_n}(g_n),F_{D_n}(d_n),$ and $F_{\Upsilon_n}\left(\gamma_n\right)$ denote the cumulative probability functions (c.d.f.s) of $G_n, D_n$, and $\Upsilon_n$, respectively. First, we derive $F_{G_n}(g_n)$. We have:
\begin{align}
    \label{cdf}
    & F_{G_n}(g_n) {=} \textrm{Pr}\left[G_n \leq g_n\right] \overset{(\textrm{E1})}{=} \textrm{Pr}\left[\phi\left(D_n\right) \leq g_n\right] \nonumber \\  \overset{(\textrm{E2})}{=} \ & \textrm{Pr}\left[D_n \geq \phi^{-1}\left(g_n\right)\right] {=} 1 - F_{D_n}(\phi^{-1}\left(g_n\right)),
\end{align}
where (E1) is due to (\ref{gn}), and (E2) is due to the monotonically decreasing property of $\phi$. Thus, we have:
\begin{align}
\label{pdfg}
p_{G_n}\left(g_n\right) {=} F_{G_n}'(g_n)\overset{(\textrm{E3})}{=}- p_{D_n}\left(\phi^{-1}\left(g_n\right)\right) \left(\phi^{-1}\right)'\left(g_n\right),
\end{align}
where (E3) is due to (\ref{cdf}) and the chain rule of derivation. By substituting (\ref{dispdf}) and (\ref{gn}) into (\ref{pdfg}), we can obtain $p_{G_n}\left(g_n\right)$ in (\ref{pg}).

Next, we derive $F_{\Upsilon_n}\left(\gamma_n\right)$. We have:
\begin{align}
\label{cdfup}
    & F_{\Upsilon_n}\left(\gamma_n\right) {=} \textrm{Pr}\left[\Upsilon_n \leq \gamma_n\right]  = \textrm{Pr}\left[A_n G_n \leq \gamma_n\right] \nonumber \\ 
    \overset{(\textrm{E4})}{=} \ & \textrm{Pr}\left[A_n G_n\leq \gamma_n | A_n = 0\right] \textrm{Pr}\left[A_n = 0 \right]  \nonumber \\ & + \textrm{Pr}\left[A_n G_n\leq \gamma_n |A_n = 1\right]\textrm{Pr}\left[A_n = 1 \right]
    \nonumber \\ 
    {=} \ &
    \begin{cases}
    1 - p_n, & \gamma_n = 0, \\
    1 - p_n +   p_nF_{G_n}(\gamma_n), & \gamma_n \in \left[g_{l,n},g_{u,n}\right],
    \end{cases}
\end{align}
where (E4) is due to the total probability theorem. Thus, we have:
\begin{align*}
    & p_{\Upsilon_n}\left(\gamma_n\right) {=} F_{\Upsilon_n}'\left(\gamma_n\right)
    \overset{(\textrm{E5})}{=} &
    \begin{cases}
    1 - p_n, & \gamma_n = 0, \\
    p_n p_{G_n}\left(\gamma_n\right), & \gamma_n \in \left[g_{l,n},g_{u,n}\right],
    \end{cases}
\end{align*}
where the first line of (E5) is from the discrete measure theory, and the second line of (E5) is due to (\ref{pdfg}) and (\ref{cdfup}).

\section{Proof of Theorem~\ref{properties}}
\label{appendix:approximate}
First, by (\ref{pg}) and (\ref{peps}), we can prove statement (i). Second, we prove statement (ii). Since $H_1\left(\gamma_n\right)$,  $H_2\left(\gamma_n\right)$, and $- \frac{2p_n}{d_{u,n}^2 - d_{l, n}^2}\phi^{-1}\left(\gamma_n\right)\left(\phi^{-1}\right)'\left(\gamma_n\right)$ are all differentiable in $[0,\epsilon)$, $\left(g_{l, n} - \epsilon, g_{l, n}\right)$, and $\left(g_{l, n} \le \gamma_n \le g_{u, n}\right]$, respectively, it remains to show the function values and gradient values at the junctions are equal, respectively, i.e., $
    H_1\left(\epsilon\right) = H_2\left(g_{l, n} - \epsilon\right) = 0, H_2\left(g_{l, n}\right) = p_{\Upsilon_n}\left(g_{l, n}\right),
   {H_1'}_{+}\left(0\right) = {H_1'}_{-}\left(\epsilon\right) =0,
  {H_2'}_{+}\left(g_{l, n} - \epsilon\right) =0, {H_2'}_{-}\left(g_{l, n}\right) = p_{\Upsilon_n}'\left(g_{l, n}\right).$
Thus, we can draw statement (ii). Third, since for all $\gamma_n \in (0, g_{l, n})$, $p^{\epsilon }_{\Upsilon_n}\left(\gamma_n\right) \to 0$, as $\epsilon \to 0$, and for all $ \gamma_n \in \{0\} \cup \left[g_{l, n}, g_{u, n}\right]$, $p^{\epsilon }_{\Upsilon_n}\left(\gamma_n\right) =  p_{\Upsilon_n}\left(\gamma_n\right)$, we can prove statement (iii). 
\color{black}
\section{Proofs of Lemma~\ref{Lem_ML}, Lemma~\ref{Lem_MAP}, Lemma~\ref{Lem_ML_U}, and Lemma~\ref{Lem_MAP_U}}
\label{appendix:solution}
In Appendix \ref{appendix:solution}, we define $f\left(x\right) = a(x-b)+c(x-b)^2$ with $x,a,b \in \mathbb{R}$, and $c>0$. To show Lemma~\ref{Lem_ML}, Lemma~\ref{Lem_MAP}, Lemma~\ref{Lem_ML_U}, and Lemma~\ref{Lem_MAP_U}, we first show the following lemmas.
\begin{Lem}
	\label{Analytical_1}
	The optimal solution of $\underset{x \in [x_l,x_u]}{\operatorname{min}}f\left(x\right)$, where $x_l,x_u \in \mathbb{R}$, is given by $x^* = \min \left\{\max \left\{b - \frac{a}{2c} ,x_l \right \},x_u \right \}$.
\end{Lem}    

\begin{IEEEproof}
Notice that $f\left(x\right)$ is a strongly convex function.  An optimal point $x^*$ satisfies one of the following conditions \cite[Proposition 3.1.1]{bertsekasnonlinear}: 
	\begin{align}
		\label{case1}
		f'\left(x^*\right) = 0, &  \text{ with } x_l<x^*<x_u, \\
		\label{case2}
		f'\left(x^*\right) \geq 0, & \text{ with } x^*=x_l, \\
		\label{case3}
		f'\left(x^*\right) \leq 0, & \text{ with } x^*=x_u,
	\end{align}
where $f'\left(x^*\right) = a + 2c(x^*-b)$. Consider the following three cases. In case 1, i.e., $ x_l < b - \frac{a}{2c} < x_u$, by (\ref{case1}), we have $x^* = b - \frac{a}{2c}$. In case 2, i.e., $ b - \frac{a}{2c} \leq x_l$, by (\ref{case2}), we have $a + 2c(x^*-b)\geq0$, implying $x^*=x_l$. In case 3, i.e.,  $ b - \frac{a}{2c}\geq x_u$, by (\ref{case3}), we have $a + 2c(x^*-b)\geq0$, implying $x^*=x_u$. Thus, we complete the proof. 
\end{IEEEproof}

\begin{Lem}
	\label{Analytical_2}
	The optimal solution of $\underset{x \in [x_l,\infty)}{\operatorname{min}}f\left(x\right)$, where $x_l\in \mathbb{R}$, is given by $x^* = \max \left\{b - \frac{a}{2c} ,x_l \right \}$.
\end{Lem}    
\begin{IEEEproof}
	Similarly, an optimal point $x^*$ satisfies (\ref{case1}) or (\ref{case2}). Thus, the proof of Lemma~\ref{Analytical_2} follows from that of Lemma~\ref{Analytical_1}. 
\end{IEEEproof}

By setting $x = \alpha_n, x_l=0, x_u=1, b=\alpha_n^{(k)},c=\frac{1}{2}\left(g_n \mathbf{s}_n^H \boldsymbol{\Sigma}_{\boldsymbol{\alpha}^{(k)}}^{-1} \mathbf{s}_n \right)^2$, and $a = \partial_{n} f_{\textrm{\rm s}}\left(\boldsymbol{\alpha}^{(k)}\right), \textrm{\rm s} \in \left\{\textrm{ML-K},\textrm{MAP-K}\right\}$, we can show Lemma~\ref{Lem_ML} and Lemma~\ref{Lem_MAP} according to Lemma~\ref{Analytical_1}. By setting $x = \gamma_n, x_l=0, x_u=g_u, a = \partial_{n} f_{\textrm{\rm MAP-UR}}\left(\boldsymbol{\gamma}^{(k)}\right),b=\gamma_n^{(k)},c=\frac{1}{2}\left(\mathbf{s}_n^H \boldsymbol{\Sigma}_{\boldsymbol{\gamma}^{(k)}}^{-1} \mathbf{s}_n \right)^2$, we can show Lemma~\ref{Lem_MAP_U}  according to Lemma~\ref{Analytical_1}. By setting $x = \gamma_n, x_l=0, a = \partial_{n} f_{\textrm{\rm ML-UD}}\left(\boldsymbol{\gamma}^{(k)}\right),b=\gamma_n^{(k)},c=\frac{1}{2}\left(\mathbf{s}_n^H \boldsymbol{\Sigma}_{\boldsymbol{\gamma}^{(k)}}^{-1} \mathbf{s}_n \right)^2$, we can show Lemma~\ref{Lem_ML_U}  according to Lemma~\ref{Analytical_2}.

\section{Proofs of Theorem~\ref{Convergence_ML} and Theorem~\ref{Convergence_MAP}}
\label{appendix:convergence}
Denote $\mathcal{S} \triangleq \left\{\textrm{ML-K},\textrm{MAP-K}\right\}$. For all $\textrm{\rm s} \in \mathcal{S}$, we show that the assumptions in \cite[Theorem 1]{razaviyayn2014parallel} are satisfied. (i) $\widetilde{f}_{\textrm{\rm s},n}(\alpha_n;\boldsymbol{\alpha}')$ is continuously differentiable and uniformly strongly convex w.r.t.  $\alpha_n \in [0,1], n \in \mathcal{N}$, for all $ \boldsymbol{\alpha}'\in [0,1]^N$. (ii) \red{$\widetilde{f}_{\textrm{\rm s},n}'(\alpha_n;\boldsymbol{\alpha}')=\partial_{n} f_{\textrm{\rm s}}\left(\boldsymbol{\alpha}'\right)$}, for all $\boldsymbol{\alpha}' \in [0,1]^N$. (iii) \red{$ \widetilde{f}_{\textrm{\rm s}}'(\alpha_n;\boldsymbol{\alpha}')$} is Lipschitz continuous on $[0,1]^N$ w.r.t. $\boldsymbol{\alpha}' \in [0,1]^N$ for all $\alpha_n\in [0,1], n \in \mathcal{N}$. Therefore, Theorem~\ref{Convergence_ML} and Theorem~\ref{Convergence_MAP} readily follow from \cite[Theorem 1]{razaviyayn2014parallel}.
\section{Computational Complexities in Table \ref{table_known} and Table \ref{table_unknown}}
\label{appendix:computational}
Let $C_1, C_2$, and $C_3$ denote some constant integers.
 
PSCA-based algorithms: For PSCA-ML-K, the dominant terms of the flop counts of \red{Steps 5}, 6, 7, \red{and 8} in Algorithm \ref{Algorithm_PSCA_General}  are $8NL^2$, $(4/3)L^3$, $32NL^2$, and $3N$, respectively. Therefore, the dominant term and order of the per-iteration flop counts of PSCA-ML-K are $40NL^2$ and $\mathcal{O}(NL^2)$, respectively. The PSCA-based algorithms share the same procedure except for Step 7.
Specifically, compared to Step 7 of PSCA-ML-K, Step 7 of PSCA-MAP-K for the general model has $2^N$ extra flops in the dominant term, and Step 7 of PSCA-MAP-K for the independent model has $C_1N$ extra flops.
Compared to Step 7 of PSCA-ML-K, Step 7 of PSCA-ML-UD has $N$ fewer flops.
Compared to Step 7 of PSCA-ML-UD, Step 7 of PSCA-MAP-UR has $C_2N$ extra flops caused by the additional terms in (\ref{gradient_MAP_U}).
Thus, we can obtain the dominant term and order of the per-iteration flop counts of each PSCA-based algorithm.

BCD-based algorithms: For BCD-ML-K, the dominant terms of the flop counts of Steps 6 and 7 are $32NL^2$ and $24NL^2$, respectively. Therefore, the dominant term and order of the per-iteration flop counts of BCD-ML-K are $56NL^2$ and $\mathcal{O}(NL^2)$, respectively. 
Specifically, compared to Step 7 of BCD-ML-K, Steps 7-14 of BCD-MAP-K for the general model has $2^N$ extra flops in the dominant term, and Steps 6-8 of BCD-MAP-K for the independent model has $C_3N$ extra flops. Thus, we can obtain the dominant term and order of the per-iteration flop counts of each BCD-based algorithm.

PG-based algorithms: For PG-ML-K, the dominant terms of the flop counts of Steps 5, 6, 7, and 8 are $8NL^2$, $(4/3)L^3$, $32NL^2$, and $4N$, respectively. Therefore, the dominant term and order of the per-iteration flop counts of PG-ML-K are $40NL^2$ and $\mathcal{O}(NL^2)$, respectively. 
The PG-based algorithms share the same procedure except for Step 8. Specifically, Step 8 of PG-ML-UD has $N$ fewer flops than Step 8 of PG-ML-K. Thus, we can obtain the dominant term and order of the per-iteration flop counts of each PG-based algorithm.
\color{black}

\end{appendices}